%% file: main.tex
\documentclass[10pt,a4paper,reqno]{amsart}

\usepackage{times}
\usepackage{mathrsfs,mathtools,bm}
\usepackage[inline]{enumitem}
\usepackage{stmaryrd}
\usepackage[usenames,dvipsnames]{xcolor}
\usepackage{amsmath}
\usepackage{amsthm}
\usepackage{amssymb}
\usepackage{amsfonts}
\usepackage{subfiles}
\usepackage{dsfont}
\usepackage{hyperref}
\usepackage[T1]{fontenc}
\usepackage[english]{babel}
\usepackage{geometry}
\usepackage{todonotes}
\usepackage{xspace}
\allowdisplaybreaks
\makeatletter \@addtoreset{equation}{section}
\def\theequation{\thesection.\arabic{equation}}
\newtheorem{theorem}{Theorem}[section]
\newtheorem*{assumption*}{\assumptionName}
	\providecommand{\assumptionName}{}
	\makeatletter
	\newenvironment{assumption}[1]
	{
		\renewcommand{\assumptionName}{Assumption ($\mathbb{#1}$)}
	 \begin{assumption*}\protected@edef\@currentlabel{$(\mathbb{#1})$}%
 	}
 	{%
  		\end{assumption*}
 	}
	\makeatother
\newtheorem{corollary}[theorem]{Corollary}
\newtheorem{example}[theorem]{Example}
\newtheorem{lemma}[theorem]{Lemma}
\newtheorem{proposition}[theorem]{Proposition}

\theoremstyle{definition}
\newtheorem{definition}[theorem]{Definition}
\newtheorem{remark}[theorem]{Remark}

\makeatletter
\def\namedlabel#1#2{\begingroup
    #2%
    \def\@currentlabel{#2}%
    \phantomsection\label{#1}\endgroup
}

\numberwithin{equation}{section}

\newcommand{\E}{\mathbb{E}}
\newcommand{\Fil}{\mathbb{G}}
\newcommand{\FilG}{\mathbb{G}}

\newcommand{\N}{\mathbb{N}}
\newcommand{\Pm}{\mathbb{P}}
\newcommand{\Rp}{\mathbb{R}_+}
\newcommand{\transp}{\top}

\newcommand{\G}{\mathcal{G}}
\newcommand{\Hcal}{\mathcal{H}}
\newcommand{\Pred}{\mathcal{P}}

\newcommand{\cadlag}{c\`adl\`ag\xspace}

\newcommand{\domain}
		   {\Omega\times\mathbb{R}_{+}\times\R  d\times\R{d\times m}\times\mathfrak H\longrightarrow\mathbb{R}^d}
\newcommand{\dt}{\ud t}
\newcommand{\ds}{\ud s}
\newcommand{\dx}{\ud x}

\newcommand{\mutorthog}
		   {\protect\mathpalette{\protect\independenT}{\perp}}
		   	\def\independenT#1#2{\mathrel{\rlap{$#1#2$}\mkern2mu{#1#2}}}

\newcommand{\mutilde}{\widetilde{\mu}}

\newcommand{\odsdx}{\pair{\omega;\ds,\dx}}

\newcommand{\ud}{\ensuremath{\mathrm{d}}}

\newcommand{\abs}[1]{\left\vert#1\right\vert}
\newcommand{\borel}[1]{\mathcal{B}\left({#1}\right)}
\newcommand{\compd}[1]{\widetilde{\mu}^{#1}\left(\ud t,\ud x \right)}

\newcommand{\dA}[1]{\ud A_{#1}}
\newcommand{\dC}[1]{\ud C_{#1}}

\newcommand{\e}[1]{\mathrm{e}^{#1}}
\newcommand{\Expect}[1]{\mathbb{E}\left[#1\right]}
\newcommand{\fnorm}[1]{\norm{\frac{f}{\alpha}}^2_{\mathbb{H}^2_{#1}}}

\newcommand{\Htwos}[1]{\mathbb{H}^{2#1}_{\beta}}
\newcommand{\Htwoshat}[1]{\mathbb{H}^{2#1}_{\hat{\beta}}}

\newcommand{\norm}[1]{\left\Vert#1\right\Vert}
\newcommand{\pair}[1]{\ensuremath{\! \left( #1 \right) }}
\newcommand{\pqc}[1]{\!\left\langle#1\right\rangle}

\newcommand{\R}[1]{\mathbb{R}^{#1}}

\newcommand{\set}[1]{\ensuremath{\! \left \{ #1 \right \} }}

\newcommand{\timegraph}[1]{\llbracket #1\rrbracket}
\newcommand{\tnorm}[1]{\left\vert\kern-0.25ex\left\vert\kern-0.25ex\left\vert #1 
    \right\vert\kern-0.25ex\right\vert\kern-0.25ex\right\vert}
\newcommand{\xinorm}[1]{\norm{\xi}^2_{\mathbb{L}^2_{#1}}}

\newcommand{\Htwo}[2]{\mathbb{H}^2_{#1}(#2)}

\newcommand{\stochint}[2]{\llbracket #1,#2\rrbracket}

\newcommand{\CGExpect}[3]{\mathbb{E}\left[\left.#1\right|\G_{#2}^{#3}\right]}

\newcommand\numberthis{\addtocounter{equation}{1}\tag{\theequation}}

\DeclareMathOperator{\tr}{Tr}

\DeclareMathOperator{\Var}{Var}

\DeclareMathAlphabet{\mathpzc}{OT1}{pzc}{m}{it}

\definecolor{red}{rgb}{0.7,0.15,0.15}
\definecolor{green}{rgb}{0,0.5,0}
\definecolor{blue}{rgb}{0,0,0.7}
\hypersetup{	
				colorlinks, linkcolor={red},
				citecolor={green}, urlcolor={blue}
			}
\setlist[enumerate]{itemindent=1cm, leftmargin=0cm}

\newcommand{\oX}{\vphantom{X}\smash{\overline X}}

\begin{document}

\setlength\parindent{0pt}

\title[Existence and uniqueness results for BSDE with jumps: the whole nine yards]{Existence and uniqueness results for BSDE with jumps:\\the whole nine yards}

\author[A. Papapantoleon]{Antonis Papapantoleon}
\author[D. Possama\"{i}]{Dylan Possama\"{i}}
\author[A. Saplaouras]{Alexandros Saplaouras}

\address{Department of Mathematics, National Technical University of Athens, Zografou Campus, 15780 Athens, Greece}
\email{papapan@math.ntua.gr}

\address{Columbia University, IEOR, 500W 120th St., New York, NY 10027, USA}
\email{dp2917@columbia.edu}

\address{Department of Mathematics, University of Michigan, East Hall, 530 Church Street, Ann Arbor, MI 48109-1043, USA}
\email{asaplaou@umich.edu}

\thanks{We thank Martin Schweizer, two anonymous referees and the associate editor for their comments that have resulted in a significant improvement of the manuscript. Alexandros Saplaouras gratefully acknowledges the financial support from the DFG Research Training Group 1845 ``Stochastic Analysis with Applications in Biology, Finance and Physics''. Dylan Possama\"i gratefully acknowledges the financial support of the ANR project {\sc Pacman}, ANR-16-CE05-0027. Moreover, all authors gratefully acknowledge the financial support from the {\sc Procope} project ``Financial markets in transition: mathematical models and challenges''.}



\date{}

\begin{abstract}
This paper is devoted to obtaining a wellposedness result for multidimensional BSDEs with possibly unbounded random time horizon and driven by a general martingale in a filtration only assumed to satisfy the usual hypotheses, \textit{i.e.} the filtration may be stochastically discontinuous. 
We show that for stochastic Lipschitz generators and unbounded, possibly infinite, time horizon, these equations admit a unique solution in appropriately weighted spaces. 
Our result allows in particular to obtain a wellposedness result for BSDEs driven by discrete--time approximations of general martingales.
\end{abstract}

\maketitle \frenchspacing

\section{Introduction}

\subfile{introduction}

\section{Preliminaries}\label{sec:2}

\subfile{preliminaries}

\subsection{Suitable spaces and associated results}\label{subsec:SuitableSpaces}
\subfile{suitable_spaces}
 
\subsection{A useful lemma for generalized inverses}

\subfile{A_useful_lemma_for_generalised_inverses}

\section{Backward stochastic differential equations driven by stochastically discontinuous martingales}\label{sec:BSDE}

\subfile{notation}

\subsection{Formulation of the problem}

\subfile{Formulation_of_the_problem}

\subsection{Existence and uniqueness: statement}

\subfile{Existence_and_uniqueness_statement}

\subsection{Comparison with the related literature}\label{subsec:Compare}
\subfile{Comparison_literature}

\subsection{\textbf{\textit{A priori}} estimates}

\subfile{A_priori_estimates}

\subsection{Proof of the main theorem}

\subfile{Proof_of_the_main_theorem}\label{subsec:BSDEMainTheoremProof}

\subsection{An alternative approach in the Lipschitz setting}\label{subsec:AlternEstim}
\subfile{New_estimates}
\section{Applications}\label{sec:4}
\subfile{Applications}

\appendix
\section{Proofs of results of Subsection \ref{subsec:OrthDec}}\label{AppendixPreliminaries}
\subfile{Auxilliary_Results_Proof_of_Propositions_Preliminaries}

\section{Proof of Lemma \ref{mainlemma}}\label{AppendixA}
\subfile{Auxilliary_Results_Proof_of_Lemma_of_generalised_inverses.tex}

\section{Proof of Lemma \ref{lem:const}}\label{AppendixB}
\subfile{Auxilliary_Results_Proof_of_Lemma_of_constants.tex}

\section{Auxiliary results on optional measures}

\subfile{Auxilliary_Results_Optional_Measures.tex}

\section{Auxiliary analysis of Subsection \ref{subsec:AlternEstim}}\label{App:AuxAnalysis}
\subfile{Auxilliary_Results_New_estimates.tex}

\bibliographystyle{plain}
\small
\bibliography{bibliographyDylan}

\end{document}

%% file: introduction.tex
A generally acknowledged fact is that backward stochastic differential equations (BSDEs for short) were introduced in their linear version by Bismut 
\cite{bismut1973analyse,bismut1973conjugate} in 1973, as an adjoint equation in the Pontryagin stochastic maximum principle. 
However, around the same time, and most probably a bit before\footnote{The authors are indebted to Sa\"id Hamad\`ene for pointing out this reference. The published version of \cite{davis1973dynamic} states that the article was received on October 27, 1971. It is also present in the bibliography of \cite{bismut1973conjugate}, though it is never referred to in the text.}, Davis and Varaiya \cite{davis1973dynamic} (see in particular their Theorem 5.1) also studied what can be considered as a prototype of a linear BSDE for characterizing the value function and the optimal controls of stochastic control problems with drift control only. 
Such linear BSDEs, still in the context of the stochastic maximum principle, were also used by Arkin and Saksonov \cite{arkin1979necessary}, Bensoussan 
\cite{bensoussan1983maximum} and Kabanov \cite{kabanov1978on}. 
The first non--linear versions of these objects were once again introduced, under the form of a Riccati equation, by Bismut \cite{bismut1978controle} and a few years later by Chitashvili \cite{chitashvili1983martingale} and Chitashvili and Mania \cite{chitashvili1987optimal,chitashvili1987optimal2}. 
Nonetheless, the first study presenting a systematic treatment of non--linear BSDEs is the seminal paper by Pardoux and Peng \cite{pardoux1990adapted}. 
Since then, and especially following the illuminating survey article of El Karoui, Peng and Quenez \cite{el1997backward}, BSDEs have become a particularly active field of research, due to their numerous potential applications to mathematical finance, partial differential equations, game theory, economics, and more generally in stochastic calculus and analysis\footnote{We emphasize that the references given below are just the tip of the iceberg, though most of them are, in our view, among the major ones of the field. Nonetheless, we do not make any claim about comprehensiveness of the following list.}.

\vspace{0.5em}
Let $T>0$ be fixed and consider a fixed filtered probability space $(\Omega, \mathcal G,\mathbb G:=(\mathcal G_t)_{0\leq t\leq T},\mathbb P)$ where $\mathbb G$ is a Brownian filtration generated by some $d-$dimensional Brownian motion $W$. 
Solving a BSDE with terminal condition $\xi$ (which is an $\mathbb R-$valued and $\mathcal G_T-$measurable random variable) and $\mathbb G-$adapted generator $f:\Omega\times[0,T]\times\mathbb R\times\mathbb R^d\longrightarrow \mathbb R$, amounts to finding a pair of processes $(Y,Z)$ which are respectively $\mathbb G-$progressively measurable and $\mathbb G-$predictable such that
$$
Y_t=\xi+\int_t^Tf(s,Y_s,Z_s)\ds-\int_t^TZ_s^\top \ud W_s, \ \ t\in[0,T],
$$
holds, $\mathbb P-a.s.$
After the work \cite{pardoux1990adapted} obtained existence and uniqueness of the solution of the above BSDE in $\mathbb L^2-$type spaces under square integrability assumptions on $\xi$ and $f(s,0,0)$, and uniform Lipschitz continuity of $f$ in $(y,z)$, generalizations of the theory have followed several different paths. 
The first one mainly aimed at weakening the Lipschitz assumption on $f$, and still considered Brownian filtrations. 
Hence, Mao \cite{mao1995adapted} considered uniformly continuous generators, Hamad\`ene \cite{hamadene1996equations} extended the result to the locally Lipschitz case, Lepeltier and San Mart\'in \cite{lepeltier1997backward} to the continuous and linear growth case in $(y,z)$, Briand and Carmona \cite{briand2000bsdes} to the case of a continuous generator Lipschitz in $z$ with polynomial growth in $y$, and Pardoux \cite{pardoux1999bsdes} to the case of a generator monotonic with arbitrary growth in $y$ and Lipschitz in $z$. 
Some authors also obtained wellposedness results in $\mathbb L^p-$type spaces, among which we mention \cite{el1997backward} for $p\geq2$, Briand, Delyon, Hu, Pardoux and Stoica \cite{briand2003p} and Briand and Hu \cite{briand2006bsde} for $p\geq 1$ (see also the papers of Fan \cite{fan2015bounded,fan20171} and Hu and Tang \cite{hu2017existence} for recent results and other references). 
Some attention has also been given to the so--called stochastic Lipschitz case, where the generator is Lipschitz continuous in $(y,z)$ but with constants which are actually random processes themselves. 
There are few papers going in this direction, among which we can mention El Karoui and Huang \cite{elkaroui1997general}, Bender and Kohlmann \cite{bender2000bsdes}, Wang, Ran and Chen \cite{wang2007lp} as well as Briand and Confortola \cite{briand2008bsdes}.  

\vspace{0.5em}
The first results going beyond the linear growth assumption in $z$, which assumed quadratic growth, were obtained independently by Kobylanski 
\cite{kobylanski1997resultats,magdalena1998quelques,kobylanski2000backward} and Dermoune, Hamad\`ene and Ouknine \cite{dermoune1999backward}, for 
bounded $\xi$ and $f$ Lipschitz in $y$. 
These results were then further studied by Eddhabi and Ouknine \cite{eddahbi2002limit}, and improved by Lepeltier and San Mart\'in \cite{lepeltier1998existence,lepeltier2002existence}, Briand, Lepeltier and San Mart\'in \cite{briand2007one} and revisited by Briand and \'Elie \cite{briand2013simple}, but still for bounded $\xi$. 
Wellposedness in the quadratic case when $\xi$ has sufficiently large exponential moments was then investigated by Briand and Hu 
\cite{briand2006bsde,briand2008quadratic}, followed by Delbaen, Hu and Richou \cite{delbaen2011uniqueness,delbaen2015uniqueness}, Essaky and Hassani \cite{essaky2011general}, and Briand and Richou \cite{briand2017uniqueness}. 
A specific quadratic setting with only square integrable terminal conditions has been considered recently by Bahlali, Eddahbi and Ouknine \cite{bahlali2013solvability,bahlali2017quadratic}, while a result with logarithmic growth was also obtained by Bahlali and El Asri \cite{bahlali2012stochastic}, and Bahlali, Kebiri, Khelfallah and Moussaoui \cite{bahlali2017one}. 
The case of a generator with super--quadratic growth in $z$ was proved to be essentially ill--posed by Delbaen, Hu and Bao \cite{delbaen2011backward} in a general non--Markovian framework, before Richou \cite{richou2012markovian}, Cheridito and Stadje \cite{cheridito2012existence} and Masiero and Richou \cite{masiero2013note} proved that wellposedness could be recovered in a Markovian setting, when $f$ has polynomial growth in $(y,z)$. 
Let us also mention the contributions by Cheridito and Nam \cite{cheridito2014bsdes}, when $\xi$ has a bounded Malliavin derivative, by Drapeau, Heyne and Kupper \cite{drapeau2013minimal} who considered minimal super--solutions of BSDEs when the generator is monotone in $y$ and convex in $z$, and by Heyne, Kupper and Mainberger \cite{heyne2014minimal} who considered lower semicontinuous generators.

\vspace{0.5em}
Most of the papers mentioned above treated the so--called one--dimensional BSDEs, that is for which the process $Y$ is $\mathbb R-$valued, but extensions to 
multidimensional settings were also explored. Hence in Lipschitz or locally Lipschitz settings with monotonicity assumptions, we can mention the works of Pardoux \cite{pardoux1999bsdes}, Bahlali \cite{bahlali2001backward,bahlali2002existence}, Bahlali, Essaky, Hassani and Pardoux \cite{bahlali2002existence2}, and Bahlali, Essaky and Hassani \cite{bahlali2010multidimensional,bahlali2015existence}. 
An early result in the case of a continuous generator in a Markovian setting was also treated by Hamad\`ene, Lepeltier and Peng \cite{hamadene1997bsdes}. 
The quadratic multidimensional case is much more involved. 
Tevzadze \cite{tevzadze2008solvability} was the first to obtain a wellposedness result in the case of a bounded and sufficiently small terminal condition. 
It was then proved by Frei and dos Reis \cite{frei2011financial} and Frei \cite{frei2014splitting} (see also Espinosa and Touzi \cite{espinosa2015optimal} for a related problem) that even in seemingly benign situations, existence of global solutions could fail. Later on, Cheridito and Nam \cite{cheridito2015multidimensional}, Kardaras, Xing and {\v{Z}}itkovi\'c \cite{kardaras2015incomplete}, Kramkov and Pulido \cite{kramkov2014stability,kramkov2016system}, Hu and Tang \cite{hu2015multi}, Jamneshan, Kupper and Luo \cite{jamneshan2014multidimensional}, or more recently Kupper, Luo and Tangpi \cite{luo2015solvability} and \'Elie and Possama\"i \cite{elie2016contracting}, all obtained some results, but only in particular instances. 
A breakthrough was then obtained by Xing and \v{Z}itkovi\'c \cite{xing2016class}, who obtained quite general existence and uniqueness results, but in a Markovian framework, while Harter and Richou \cite{harter2016stability} and Jamneshan, Kupper and Luo \cite{jamneshan2016solvability} have proved positive results in the general setting.

\vspace{0.5em}
A second possible generalization of these results consisted in extending them to the case where $T$ is assumed to be a, possibly unbounded, stopping time. 
Hence, Peng \cite{peng1991probabilistic}, Darling and Pardoux \cite{darling1997backwards}, Briand and Hu \cite{briand1998stability}, Bahlali, Elouaflin and N'zi \cite{bahlali2004backward}, Royer \cite{royer2004bsdes}, Hu and Tessitore \cite{hu2007bsde} and Briand and Confortola \cite{briand2008quadratic2} all studied this problem, applying it to homogenization or representation problems for elliptic PDEs and stochastic control in infinite horizon. 
This theory was recently revisited by Lin, Ren, Touzi and Yang \cite{lin2018second} in the context of second--order BSDEs with random horizon.

\vspace{0.5em}
Another avenue of generalization concerned the underlying filtration itself, which could be assumed to no longer be Brownian, as well as the driving 
martingale, which could also be more general than a Brownian motion. 
In such cases, the predictable (martingale) representation property may fail to hold, and one has in general to add another martingale to the definition of a 
solution. 
Hence, for a given martingale $M$, the problem becomes to find a triplet of processes $(Y,Z,N)$ such that $N$ is orthogonal to $M$ and
$$
Y_t=\xi+\int_t^Tf(s,Y_s,Z_s)\ud C_s-\int_t^TZ_s^\top \ud M_s-\int_t^T\ud N_s,\ 
t\in[0,T],\ \mathbb P-a.s.,
$$
where the non-decreasing process $C$ is absolutely continuous with respect to 
$\langle M\rangle$.

\vspace{0.5em}
As far as we know, the first paper where such BSDEs appeared is due to Chitashvili \cite{chitashvili1983martingale} (see in particular the corollary at the end of page 91). 
Then, results on BSDEs driven by a general c\`adl\`ag martingale were obtained by Buckdahn \cite{buckdahn1993backward}, El Karoui {\it et al.} \cite{el1997backward}, as well as El Karoui and Huang \cite{elkaroui1997general}, Briand, Delyon and M\'emin \cite{briand2002robustness} and Carbone, Ferrario and Santacroce \cite{carbone2008backward}, in Lipschitz type settings. 
The case of generators with quadratic growth has also been investigated by Tevzadze \cite{tevzadze2008solvability}, Morlais \cite{morlais2009quadratic}, R\'eveillac \cite{reveillac2012orthogonal}, Imkeller, R\'eveillac and Richter \cite{imkeller2012differentiability}, Mocha and Westray \cite{mocha2012quadratic} and Barrieu and El Karoui \cite{barrieu2013monotone}. 
More general versions of these equations, coined semimartingale BSDEs, were also studied in depth in the context of financial applications, especially utility maximization, see Bordigoni, Matoussi and Schweizer \cite{bordigoni2007stochastic}, as well as Hu and Schweizer \cite{hu2011some}, and mean--variance hedging, see Bobrovnytska and Schweizer \cite{bobrovnytska2004mean}, Mania and Tevzadze \cite{mania2000semimartingale,mania2003semimartingale2,mania2003unified}, Mania, 
Santacroce and Tevzadze \cite{mania2002semimartingale,mania2003semimartingale}, Mania and Schweizer \cite{mania2005dynamic} as well as Jeanblanc, Mania, Santacroce and Schweizer \cite{jeanblanc2012mean}. 

\vspace{0.5em}
When one has more information on the filtration, it may be possible to specify the orthogonal martingale $N$ in the definition of the solution. 
For instance, if the filtration is generated by a Brownian motion and an orthogonal Poisson random measure, one ends up with the so--called BSDEs with jumps, which were introduced first by Tang and Li \cite{tang1994necessary}, followed notably by Buckdahn and Pardoux \cite{buckdahn1994bsde}, Barles, Buckdahn and Pardoux \cite{barles1997backward}, Situ \cite{rong1997solutions}, Royer \cite{royer2006backward}, Becherer \cite{becherer2006bounded}, Morlais \cite{morlais2009utility,morlais2010new}, Ankirchner, Blanchet-Scalliet and Eyraud-Loisel \cite{ankirchner2010credit}, Lim and Quenez \cite{lim2011exponential}, Jenablanc, Matoussi and Ngoupeyou \cite{jeanblanc2012robust}, Kharroubi, Quenez and Sulem \cite{quenez2013bsdes}, Lim and Ngoupeyou \cite{kharroubi2013mean}, Kharroubi and Lim \cite{kharroubi2014progressive}, Laeven and Stadje \cite{laeven2014robust}, Richter \cite{richter2014explicit}, Jeanblanc, Mastrolia, Possama\"i and R\'eveillac \cite{jeanblanc2015utility}, Kazi-Tani, Possama\"i and Zhou \cite{kazi2015quadratic1,kazi2015quadratic}, Fujii and Takahashi \cite{fujii2015quadratic}, Dumitrescu, Quenez and Sulem \cite{dumitrescu2016bsdes}, and El Karoui, Matoussi and Ngoupeyou 
\cite{karoui2016quadratic}, while the specific case of L\'evy processes was treated by Nualart and Schoutens \cite{nualart2001backward} and later Bahlali, Eddahbi and Essaky \cite{bahlali2003bsde}. 
A general presentation has been proposed recently by Kruse and Popier \cite{kruse2015bsdes,kruse2017lp}, to which we refer for more references (see also the recent paper of Yao \cite{yao2017lp}).

\vspace{0.5em}
One point that is actually shared by all the above references, is that the underlying filtration is assumed to be quasi-left continuous, which for instance rules out the possibility that any of the involved processes has jumps at predictable, and {\it a fortiori} deterministic times. 
The important simplification that arises is that the process $C$ is then necessarily continuous in time. 
As far as we know, the first articles that went beyond this assumption were developed in a very nice series of papers by Cohen and Elliott \cite{cohen2012existence} and Cohen, Elliott and Pearce \cite{cohen2010general}, where the only assumption on the filtration is that the associated $\mathbb L^2$ space is separable, so that a very general martingale representation result due to Davis and Varaiya \cite{davis1974multiplicity}, involving countably many orthogonal martingales, holds. 
In these works, the martingales driving the BSDE are actually imposed by the filtration, and not chosen {\it a priori}, and the non--decreasing process $C$ is not necessarily related to them, but has to be deterministic and can have jumps in general, though they have to be small for existence to hold (see \cite[Theorem 5.1]{cohen2012existence}). 
A similar approach is taken by Hassani and Ouknine in \cite{hassani2002general}, where a form of BSDE is considered using generic maps from a space of semimartingales to the spaces of square--integrable martingales and of finite--variation processes integrable with respect to a given continuous increasing process. 
Similarly, Bandini \cite{bandini2015existence} obtained wellposedness results in the context of a general filtration allowing for jumps, with a fixed driving martingale and associated random process $C$, which must have again small jumps, see \cite[Equation (1.1)]{bandini2015existence}. 
Let us also mention the work by Confortola, Fuhrman and Jacod \cite{confortola2014backward} which concentrates on the pure--jump general case and gives in particular counterexamples to existence. 
Finally, Bouchard, Possama\"i, Tan and Zhou \cite{bouchard2015unified} provided a general method to obtain {\it a priori} estimates in general filtrations when the martingale driving the equation has quadratic variation absolutely continuous with respect to the Lebesgue measure. 

\vspace{0.5em}
In this paper, we improve the general result on existence and uniqueness of solutions of backward stochastic differential equations given by El Karoui and Huang in \cite{elkaroui1997general} to the case where the martingale $M$ driving the equation is possibly stochastically discontinuous. 
In other words, our framework includes as driving martingales discrete--time approximations of general martingales as well as $K-$almost quasi--left--continuous martingales, {\it i.e.} processes whose compensator has jumps which are almost surely bounded by some constant $K$.
Unlike all the related papers mentioned above (with the notable exception of \cite{cohen2012existence}, see their Theorem 6.1, albeit with a deterministic Lipschitz constant), this bound $K$ can actually be arbitrarily large.
However, the product of this bound and the maximum (of functionals) of the Lipschitz constants needs to be small, which is in line with the previous literature.
Otherwise, we remain in the same relaxed framework regarding the generator, that is to say we assume that it satisfies a stochastic Lipschitz property, and do not assume that the martingale possesses the predictable representation property. 
Furthermore, we work in a setting with random horizon. 
This result enables us to treat under the same framework continuous--time as well as discrete--time BSDEs. 
The method of proof is somehow similar to the one given in \cite{elkaroui1997general}, but the required estimates are much harder to prove in our setting due to the possible jumps of the non--decreasing process $C$. 
We also emphasize that this wellposedness result will be of fundamental importance in a related forthcoming work, where we will use it to study robustness properties of general BSDEs, extending well--known results on stability of semimartingale decompositions with respect to the extended convergence.

\vspace{0.5em}
This paper is structured as follows: in Section \ref{sec:2} we introduce the notation and several results that will be useful in the analysis. In Section 
\ref{sec:BSDE} we prove \textit{a priori} estimates for the considered class of BSDEs and provide the existence and uniqueness results. 
Finally, Section \ref{sec:4} discusses some applications of the main results, while the Appendices contain proofs and auxiliary results.

\subsection*{Notation} 

Let $\mathbb R_+$ denote the set of non-negative real numbers. 
For any positive integer $\ell$, and for any $(x,y)\in \R \ell\times\R \ell,$ $\abs x$ will denote the Euclidean norm of $x$. 
For any additional integer $q$, a $q\times \ell-$matrix with real entries will be considered as an element of $\R {q\times \ell}$. 
For any $z\in\R {q\times \ell}$, its transpose will be denoted by $z^\top\in\R {\ell\times q}$. 
We endow $\R {q\times \ell}$ with the norm defined for any $z\in\R {q\times \ell}$ by $\norm z^2:=\tr [z^\transp z]$ and remind the reader that this norm derives from the inner product defined for any $(z,u)\in\R {q\times \ell}\times\R {q\times \ell}$ by $\tr[zu^{\transp}]$. 
We abuse notation and denote by $0$ the neutral element in the group $(\mathbb R^{q\times\ell},+)$. 
Furthermore, for any finite dimensional topological space $E$, $\mathcal B(E)$ will denote the associated Borel $\sigma-$algebra. 
In addition, for any other finite dimensional space $F$, and for any non-negative measure $\nu$ on $(\mathbb R_+\times E,\mathcal B(\mathbb R_+)\otimes\mathcal B(E))$, we will denote indifferently Lebesgue--Stieltjes integrals of any measurable map $f:(\mathbb R_+\times E,\mathcal B(\mathbb R_+)\otimes\mathcal B(E))\longrightarrow (F,\mathcal B(F))$, by
 $$\int_{(u,t]\times A}f(s,x)\nu(\ud s,\ud x), \text{ for any }  (t,A)\in\mathbb R_+\times\mathcal B(E),$$
 $$\int_{(u,\infty)\times A}f(s,x)\nu(\ud s,\ud x),  \text{ for any }   A\in\mathcal B(E),$$
where the integrals are to be understood in a component--wise sense. 
Finally, the letters $p,q, d,m$ and $n$ are reserved to denote arbitrary positive integers.
The reader may already keep in mind that $m$ will denote the dimension of the state space of an It\=o integrator, $n$ will denote the dimension of the state space of a process associated to an integer--valued random measure and $d$ will denote the dimension of the state space of a stochastic integral.

%% file: preliminaries.tex
\subsection{The stochastic basis}\label{subsec:StochBasis}

Let $(\Omega,\mathcal G,\mathbb G,\mathbb P)$ be a complete stochastic basis in the sense of Jacod and Shiryaev \cite[Definition I.1.3]{jacod2003limit}.
Expectations under $\mathbb P$ will be denoted by $\mathbb E[\cdot]$. 
We will then denote the set of $\mathbb R^p-$valued, square--integrable $\mathbb G-$martingales by $\Hcal^2(\mathbb R^p)$, \emph{i.e.}
\begin{align*}
\mathcal H^2(\mathbb R^p):=\big\{ X:\Omega\times\mathbb R_+\to \mathbb R^p, X \text{ is a $\mathbb G-$martingale with }\sup_{t\in\mathbb R_+}\mathbb E[\vert X_t\vert^2]<\infty \big\}. 
\end{align*}
Let $X\in\Hcal^2(\mathbb R^p)$, then its norm will be defined by $$\Vert X\Vert^2_{\mathcal H^2(\mathbb R^p)}:= \mathbb E[ \vert X_\infty \vert^2]= \Expect{{\rm Tr}[\langle X\rangle_\infty]}.$$

\vspace{0.5em}
In the sequel we will say that $M,N\in\Hcal^2(\mathbb R)$ are \emph{$($mutually$)$ 
orthogonal}, denoted by $M \mutorthog N$, if $MN$ is a martingale, see 
\cite[Definition I.4.11.a, Lemma I.4.13.c, Proposition I.4.15]{jacod2003limit} for equivalent definitions.

\vspace{0.5em}
For a subset $\mathcal N$ of $\mathcal H^2(\mathbb R^p)$, we denote the space of martingales orthogonal to each component of every element of $\mathcal N$ by $\mathcal N^\perp$, \emph{i.e.}
\begin{align*}
\mathcal N^\perp 			
  &:=\{M\in\mathcal{H}^2(\mathbb R^p),\   \langle M,N\rangle=0\text{ for every } N\in\mathcal N\},
\end{align*} 
where we suppress the explicit dependence in the state space in the notation. 
Observe, however, that in the above definition the predictable quadratic covariation is an $\mathbb R^{p\times p}-$valued process.
A martingale $M\in\Hcal^2(\mathbb R^p)$ will be called a \emph{purely discontinuous} martingale if $M_0=0$ and if each of its components is orthogonal to all continuous martingales of $\Hcal^2(\mathbb R).$ 
Using \cite[Corollary I.4.16]{jacod2003limit} we can decompose $\Hcal^2(\mathbb R^p)$ as follows
\begin{align}\label{SqIntMartDecomp}
\mathcal H^{2}(\mathbb R^p)=\mathcal H^{2,c}(\mathbb R^p)\oplus\mathcal H^{2,d}(\mathbb R^p),
\end{align}
where $\Hcal^{2,c}(\mathbb R^p)$ is the subspace of $\Hcal^2(\mathbb R^p)$ consisting of continuous square--integrable martingales and $\Hcal^{2,d}(\mathbb R^p)$ is the subspace of $\Hcal^2$ consisting of all purely discontinuous square--integrable martingales. 
It follows then from \cite[Theorem I.4.18]{jacod2003limit}, that any $\mathbb G-$martingale $X\in\Hcal^{2}(\mathbb R^p)$ admits a unique (up to $\mathbb P-$indistinguishability) decomposition
\begin{align*}
X_\cdot=X_0+X^c_\cdot+X^d_\cdot,
\end{align*}
where $X^c_0=X^d_0=0$. 
The process $X^c\in\Hcal^{2,c}(\mathbb R^p)$ will be called the \emph{continuous martingale part of $X$}, while the process $X^d\in\Hcal^{2,d}(\mathbb R^p)$ will be called the \emph{purely discontinuous martingale part of $X$}. The pair $(X^c, X^d)$ will be called the \emph{natural pair of $X$ $($under $\mathbb G )$.}

\subsection{Stochastic integrals}
Let $X\in\mathcal H^2(\mathbb R^m)$ and $C$ be a predictable, non--decreasing and \cadlag process such that 
\begin{align}\label{FactorizationPQC}
\langle X\rangle = \int_{(0,\cdot]} \frac{\ud \langle X\rangle_s}{\ud C_s} \ud C_s,
\end{align} 
where the equality is understood componentwise.
That is to say, $\frac{\ud \langle X\rangle}{\ud C}$ is a predictable process with values in the set of  all symmetric, non--negative definite $m\times m$ matrices. 
In the next lines, we follow closely \cite[Section III.6.a]{jacod2003limit}. 
We start by defining
\begin{align*}
\mathbb H^2(X) &:=\left\{ Z:(\Omega\times\Rp,\Pred) 
\longrightarrow (\R{d\times m},\mathcal B(\R{d\times m})),\  
\Expect{\int_0^\infty {\rm Tr}\left[Z_t \frac{\ud \langle X\rangle_s}{\ud C_s}Z_t^\transp\right]\ud C_t}<\infty \right\},
\end{align*}
where $\mathcal{P}$ denotes the $\mathbb G-$predictable $\sigma-$field on $\Omega\times\mathbb R_+$; see \cite[Definition I.2.1]{jacod2003limit}.
Let $Z\in\mathbb H^2(X)$, then the It\=o stochastic integral of $Z$ with respect to $X$ is well defined and is an element of $\Hcal^{2}(\mathbb R^d)$, see \cite[Theorem III.6.4]{jacod2003limit}. 
It will be denoted by $\int_0^{\cdot}Z_s\ \ud X_s$ or $Z\cdot X$ interchangeably, and we will also use the same notation for any Stieltjes--type integral. 
Moreover, by \cite[Theorem III.6.4.c)]{jacod2003limit} we have that $(Z \frac{\ud \langle X\rangle}{\ud C} Z^\transp)\cdot C = \langle Z\cdot X\rangle$, hence the following equality holds
\begin{align*}
\norm{Z}_{\mathbb H^2(X)}^2:=\Expect{\int_0^\infty {\rm Tr}\left[Z_t \frac{\ud \langle X\rangle_s}{\ud C_s}Z_t^\transp\right]\ud C_t}=\Expect{\tr[\,\pqc{Z\cdot X}_\infty]}.
\end{align*}
We will denote the space of It\=o stochastic integrals of processes in 
$\mathbb H^2(X)$ with respect to $X$ by $\mathcal{L}^2(X)$.
In particular, for $X^c\in\mathcal H^{2,c}(\mathbb R^m)$ we remind the reader that, by \cite[Theorem III.4.5]{jacod2003limit}, $Z\cdot X^c\in\mathcal H^{2,c}(\mathbb R^d)$ for every $Z\in\mathbb H^2(X^c)$, 
\emph{i.e.} $\mathcal{L}^2(X^c)\subset\mathcal H^{2,c}(\mathbb R^d)$.

\vspace{0.5em}
Let us define the space
\[
\big(\widetilde{\Omega},\widetilde{\Pred}\big):=
\big(\Omega\times\Rp\times\R {n},\Pred\otimes \borel{\R {n}}\big).
\] 
A measurable function 
$U:\big(\widetilde{\Omega},\widetilde{\Pred}\big)\longrightarrow 
\left(\R {d},\borel{\R {d}} \right)$ is called 
\emph{$\widetilde{\Pred}-$measurable function} or 
\emph{$\mathbb G-$predictable function}.

\vspace{0.5em}
Let $\mu := \set{\mu\pair{\omega;\dt,\dx}}_{\omega\in\Omega}$ be a random 
measure on $\Rp\times\R{n}$, that is to say a family of non--negative measures 
defined on $\pair{\Rp\times\R{n},\borel{\Rp}\otimes\borel{\R{n}}}$ satisfying 
$\mu\pair{\omega;\set{0}\times\R{n}}=0$, identically. For a $\mathbb 
G-$predictable function $U$, we define the process
$$
U\star\mu_\cdot(\omega) :=
			\begin{cases}
			\displaystyle \int_{(0,\cdot]\times\R{n}} 
U\pair{\omega,s,x} \mu\odsdx,\textrm{ if } \int_{(0,\cdot]\times\R{n}} 
\abs{U\pair{\omega,s,x}} \mu\odsdx<\infty,\\
		\displaystyle	\infty,\textrm{ otherwise}.
			\end{cases}
$$
Let us now consider some $X\in\mathcal H^{2,d}(\mathbb R^n)$. 
We associate to $X$ the $\mathbb G-$optional integer--valued random measure $\mu^{X}$ on $\Rp\times\R{m}$ defined by
\[
\mu^X\pair{\omega;\dt,\dx} 
  := \sum_{s>0} \mathds{1}_{\set{\Delta X_s(\omega) \neq 0}} 
     \delta_{(s,\Delta X_s(\omega))}\pair{\dt,\dx},
\]
see \cite[Proposition II.1.16]{jacod2003limit}, where, for any $z\in\Rp\times\R{n} $, $\delta_z$ denotes the Dirac measure at the point $z$. 
Notice that $\mu^X(\omega;\mathbb R_+\times\{0\})=0.$
Moreover, for a $\mathbb G-$predictable stopping time $\sigma$ we define the random variable 
\begin{align*} 
\int_{\R{n}}U(\omega,\sigma,x)\mu^{X}(\omega;\set{\sigma}\times\dx) :=&\ 
 U(\omega,\sigma(\omega),\Delta X_{\sigma(\omega)}(\omega))\mathds{1}_{\{\Delta 
X_\sigma\neq 0, |U(\omega,\sigma(\omega),\Delta 
X_{\sigma(\omega)}(\omega))|<\infty\}} \\ 
&+			\infty{\mathds 1}_{\{|U(\omega,\sigma(\omega),\Delta 
X_{\sigma(\omega)}(\omega)|=\infty\}}.
\end{align*}

Since $X\in\mathcal H^2(\mathbb R^n)$, the compensator of $\mu^X$ under $\Pm$ exists, see 
\cite[Theorem II.1.8]{jacod2003limit}. This is the unique, up to a $\Pm-$null set, $\mathbb 
G-$predictable random measure $\nu^X$ on $\Rp\times\R{n}$, for which
\begin{align*}
\Expect{U\star\mu^X_\infty}=\Expect{U\star\nu^X_\infty},
\end{align*}
holds for every non--negative $\mathbb G-$predictable function $U.$

\vspace{0.5em}
For a non-negative $\mathbb G-$predictable function $U$ and a $\mathbb G-$predictable time $\sigma$, whose graph is denoted by $\timegraph{\sigma}$ (see \cite[Notation I.1.22]{jacod2003limit} and the comments afterwards), we define the random variable 
\[
\int_{\R{n}}U\pair{\omega,\sigma,x}\nu^X 
\pair{\omega;\set{\sigma}\times\dx}:=
	\int_{\Rp\times\R{n}} 
U\pair{\omega,\sigma(\omega),x}\mathds{1}_{\timegraph{\sigma}}\ \nu^X\odsdx,
\]
if  $\int_{\Rp\times\R{n}} 
\abs{U\pair{\omega,\sigma(\omega),x}}\mathds{1}_{\timegraph{\sigma}}\ 
\nu^X\odsdx<\infty$, otherwise we define it to be equal to $\infty$. By 
\cite[Property II.1.11]{jacod2003limit}, we have
\begin{align}\label{PropertyII-1-11}
\int_{\R{n}}U\pair{\omega,\sigma,x}\nu^{X}\pair{\omega;\set{\sigma}\times\dx} 
= \CGExpect{ 
\int_{\R{n}}U\pair{\omega,\sigma,x}\mu^{X}\!\pair{\omega;\set{\sigma}\times\dx}}
{{\sigma-}}{}.
\end{align}
In order to simplify notations further, let us denote for any $\mathbb G-$predictable time 
$\sigma$
\begin{align}
\widehat{U}_{\sigma}^{X}&:=\int_{\R{n}}U\pair{\omega,\sigma,x}\nu^{X}\pair{\omega;\set{\sigma}\times\dx}.
\label{notation:U-hat}
\shortintertext{In particular, for $U=\mathds{1}_{\mathbb R^n}$ we define}
\zeta_\sigma^X&:=\int_{\mathbb R^n}\nu^X(\omega;\{\sigma\}\times\ud x)
\label{notation:1-hat}
\end{align}

In order to define the stochastic integral of a $\mathbb G-$predictable function $U$ with respect to the \emph{compensated integer--valued random measure 
\mbox{$\mutilde^X\!:=\mu^X-\nu^X$}}, we will need to consider the following 
class
\begin{align*}
G_2(\widetilde\mu^X)=\set{
	U:\big(\widetilde{\Omega},\widetilde{\Pred}\big)\longrightarrow 
\big(\R{d},\mathcal B(\R{d})\big),\, 
	\Expect{\sum_{t>0} \pair{U\pair{t,\Delta X_t}\mathds{1}_{\set{\Delta 
X_t\neq0}}-\widehat{U}_t^{X}}^2}
	<\infty
}.
\end{align*}
Any element of $G_2(\widetilde\mu^X)$ can be associated to an element of $\Hcal^{2,d}$, uniquely up to $\mathbb P-$indistinguishability via 
\begin{align*}
G_2\pair{\widetilde\mu^X}\ni U\longmapsto U\star\mutilde^{X}\in \Hcal^{2,d},
\end{align*}
see \cite[Definition II.1.27, Proposition II.1.33.a]{jacod2003limit} and \cite[Theorem XI.11.21]{he1992semimartingale}. We call $U\star\mutilde^{X}$ the \emph{stochastic integral of $U$ with respect to $\mutilde^{X}$}. 
We will also make use of the following notation for the space of stochastic 
integrals with respect to $\mutilde^{X}$ which are square integrable martingales
\begin{align*}
\mathcal{K}^2(\widetilde\mu^X)
  := \set{U\star\mutilde^{X},\  U\in G_2(\widetilde\mu^X)}. 
\end{align*}
Moreover, by \cite[Theorem II.1.33]{jacod2003limit} or \cite[Theorem 11.21]{he1992semimartingale}, we have
\begin{align*}
\Expect{\,\pqc{U\star\mutilde^{X}}_{\infty}}<\infty ,
   \textrm{ if and only if } 
U\in G_2\pair{\widetilde\mu^X},
\end{align*}
which enables us to define the following more convenient space
\begin{align*}
\Htwo{}{X}
 &:=\left\{ U:\big(\widetilde{\Omega},\widetilde{\Pred}\big) 
  \longrightarrow \big(\R{d},\mathcal B(\R{d})\big),\ 
  \Expect{\int_0^{\infty} \ud {\rm Tr} 
  \left[\,\pqc{U\star\mutilde^{X}}_t\right]}<\infty \right\},
\end{align*}
and we emphasize that we have the direct identification
$$\mathbb H^2(X)=G_2(\widetilde\mu^X).$$

\subsubsection{Orthogonal decompositions}\label{subsec:OrthDec}

We close this subsection with a discussion on orthogonal decompositions of square integrable martingales.
We associate the measure $M_{\mu}:(\widetilde{\Omega}, \G\otimes \borel{\Rp}\otimes\borel{\R{n}})\longrightarrow \Rp$ to a random measure $\mu$, which is defined as $M_{\mu}(B)=\E[ \mathds{1}_B\star\mu_{\infty}]$. 
We will refer to $M_{\mu}$ as the \emph{Dol\'eans measure associated to $\mu.$}
If there exists a $\mathbb G-$predictable partition $(A_k)_{k\in\N}$ of $\widetilde{\Omega}$ such that $M_{\mu}(A_k)<\infty,$ for every $k\in\N,$ then we will say that $\mu$ is $\mathbb G-$predictably $\sigma-$integrable and we will denote it by $\mu\in\widetilde{\mathcal{A}}_{\sigma}.$ 
For a sub--$\sigma$--algebra $\mathcal{A}$ of $ \G\otimes \borel{\Rp}\otimes\borel{\R{n}}$, the restriction of the measure $M_{\mu}$ to $(\widetilde{\Omega}, \mathcal{A})$ will be denoted by $M_{\mu}|_{\mathcal{A}}.$
Moreover, for $W:(\widetilde{\Omega},\G\otimes \borel{\Rp}\otimes\borel{\R{n}})\longrightarrow (\R{},\borel{\R{}})$, we define the random measure 
$W\mu$ as follows
$$(W\mu)(\omega;\ds,\dx):= W(\omega,s,x)\mu(\omega;\ds,dx).$$ 

\begin{definition}\label{CondFPredProj}
Let $\mu\in\widetilde{\mathcal{A}}_{\sigma}$ and 
$W:(\widetilde{\Omega},\G\otimes \borel{\Rp}\otimes\borel{\R{n}})\longrightarrow (\R{p},\borel{\R{p}})$ be such that $|W^i|\mu\in\widetilde{\mathcal{A}}_{\sigma},$ for every $i=1,\ldots,p,$ where $W^i$ denotes the $i-$th component of $W$.
Then, the \emph{conditional $\mathbb G-$predictable projection of $W$ on $\mu$}, denoted by $M_{\mu}\big[W|\widetilde{\mathcal{P}}\big],$ is defined componentwise as follows
$$M_{\mu}\big[W|\widetilde{\mathcal{P}}\big]^i:=
  \frac{\ud M_{W^i\mu}|_{\widetilde{\mathcal{P}}}}{\ud M_{\mu}|_{\widetilde{\mathcal{P}}}}, \text{ for }i=1,\ldots,p.$$ 
\end{definition}

\begin{definition}\label{def:OrthogDecomp}
Let $(X^\circ,X^\natural)\in\mathcal H^2(\mathbb R^m)\times\mathcal H^{2,d}(\mathbb R^n)$ and $Y\in\mathcal H^2(\mathbb R^d)$.
The decomposition
\begin{align*}
Y= Y_0 + Z\cdot X^\circ + U\star \widetilde{\mu}^{X^\natural} + N,
\end{align*}
where the equality is understood componentwise, will be called the \emph{orthogonal decomposition of $Y$ with respect to $(X^\circ,X^\natural)$} if
\begin{enumerate}[label=(\roman*)]
  \item $Z\in\mathbb H^2(X^\circ)$ and $U\in\mathbb H^2(\mu^{X^\natural}),$
  \item $ Z\cdot X^\circ \mutorthog U\star \widetilde{\mu}^{X^\natural}$, 
  \item $N\in\mathcal H^2(\mathbb R^d)$ with $\langle N,X^{\circ}\rangle=0$ and $M_{\mu^{X^\natural}}[\Delta N|\widetilde{\mathcal P}]=0.$
\end{enumerate}
\end{definition}

Let $X\in\mathcal H^2(\mathbb R^m)$.
Then \cite[Lemma III.4.24]{jacod2003limit}, which is restated in the coming lines, provides the orthogonal decomposition of a martingale $Y$ with respect to  $(X^c,X^d)$, \emph{i.e.} the natural pair of $X$.

\begin{lemma}
Let $Y\in\mathcal H^2(\mathbb R^d)$ and $X\in\mathcal H^2(\mathbb R^m)$.
Then, there exists a pair $(Z,U)\in\mathbb{H}^2({X}^c)\times\mathbb{H}^2(\mu^X)\footnotemark$ \footnotetext{We assume that $m=n$.} and $N\in\mathcal H^2(\mathbb R^d)$ such that
\begin{align*}
Y= Y_0 + Z\cdot {X}^c + U\star\mutilde^{X} + {N},
\end{align*}
where the equality is understood componentwise, with $\langle {X}^{c},N^{c}\rangle = 0$ and $M_{\mu^{X}}\big[ \Delta {N} |\widetilde{\mathcal{P}}^{\mathbb F} \big]=0$.
Moreover, this decomposition is unique, up to indistinguishability. 
\end{lemma}

\label{Notation:partsofmart}
In the rest of this subsection, we will provide some useful results, which allow us to obtain the orthogonal decomposition as understood in Definition \ref{def:OrthogDecomp}.
Their proofs  are relegated to Appendix \ref{AppendixPreliminaries}.
We also need to introduce at this point some further helpful notation.
\begin{enumerate}[label=\raisebox{1.5pt}{\tiny$\bullet$}]
  \item For a multidimensional process $L$, resp. random variable $\psi$, its $i-$component will be denoted by $L^i$, resp $\psi^i$.
  \item The continuous part of the martingale $X^\circ$ will be denoted by $X^{\circ,c}$. 
  \item The purely discontinuous part of the martingale $X^\circ$ will be denoted by $X^{\circ,d}$. 
  \item $X^{\circ,i}$ denotes the $i-$component of $X^\circ$.
  \item $X^{\circ,c,i}$ denotes the $i-$component of the continuous part of $X^\circ$.
  \item $X^{\circ,d,i}$ denotes the $i-$component of the purely discontinuous part of $X^\circ$.
  \item $X^{\natural,j}$ denotes the $j-$component of $X^\natural.$
\end{enumerate}

\begin{lemma}\label{lem:StochIntOrthog}
Let $(X^\circ,X^\natural)\in\mathcal H^2(\mathbb R^m)\times\mathcal H^{2,d}(\mathbb R^n)$ with $M_{\mu^{X^\natural}}[\Delta X^\circ|\widetilde{\mathcal P}]=0$.
Then, for every $Y^\circ\in\mathcal L^2(X^\circ)$, $Y^\natural\in\mathcal K^2(\mu^{X^\natural})$, we have $\langle Y^\circ,Y^\natural\rangle=0.$
In particular, $\langle X^\circ, X^\natural\rangle=0$.
\end{lemma}

In view of Lemma \ref{lem:StochIntOrthog}, we can provide in the next proposition the desired orthogonal decomposition of a martingale $Y$ with respect to a pair $(X^\circ,X^\natural)\in\mathcal H^2(\mathbb R^m)\times\mathcal H^{2,d}(\mathbb R^n)$, \emph{i.e.} we do not necessarily use the natural pair of the martingale $X$. 
Observe that in this case we do allow the first component to have jumps.
This is particularly useful when one needs to decompose a discrete--time martingale as a sum of an It\=o integral, a stochastic integral with respect to an integer--valued measure and a martingale orthogonal to the space of stochastic integrals. 
\begin{proposition}\label{prop:OrthogDecomp}
Let $Y\in\mathcal H^2(\mathbb R^d)$ and $(X^\circ,X^\natural)\in\mathcal H^2(\mathbb R^m)\times\mathcal H^{2,d}(\mathbb R^n)$ with $M_{\mu^{X^\natural}}[\Delta X^\circ|\widetilde{\mathcal P}]=0$, where the equality is understood componentwise.
Then, there exists a pair $(Z,U)\in\mathbb{H}^2({X}^\circ)\times\mathbb{H}^2(\mu^{X^\natural})$ and $N\in\mathcal H^2(\mathbb R^d)$ such that
\begin{align}\label{eq:orth-deco}
Y= Y_0 + Z\cdot {X}^\circ + U\star\mutilde^{X^\natural} + {N},
\end{align}
where the equality is understood componentwise, with $\langle {X}^{\circ},N\rangle = 0$ and $M_{\mu^{X^\natural}}\big[ \Delta {N} |\widetilde{\mathcal{P}}^{\mathbb F} \big]=0$. 
Moreover, this decomposition is unique, up to indistinguishability. 
\end{proposition}

In other words, the orthogonal decomposition of $Y$ with respect to the pair $(X^\circ,X^\natural)$ is well--defined under the above additional assumption on the jump parts of the martingales $X^\circ$ and $X^\natural$.

\vspace{0.5em}
We conclude this subsection with some useful results. 
Let $\oX:=(X^\circ,X^\natural)\in\mathcal H^2(\mathbb R^m)\times\mathcal H^{2,d}(\mathbb R^n)$ with $M_{\mu^{X^\natural}}[\Delta X^\circ|\widetilde{\mathcal P}]=0$.
Then we define
\begin{align*}
\mathcal{H}^2(\oX^\perp)  :=\big(\mathcal{L}^2(X^\circ)\oplus\mathcal{K}^2(\mu^{X^\natural})\big)^\perp.
\end{align*}  
If $(X^c,X^d)$ is the natural pair of $X\in\mathcal H^2(\mathbb R^m)$, then we define $\mathcal H^2(X^\perp):=\mathcal{H}^2((X^c,X^d)^\perp)$.
In view of the previous definitions, we will abuse notation and we will denote the natural pair of $X$ by $X$ as well.

\begin{proposition}\label{prop:CharacterOrthogSpace}
Let $\oX:=(X^\circ,X^\natural)\in\mathcal H^2(\mathbb R^m)\times\mathcal H^{2,d}(\mathbb R^n)$ with $M_{\mu^{X^\natural}}[\Delta X^\circ|\widetilde{\mathcal P}]=0$.
Then,
\begin{align*}
\mathcal{H}^2(\oX^\perp) = \big\{ L\in\mathcal H^2(\mathbb R^d), \langle X^{\circ}, L\rangle=0 \text{ and } M_{\mu^{X^\natural}}[\Delta L|\widetilde{\mathcal P}]=0\big\}.
\end{align*}
Moreover, the space $\big(\mathcal{H}^2(\oX^\perp), \Vert \cdot\Vert_{\mathcal H^2(\mathbb R^d)}\big)$ is closed.
\end{proposition}

\begin{corollary}\label{cor:SpaceDecomposition}
Let $\oX:=(X^\circ,X^\natural)\in\mathcal H^2(\mathbb R^m)\times\mathcal H^{2,d}(\mathbb R^n)$ with $M_{\mu^{X^\natural}}[\Delta X^\circ|\widetilde{\mathcal P}]=0$.
Then, 
\begin{align*}
\mathcal H^2(\mathbb R^p) = \mathcal L^2(X^\circ) \oplus \mathcal K^2(\mu^{X^\natural})\oplus \mathcal H^2(\oX^\perp),
\end{align*}
where each of the spaces appearing in the above identity is closed.
\end{corollary}

\begin{remark}
The aim of this paper is to provide a general result for the existence and uniqueness of the solution of a BSDE.
In other words, our result should also cover the case where the underlying filtration is not quasi--left--continuous, see \cite[Definition 3.39]{he1992semimartingale} and \cite[Theorem 3.40, Theorem 4.35]{he1992semimartingale}.
In the same vein, we should be able to choose the It\=o integrator to be an arbitrary square integrable martingale, \emph{i.e.} not necessarily quasi--left--continuous, not to mention continuous.
On the other hand, it is well--known that the orthogonal decomposition of martingales is hidden behind the definition of the contraction mapping used to prove wellposedness of solutions to BSDEs. 
Therefore, we have by Proposition \ref{prop:OrthogDecomp} a sufficient condition on the jumps of the two stochastic integrators $X^\circ$ and $X^\natural,$ so that we can obtain the orthogonal decomposition in such a general framework.
\end{remark}

%% file: suitable_spaces.tex

Let us first define the maps $\mathbb R^{n}\ni x\overset{\mathrm q}{\longmapsto}xx^\top\in\mathbb R^{n\times n}$ and $\mathbb R^{n}\ni x\overset{\mathrm I}{\longmapsto}x\in\mathbb R^{n}$.
Next, we provide a result which justifies the validity of Assumption \ref{assumptionC} below, which is based on \cite[Property II.1.2 and Proposition II.2.9]{jacod2003limit}.

\begin{lemma}\label{prop:II.2.9-jacod2003limit}
Let $\oX:=(X^{\circ}, X^{\natural})\in\mathcal H^2(\mathbb R^m)\times\mathcal H^{2,d}(\mathbb R^n)$. 
Then, there exists a predictable, non-decreasing and \cadlag process $C^{\oX}$ such that 
\begin{enumerate}[label=$(\roman*)$]
\item\label{prop:II.2.9-jacod2003limit-PQC} Each component of $\langle X^{\circ} \rangle$ is absolutely continuous with respect to $C^{\oX}$.  
  In other words, there exists a predictable, positive definite and symmetric $m\times m-$matrix $\frac{\textup{d} \langle X^{\circ}\rangle}{\textup{d} C^{\oX}}$ such that for any $1\leq i,j\leq m$
  \begin{align*}
    \langle X^{\circ} \rangle^{ij}_\cdot  =   \int_{(0,\cdot]} \frac{\textup{d}\langle X^{\circ}\rangle_s^{ij}}{\textup{d} C^{\oX}_s}\, \textup{d} C^{\oX}_s.
  \end{align*}
\item\label{prop:II.2.9-jacod2003limit-Disintegration} The \emph{disintegration property}\index{Disintegration property} given $C^{\oX}$ holds for the compensator $\nu^{X^{\natural}}$, 
  \emph{i.e.} there exists a transition kernel $K^{\oX}:(\Omega\times \mathbb R_+,\mathcal P)\longrightarrow \mathcal{R}(\mathbb R^n,\mathcal B(\mathbb R^n))$, 
  where $\mathcal{R}\big(\mathbb R^n,\mathcal B(\mathbb R^n)\big)$ is the space of Radon measures on $\big(\mathbb R^n,\mathcal B(\mathbb R^n)\big)$, such that
  \begin{align*}
    \nu^{X^{\natural}}(\omega;\textup{d} t,\textup{d} x)  =   K_t^{\oX}(\omega;\textup{d} x)\, \textup{d} C^{\oX}_t.
  \end{align*}
\item\label{prop:II.2.9-jacod2003limit-QLC} $C^{\oX}$ can be chosen to be continuous if and only if $\oX$ is $\mathbb G-$quasi-left-continuous. 
\end{enumerate}
\end{lemma}
\begin{proof}
We can follow exactly the same arguments as in \cite[Proposition II.2.9]{jacod2003limit} for the process
\begin{equation*}
  C^{\oX}_\cdot:= \sum_{i,j=1}^m\textup{Var}\big(\langle X^{\circ}\rangle^{ij}\big)_\cdot + (\vert\mathrm{I}\vert^2\wedge 1)\star\nu^{X^\natural}_\cdot,
\end{equation*}
where $\textup{Var}(A)$ denotes the total variation process of the finite variation process $A.$
We need to underline that under our framework the process $\langle X^\circ\rangle$ is not necessarily continuous.
However, we can indeed follow the same arguments as in \cite[Proposition II.2.9]{jacod2003limit}.
\end{proof}
\begin{assumption}{C}\label{assumptionC}
Let $\oX:=(X^\circ,X^\natural)\in\mathcal H^2(\mathbb R^m)\times\mathcal H^2(\mathbb R^n)$ and $C$ be a predictable, c\`adl\`ag and increasing process. 
The pair $(\oX,C)$ satisfies Assumption \ref{assumptionC} if 
each component of $\langle X^\circ \rangle$ is absolutely continuous with respect to $C$ and if the disintegration property given $C$ holds for the compensator $\nu^{X^\natural}$. 
\end{assumption}

\begin{remark}
Let $X\in\mathcal H^2(\mathbb R^p)$.
Recall that we have abused notation and we denote its natural pair $(X^c,X^d)$ by $X$ as well.
Then there exist several possible choices for $C^X$ such that Assumption \ref{assumptionC} is satisfied. 
In {\rm \cite[Proposition II.2.9]{jacod2003limit}}, for example, the following process is used
\[
	\widetilde C^{X}:=\sum_{i,j=1}^n 
\Var\left(\,\pqc{X^{c,i},X^{c,j}}\right) + (\vert\mathrm I\vert\wedge 1)^2\star\nu^{X},
\]
while one could also take
\[
	\overline C^{X}:=\tr[\,\pqc{X^{c}}] + \vert \mathrm I \vert^2\star\nu^{X}.
\]
\end{remark}

\begin{remark}
Let $\oX:=(X^\circ,X^\natural)\in\mathcal H^2(\mathbb R^m)\times\mathcal H^2(\mathbb R^n)$ and consider a pair $(\oX,C^{\oX})$ which satisfies Assumption \ref{assumptionC}.
Then, the Radon-Nikod\'ym derivative  $\frac{\ud \langle X^\circ\rangle}{\ud C^{\oX}}$ is $\mathbb G-$predictable, positive definite and symmetric. 
Indeed, the predictability and the positive definiteness follows from {\rm\cite[Statement~III.6.2]{jacod2003limit}}, while the symmetry is immediately inherited from the symmetry of $\langle X^\circ\rangle$. 
The above properties enable us to define the following $\mathbb G-$predictable process $c^{\oX}$ 
\begin{equation}\label{srrn}
c^{\oX} := \left( \frac{\ud \langle X^\circ\rangle}{\ud C^{\oX}}\right)^{\frac{1}{2}}.
\end{equation} 
In addition, if we define the random measure $\mu^{\oX}_{\Delta}$ on 
$(\Rp,\borel{\Rp})$, for any $t\ge0$, via 
\begin{equation}\label{jumprn}
\mu^{\oX}_{\Delta}\left(\omega;[0,t]\right)
  :=\sum_{0<s\leq t} \left( \Delta C^{\oX}_s(\omega)\right)^2  
 \text{ then it holds }
\frac{\ud \mu^{\oX}_{\Delta}}{\ud C^{\oX}}(t)= \Delta C^{\oX}_t.
\end{equation}
\end{remark}

Assume now that ${\oX}:=(X^\circ,X^\natural)\in\mathcal H^2(\mathbb R^n)\times\mathcal H^2(\mathbb R^m)$ with $M_{\mu^{X^\natural}}[\Delta X^\circ|\widetilde{\mathcal P}]=0$. 
Assume, moreover, that there exists a process $C^{\oX}$ such that $({\oX},C^{\oX})$ satisfies Assumption \ref{assumptionC}.
The process $\oX$ is not assumed to be quasi--left--continuous, hence it is possible that it has fixed times of discontinuities. 
Using \cite[Proposition~II.2.29.b]{jacod2003limit}, we 
have that $\oX\in\mathcal H^2(\mathbb R^m)\times\mathcal H^2(\mathbb R^n)$ if and only if the following holds
$$\Expect{\textrm{Tr}\big[\langle X^\circ\rangle_\tau\big] + \norm{\mathrm q}^2\star\nu_\tau^{X^\natural}}<\infty.$$

Take into account, now, that $X^\natural \in\mathcal H^2(\mathbb R^m)$, the fact that the predictable projection of $\Delta X^\natural$ is indistinguishable from the zero process, see \cite[Corollary I.2.31]{jacod2003limit}, and 
Property \eqref{PropertyII-1-11}. 
All the above yield 
\begin{align*}
\int_{\mathbb R^n}x\nu^{X^\natural}(\{s\}\times \ud x) = 0.
\end{align*} 
Therefore, the predictable quadratic variation of $\oX$ admits, by Lemma \ref{lem:StochIntOrthog} and \cite[Theorem II.1.33]{jacod2003limit}, the following 
representation
\begin{align*}
\pqc{\oX}_\cdot 
&=
\begin{bmatrix}
 \pqc{X^\circ}_\cdot &  0\\
  0 & \displaystyle{\mathrm q}\star\nu^{X^\natural}_\cdot 
\end{bmatrix}.
\end{align*}
However, the reader should keep in mind that for the arbitrary element $W\star \widetilde{\mu}^{X^\natural}$ of $\mathcal K^2(\mu^{X^\natural})$ its predictable quadratic variation is represented as
\begin{equation*}
\langle W\star \widetilde{\mu}^{X^\natural} \rangle_\cdot = \{(W-\widehat{W}^{X^\natural})(W-\widehat{W}^{X^\natural})^\top\}\star\nu^{X^\natural}_\cdot +\sum_{s\le \cdot} \{(1-\zeta^{X^\natural}_s)\widehat{W}^{X^\natural}\big(\widehat{W}^{X^\natural}\big)^\top\};
\end{equation*}
use the polarization identity and \cite[Theorem II.1.33]{jacod2003limit}.\footnote{The reader may recall \eqref{notation:U-hat} and \eqref{notation:1-hat} for the definition of the process $\zeta^{X^\natural}$ and of $\widehat{W}^{X^\natural}$.}
If, in addition,
 $\mathbb E\big[ \sum_{s\ge 0} (W(s,\Delta X^\natural_s)\big)^2\big]<\infty$, then 
\begin{equation*}
\langle W\star \widetilde{\mu}^{X^\natural} \rangle_\cdot = (WW^\top)\star\nu^{X^\natural}_\cdot -\sum_{s\le \cdot} \left(\widehat{W}^{X^\natural}\big(\widehat{W}^{X^\natural}\big)^\top\right);
\end{equation*}
see \cite[Theorem 11.21]{he1992semimartingale} or \cite[Theorem 13.3.16]{cohen2015stochastic}, where we use again the polarization identity in order to conclude in the multidimensional case.

\vspace{0.5em}
Let $U$ be a $\mathbb G-$predictable function taking values in $\mathbb R^{d}$.
Then we define, abusing notations slightly,
\[
\widehat K_{t}^{\oX}(U_t(\omega;\cdot))(\omega)
  := \int_{\R n}U_t(\omega;x) K_t^{\oX}(\omega;\ud x),\ t\geq 0,
\]
where $K^{\oX}$ is the transition kernel from Assumption \ref{assumptionC}.
Using Assumption \ref{assumptionC} and \eqref{jumprn}, we get that
\[
\widehat{U}^{X^\natural}_t(\omega)=\int_{\R n}U(\omega,t,x)\nu^{X^\natural}(\omega;\{t\}\times \ud x)
  = \widehat K_{t}^{\oX}(U_t(\omega;\cdot))
    (\omega)\Delta C^{\oX}_t(\omega),\ t\geq 0.
\]

Using the previous definitions and results, we can rewrite $\langle \oX \rangle$ as 
follows
\begin{align*}				
\langle \oX\rangle_\cdot	
&=
\begin{bmatrix}
\displaystyle\int_0^{\cdot} c^{\oX}_s(c_s^{\oX})^{\transp} \ud C^{\oX}_{s}& 0\\ 
0 &
\displaystyle\int_{(0,\cdot]\times\R n}xx^{\transp}K_s^{\oX}(\ud x) \ud C^{\oX}_{s}
\end{bmatrix}
=
\begin{bmatrix}
\displaystyle \int_0^{\cdot}  c^{\oX}_s (c_s^{\oX})^{\transp} \ud C^{\oX}_{s}
&0\\ 
0&\displaystyle\int_0^\cdot  \widehat K_{s}^{\oX}({\mathrm q})  
\ud C^{\oX}_{s}
\end{bmatrix}.
\end{align*}
On the other hand, for the predictable quadratic variation of $Z\cdot X^\circ\in\mathcal L^2(X^\circ)$ we have by \cite[Theorem III.6.4]{jacod2003limit}
\begin{align*}
\langle Z\cdot X^\circ\rangle_\cdot 
  &= \int_0^\cdot Z_s \frac{\ud \langle X^\circ \rangle_s}{\ud C^{\oX}_s}Z_s^\top \, \ud C^{\oX}_s
   = \int_0^\cdot  (Z_sc^{\oX}_s) (Z_sc^{\oX}_s)^\top\, \ud C^{\oX}_s
\end{align*}
and for the predictable quadratic variation of $W\star\widetilde{\mu}^{X^\natural}\in\mathcal K^2(\mu^{X^\natural})$ we have by the definitions and comments above that 
\begin{align*}
\langle W\star \widetilde{\mu}^{X^\natural} \rangle_\cdot 
  &= \{(W-\widehat{W}^{X^\natural})(W-\widehat{W}^{X^\natural})^\top\}\star\nu^{X^\natural}_\cdot +\sum_{s\le \cdot} \{(1-\zeta^{X^\natural}_s)\widehat{W}^{X^\natural}\big(\widehat{W}^{X^\natural}\big)^\top\}\\
  &= \int_0^\cdot \widehat{K}_s^{\oX}\big((W_s(\cdot)-\widehat{W}^{X^\natural}_s)(W_s(\cdot)-\widehat{W}_s^{X^\natural})^\top\big) \ud C_s^{\oX} 
  +\sum_{s\le \cdot} \Big\{(1-\zeta^{X^\natural}_s)\widehat{K}_s^{\oX}(W_s(\cdot))\Big(\widehat{K}_s^{\oX}(W_s(\cdot))\Big)^\top\big(\Delta C^{\oX}_s\big)^2\Big\}\\
    &= \int_0^\cdot   \Big\{\widehat{K}_s^{\oX}\big((W_s(\cdot)-\widehat{W}^{X^\natural}_s)(W_s(\cdot)-\widehat{W}^{X^\natural}_s)^\top\big)  
  + (1-\zeta^{X^\natural}_s)\Delta C_s^{\oX}\widehat{K}_s^{\oX}(W_s(\cdot))\Big(\widehat{K}_s^{\oX}(W_s(\cdot))\Big)^\top\Big\} \ud C^{\oX}_s.
  \numberthis\label{repr:pqc-PDM}
\end{align*}

We proceed now by defining the spaces that will be necessary for our analysis, see also \cite{elkaroui1997general}. 
Let $\beta\geq0$ and $A:(\Omega\times \Rp,\mathcal{G}\otimes 
\borel{\Rp})\longrightarrow \Rp$ be a \cadlag, increasing and measurable 
process. 
We then define the following spaces, where the dependence on $A$ is suppressed 
for ease of notation:\label{def-spaces}
\begin{align*}
\mathbb{L}^{2}_{\beta}(\G_\tau)
  & : =   \set{\xi,\ \textrm{$\mathbb R^d-$valued, \ $\G_\tau-$measurable},
\ \norm \xi _{\mathbb{L}^{2}_{\beta}(\G_\tau;\R{d})}^2 
  := \E\left[ \e {\beta A_\tau} \abs \xi^2\right] <\infty}, \\
\mathcal{H}^{2}_{\beta}
  &:=\left\{ M\in\mathcal{H}^{2}, \ 
  \norm M_{\mathcal{H}^{2}_{\beta}}^2
  :=\mathbb{E} \left[\int_0^\tau \e {\beta A_t} \ud {\rm Tr} 
    \left[\langle M\rangle_t\right]\right]<\infty \right\}, \\
\mathbb H^{2}_{\beta}(\oX) &:= 
\Bigg\{\phi \textrm{ is an $\mathbb R^d-$valued 
  $\mathbb G-$optional semimartingale with c\`adl\`ag paths and}\\
  &\qquad\quad \norm \phi _{\mathbb H^{2}_{\beta}(\oX)}^2:= 
  \Expect{ \int_0^\tau\e{\beta A_t} \abs{\phi_t}^2 \ud C^{\oX}_{t}}< \infty\Bigg\},\\
\mathcal S^{2}_{\beta}(\oX) &:= 
  \Bigg\{\phi \textrm{ is an $\mathbb R^d-$valued 
  $\mathbb G-$optional semimartingale with c\`adl\`ag paths and}\\
  &\qquad\quad \norm {\phi}_{\mathcal S^{2}_{\beta}(\oX)}^2:= 
    \Expect{ \sup_{t\in\stochint{0}{\tau}}\e {\beta A_t} \abs 
    {\phi_t}^2}<\infty\Bigg\},\\
\mathbb H^{2}_{\beta}(X^\circ) 
  &:= \set{Z\in \mathbb H^{2}(X^\circ), \ 
    \norm Z_{\mathbb H^{2}_{\beta}(X^\circ)}^2 
  :=\Expect{\int_0^\tau \e{\beta A_t} \ud {\rm Tr}[\,\pqc{Z\cdot X^\circ}_t]} 
  <\infty}, \\
\mathbb H^{2}_{\beta}(X^\natural) 
  &:= \bigg\{U\in \mathbb H^{2}(X^\natural) , \ 
    \norm{U}_{\mathbb H^{2}_{\beta}(X^\natural) }<\infty,\text{ with } \norm{U}_{\mathbb 
  H^{2}_{\beta}(X^\natural) }:=\Expect{\int_0^{\tau}  \e{\beta 
  A_t}\ud {\rm Tr}[\langle U\star\mutilde^{X^\natural}\rangle_t]}\bigg\},\\
\mathcal{H}^2_{\beta}(\oX^\perp)	
  &:=\set{ M\in\mathcal{H}^2(\oX^\perp), \ 
  \norm M^2_{\mathcal{H}^2_{\beta}(\oX^\perp)}:= \Expect{\int_0^\tau \e 
{\beta A_t} \ud {\rm Tr}[\,\pqc{M}_t]}<\infty}.
\end{align*}

\vspace{0.5em}

Finally, for $(Y,Z,U,N)\in 	\mathbb H^{2}_{\beta}(\oX)\times
				\mathbb H^{2}_{\beta}(X^\circ)\times
				\mathbb H^{2}_{\beta}(X^\natural)\times
				\mathbb H^{2}_{\beta}(\oX^\perp)$
and assuming that $A=\alpha^2 \cdot C^{\oX}$ for a measurable process $\alpha:(\Omega \times \Rp, \mathcal{G} \otimes \borel{\Rp})\longrightarrow \mathbb R$,
we define				\begin{align*}
					\norm {(Y,Z,U,N)} _{\beta,\oX}^2 &:= 
					\norm {\alpha Y}_{\mathbb 
H^{2}_{\beta}(\oX)}^2 
					+	\norm{Z}_{\mathbb 
H^{2}_{\beta}(X^\circ)}^2 
					+ 	\norm{U}_{\mathbb 
H^{2}_{\beta}(X^\natural)}^2 
					+	\norm N _{\mathbb 
H^{2}_{\beta}({\oX}^\perp)}^2,
				\end{align*}
			and for
				$(Y,Z,U,N)\in 	\mathcal 
S^{2}_{\beta}(\oX)\times
				\mathbb H^{2}_{\beta}(X^\circ)\times
				\mathbb H^{2}_{\beta}(X^\natural)\times
				\mathbb H^{2}_{\beta}(\oX^\perp)$, 
we define
				\begin{align*}
					\norm {(Y,Z,U,N)} 
_{\star,\beta,\oX}^2 &:=
					\norm {Y}_{\mathcal 
S^{2}_{\beta}(\oX)}^2 
					+	\norm{Z}_{\mathbb 
H^{2}_{\beta}(X^\circ)}^2 
					+ 	\norm{U}_{\mathbb 
H^{2}_{\beta}(X^\natural)}^2 
					+	\norm N _{\mathbb 
H^{2}_{\beta}(\oX^\perp)}^2.
				\end{align*}

The next lemma will be useful for future computations and, in addition, 
justifies the definition of the norms on the spaces provided above.

\begin{lemma}\label{lemma1}
Let $(Z,U)\in \mathbb H^{2}_{\beta}(X^\circ)\times \mathbb H^{2}_{\beta}(X^\natural)$. 
Then 
\begin{align}		
		\norm Z_{\mathbb H^{2}_{\beta}(X^\circ)}^2 
      &= \E\left[ \int_0^\tau \e {\beta A_t} \norm {c_t Z_t}^2\ud C^{\oX}_{t}\right], \label{normeq1}\\		
	  \norm{U}_{\mathbb H^{2}_{\beta}(X^\natural)}^2	
      &= \E\left[ \int_0^\tau\e {\beta A_t} \big(\tnorm {U_t(\cdot)}_t^{\oX}\big)^2 \ud C^{\oX}_{t}\right], \label{normeq2}
\end{align}
where for every $(t,\omega)\in\mathbb R_+\times\Omega$
$$\left(\tnorm{U_t(\omega;\cdot)}^{\oX}_t(\omega)\right)^2 
  := \widehat K_{t}^{\oX}(|U_t(\omega;\cdot) - \widehat{U}^{X^\natural}_t(\omega)|^2)(\omega) 
  +  (1-\zeta^{X^\natural}_t)\Delta C^{\oX}_t(\omega) |\widehat K_{t}^{\oX}(U_t(\omega;\cdot))(\omega)|^2\geq 0.\footnotemark$$
Furthermore \footnotetext{The process $\zeta$ has been defined in \eqref{notation:1-hat}.}
\begin{align*}
\big\Vert Z\cdot X^{\circ} + U\star\widetilde{\mu}^{X^\natural}\big\Vert_{\Hcal^2_\beta}^2 
  &=
	\norm{Z}_{\mathbb H^{2}_{\beta}(X^{\circ}) } +
	\norm{U}_{\mathbb H^{2}_{\beta}(X^{\natural}) }.
\end{align*}
\end{lemma}

\begin{proof}
Let $Z\in \mathbb H^{2}_{\beta}(X^{\circ})$, then using \cite[Theorem~III.6.4]{jacod2003limit}, we get that $Z\cdot X^{\circ}\in \mathcal{H}^{2}$  with 
\begin{equation}\label{pqceq1}
\pqc{Z\cdot X^{\circ}}
  = Z\frac{\ud \langle X^\circ \rangle}{\ud C^X}Z^{\transp}\cdot C^X 
  = Zc^{\oX}(c^{\oX})^{\transp}Z^{\transp} \cdot C^X,
\end{equation}
for $c^{\oX}$ as introduced in \eqref{srrn}. The first result is then obvious.

\vspace{0.5em}
Now, let $U\in\mathbb H^{2}_{\beta}(X^{\natural})$.
Then by the previous computations, we have 
\begin{align*}
\langle U\star\mutilde^{X^\natural}\rangle_\cdot
  &= \int_0^{\cdot} \Big\{\widehat{K}_s^{\oX}\big((U_s(\cdot)-\widehat{U}^{X^\natural}_s)(U_s(\cdot)-\widehat{U}^{X^\natural}_s)^\top\big)  
  + (1-\zeta^{X^\natural}_s)\Delta C_s^{\oX}\widehat{K}_s^{\oX}(U_s(\cdot))\Big(\widehat{K}_s^{\oX}(U_s(\cdot))\Big)^\top\Big\} 
  \ud C^{\oX}_{s},\numberthis \label{pqceq2}
	\end{align*}
from which the second result is also clear. 
Moreover, we have 
\begin{align*}
\norm{Z\cdot X^{\circ} + U\star\widetilde{\mu}^{X}}_{\Hcal^2_\beta}^2 
&= \Expect{\int_0^\tau \e {\beta A_t}  \ud {\rm Tr}\left[\big\langle Z\cdot X^{\circ}+U\star\widetilde{\mu}^{X^\natural}\big\rangle_t\right]}\\
	&=	\Expect{\int_0^\tau \e {\beta A_t}  \ud {\rm 
Tr}\left[\langle Z\cdot X^{\circ}\rangle _t\right]+\int_0^\tau \e {\beta A_t}  \ud {\rm 
Tr}\left[\big\langle U\star\widetilde{\mu}^{X^\natural}\big\rangle_t\right]},
	\end{align*}
where the second equality holds due to Lemma \ref{lem:StochIntOrthog}.

\vspace{0.5em}
Notice finally that the process ${\rm Tr}[\langle U\star\mutilde^{X^\natural}\rangle]$ is 
non-decreasing, and observe that 
\begin{equation}\label{eq:norm}
\Delta {\rm Tr}\left[\langle U\star\mutilde^{X^\natural}\rangle_t\right]
  = \Big(\tnorm{U_t(\cdot)}_t^{\oX}\Big)^2 \Delta C^{\oX}_t,\ t\geq 0.
\end{equation}
Since $C^{\oX}$ is non-decreasing, we can deduce that $\tnorm{U_t(\cdot)}_t^{\oX}\geq 0$.
\end{proof}

We conclude this section with the following convenient result. 
Define the following space for $\ud C^{\oX}\otimes \ud \mathbb P-a.e.$ 
$(t,\omega)\in\mathbb R_+\times\Omega$ 
$$\mathfrak H_{t,\omega}^{\oX}
  :=\overline{\left\{\mathcal U:(\mathbb R^n,\mathcal B(\mathbb 
    R^n))\longrightarrow (\mathbb R^d,\mathcal B(\mathbb R^d)),\ 
    \tnorm{\ \mathcal U(\cdot)}^{\oX}_t(\omega)<\infty\right\}}.$$
Define also
$$\mathfrak H^{\oX}:=\left\{U:[0,T]\times\Omega\times\mathbb 
  R^n\longrightarrow\mathbb R^d,\ U_t(\omega;\cdot)\in \mathfrak 
  H_{t,\omega}^{\oX}, \text{ for $\ud C^{\oX}\otimes \ud \mathbb P-a.e.$ 
  $(t,\omega)\in\mathbb R_+\times\Omega$}\right\}.$$

\begin{lemma}\label{lem:polish}
The space $\left(\mathfrak H_{t,\omega}^{\oX},\tnorm{\cdot}^{\oX}_t(\omega)\right)$ 
is Polish, for $\ud C^{\oX}\otimes \ud\mathbb P-a.e.$ $(t,\omega) \in \mathbb R_+\times\Omega$.
\end{lemma}

\begin{proof}
The fact that $\tnorm{\cdot}^{\oX}_t(\omega)$ is indeed a norm is immediate from \eqref{eq:norm} and Kunita--Watanabe's inequality. 
Then, the space is clearly Polish since the measure $K_t^{\oX}(\omega;\ud x)$ is regular; it integrates $x\longmapsto |x|^2$ for $\ud C^{\oX}\otimes \ud \mathbb P-a.e.$ 
$(t,\omega)\in\mathbb R_+\times\Omega$.
\end{proof}

%% file: A_useful_lemma_for_generalised_inverses.tex
In the following sections we will need a result on generalized inverses which is stated as a corollary of the following lemma. The proof is presented in Appendix \ref{AppendixA}.

\begin{lemma}\label{mainlemma}
	Let $g$ be a non-decreasing sub-multiplicative function on $\mathbb{R}_+$, that is to say 
	$$g(x+y)\leq \ell g(x)g(y),$$ for some $\ell>0$ and for every $x,y\in\mathbb{R}_+$. 
	Let $A$ be a c\`adl\`ag  and non-decreasing function and define its left-continuous inverse $L$ by
	\[
		L_s:=\inf\left\{t\geq 0,\  A_t \geq s  \right\}.
	\]
	Then it holds that
	\[
		\int_0^t g(A_s) \dA s \leq \ell g\left(\max_{\set{s,\, L_s<\infty}}\Delta A_{L_s}\right)\int_{A_0}^{A_t} g(s) \ds.
	\]
\end{lemma}

\begin{corollary}\label{maincorollary}
Let $A$ and $g$ as in Lemma \ref{mainlemma} with the additional assumption that $A$ has uniformly bounded jumps, say by $K$. 
Then there exists a universal constant $K'>0$ such that
\[
	\int_0^t g(A_s) \dA s \leq K' \int_{A_0}^{A_t} g(s) \ds.
\]
The constant $K'$ equals $\ell g(K)$, where $\ell$ is the sub-multiplicativity constant of $g$.
\end{corollary}

%% file: notation.tex

In this section, we will work on the complete stochastic basis $(\Omega, \mathcal G, \mathbb G, \mathbb P)$ and fix throughout 
\begin{enumerate}[label=\textbullet]
\item a $\mathbb G-$stopping time $T$,
\item an $\mathbb R^{m+n}-$valued, $\mathcal G\otimes\mathcal B(\mathbb R_+)-$measurable process $\oX:=(X^\circ,X^\natural)$ such that 
$$\oX^T\footnotemark\in\mathcal H^2(\mathbb R^m)\times\mathcal H^{2,d}(\mathbb R^n)
\footnotetext{As usual, for a measurable process $X$, the corresponding process stopped at $T$, denoted by $X^T$, is defined by $X^T_t:=X_{t\wedge T}$, $t\geq 0$.} 
\text{ with } M_{\mu^{(X^\natural)^T}}\big[\Delta \big((X^\circ)^T\big)\big|\widetilde{\mathcal P}\big]=0,$$ 
\item a non--decreasing, predictable and c\`adl\`ag process $C^{\oX}$ such that the pair $\big({\oX}^T,(C^{\oX})^T\big)$ satisfies Assumption \ref{assumptionC}.
\end{enumerate}
Abusing notation, we will refer to the stopped processes $X^T$ and $(C^{\oX})^T$ as simply the processes $\oX$ and $C^{\oX}$, since $T$ is given. 
Hence, the time interval on which we will be working throughout this section will always be the stochastic time interval $\llbracket0,T\rrbracket$. 
In addition, a non-decreasing process $A$ will be fixed below, see \eqref{eq:A}. 
In order to simplify notation, and since there is no danger of confusion, we will omit $\oX$ from our spaces and norms. Therefore they become, for any $\beta\geq 0$ and $(t,\omega)\in\mathbb R_+\times\Omega$, 
$$
\mathbb{L}^2_{\beta}:=\mathbb{L}^{2}_{\beta}(\mathcal G_{T}),\
\mathbb{H}^{2}_{\beta}:=\mathbb{H}^{2}_{\beta}(\oX), \ 
\mathcal S^2_\beta:=\mathcal S^2_\beta(\oX), \
\mathbb{H}^{2,\circ}_{\beta}:=\mathbb H^2_\beta(X^\circ), \
\mathbb{H}^{2,\natural}_{\beta}:=\mathbb H^2_\beta(X^\natural),
$$
$$
\mathcal H^{2,\perp}_{\beta}:=\mathcal H^2_\beta(\oX^\perp), \ 
\mathfrak H_{t,\omega}:=\mathfrak H_{t,\omega}^{\oX}, \ 
\mathfrak H:=\mathfrak H^{\oX},$$ 
$$
\norm{\cdot}_\beta:=\norm{\cdot}_{\beta,\oX}, \ 
\norm{\cdot}_{\star,\beta}:=\norm{\cdot}_{\star,\beta,\oX}, \ 
\tnorm{\cdot}_t:=\tnorm{\cdot}_t^{\oX},
$$ 
$$C:=C^{\oX}, \ c:=c^{\oX}, \ \mu^\natural:=\mu^{X^\natural}, \ \nu^\natural:=\nu^{X^\natural}, \ \widetilde \mu^\natural:=\widetilde \mu^{X^\natural}$$
When $\beta=0$, we also suppress it from the notation of the previous spaces. 

\vspace{1em}
We are interested in proving existence and uniqueness of the solution of a backward stochastic differential equation of the form
\begin{equation}\label{BSDEintegralform}
	Y_t = \xi 	+ \int_t^T f(s,Y_s,Z_s,U_s(\cdot)) \dC s 
				- \int_t^T Z_s \ud X^\circ_s 
				-\int_t^T\int_{\R{n}}U_s(x) \widetilde\mu^\natural(\ud s,\ud x)
				- \int_t^T \ud N_s,
\end{equation}
which means that, given the data $(\oX,\Fil,T,\xi,f,C)$, we seek a quadruple $(Y,Z,U,N)$ that satisfies equation \eqref{BSDEintegralform}, $\Pm-a.s.$ 
The martingale $\oX$ is not assumed quasi-left-continuous and may have stochastic discontinuities. 
As a result, the process $C$ may also have discontinuities.
In other words, we consider BSDEs with jumps that are driven both by continuous--time and by discrete--time martingales in a unified framework.

%% file: Formulation_of_the_problem.tex
The data of the BSDE should satisfy the following conditions:
\begin{enumerate}[label=\bf{(F\arabic*)},leftmargin=*,itemindent=0.0cm]
	\item\label{Fone} The martingale $\oX$ belongs to $\mathcal H^2(\mathbb R^m)\times \mathcal H^2(\mathbb R^n)$ and $(\oX,C)$ satisfies Assumption \ref{assumptionC}.
	\item\label{Ftwo} The terminal condition satisfies $\xi\in \mathbb{L}^2_{\hat{\beta}}$ for some $\hat\beta>0$. 
	\item\label{Fthree} The generator\footnote{This is also called \textit{driver} of the BSDE.} of the equation $f:\domain $ is such that for any 
				$(y,z,u)\in\R d\times \R {d\times m}\times \mathfrak H$, 
				the map 
				$$(t,\omega)\longmapsto f(t,\omega,y,z,u_t(\omega;\cdot))\ \text{is $\mathcal F_t\otimes\mathcal B([0,t])-$ measurable.}$$ 
				Moreover, $f$ satisfies a stochastic Lipschitz condition, that is to say there exist 
				\[
					r:\left(\Omega\times\mathbb{R}_+,\Pred\right)\longrightarrow(\Rp,\borel{\Rp})\ \textrm{ and }\
					\vartheta =(\theta^\circ,\theta^\natural):\left(\Omega\times\mathbb{R}_+,\Pred\right)\longrightarrow(\Rp^2,\borel{\Rp^2}),
				\] such that, for $\ud C \otimes \ud \Pm-a.e.$ $(t,\omega)\in\mathbb R_+\times\Omega$
				\begin{align}\label{fstochLip}
				\begin{split}
					& \quad \left|f(t,\omega,y,z,u_t(\omega;\cdot))-f(t,\omega,y',z',u_t'(\omega;\cdot))\right|^2\\
							&\hspace{2em}\leq 	r_t(\omega) |y-y'|^2 	+ \theta^\circ_t(\omega)\Vert c_t(\omega) (z-z')\Vert^2 + \theta^\natural_t(\omega) \left(\tnorm{u_t(\omega;\cdot)-u'_t(\omega;\cdot)}_t(\omega)\right)^2.		
				\end{split}
				\end{align}
	\item\label{Ffour} Let\footnote{We assume, without loss of generality, that $\alpha_t>0$, $\ud C\otimes \ud\Pm-a.e.$} 
				$\alpha^2_\cdot:= \max\{\sqrt{r_\cdot}, \theta^\circ_\cdot,\theta^\natural_\cdot\}$
				and define the increasing, $\mathbb G-$predictable and \cadlag \ process
				\begin{equation}\label{eq:A}
					A_\cdot :=\int_0^{\cdot} \alpha_s^2 \dC{s}.
				\end{equation}
				\noindent Then there exists $\Phi>0$ such that
				\begin{equation}\label{Ffouruniformjumps}
					\Delta A_t(\omega) \leq \Phi, \ \text{for $\ud C \otimes \ud \Pm-a.e.$ $(t,\omega)\in\mathbb R_+\times\Omega$}.
				\end{equation}
	\item\label{Ffive} We have for the same $\hat\beta$ as in \ref{Ftwo}
 		 	\[
 				\E\left[\int_0^T\e{\hat{\beta} A_t} \frac{\abs{f(t,0,0,\mathbf 0)}^2}{\alpha^2_t} \ \dC{t}\right] <\infty,
 			\]
 	where $\mathbf 0$ denotes the null application from $\mathbb R^n$ to $\mathbb R$.		
\end{enumerate}

\begin{remark}\label{rem:Ys-or-Ys-}
In the case where the integrator $C$ of the Lebesgue--Stieltjes integral is a continuous process, we can choose between the integrands 
$$\big(f(t,Y_{t},Z_t,U_t(\cdot))\big)_{t\in\llbracket 0,T\rrbracket} 
	\ \ \text{ and } \ \ 
  \big(f(t,Y_{t-},Z_t,U_t(\cdot))\big)_{t\in\llbracket 0,T\rrbracket},$$ 
and we still obtain the same solution, as they coincide outside of a $\ud C \otimes \ud\Pm-$null set. 
However, in the case where the integrator $C$ is \cadlag, the corresponding solutions may differ. 
In the formulation of the problem we have chosen the first one, while the \emph{a-priori} estimates can readily be adapted to the second case as well. 
However, in order to obtain the unique solution in the second case we will need an additional property to hold for the integrator $C$; see Condition \ref{Hfive} in Subsection \ref{subsec:AlternEstim}. 
Moreover, in Subsection \ref{subsec:Compare} we will see that, in special cases, the conditions for existence and uniqueness of a solution in these two cases can differ significantly.   
\end{remark}

In classical results on BSDEs, the pair $(\xi,f)$ is called \emph{standard data}. 
In our case, we generalize the last term and say that the sextuple $(X,\Fil, T,\xi,f,C)$ is the \emph{standard data under $\hat{\beta}$}, whenever its elements satisfy Assumptions \ref{Fone}--\ref{Ffive} for this specific $\hat{\beta}$. 

\begin{definition}\label{solutiondef}
A \emph{solution of the {\rm BSDE} \eqref{BSDEintegralform} with standard data $(X,\mathbb G, T,\xi,f,C)$} under $\hat\beta>0$ is a quadruple of processes 
\begin{align*}
(Y,Z,U,N)&\in 		\Htwos{}	\times	\Htwos{,\circ} \times
					\Htwos{,\natural} 	\times	\mathcal{H}^{2,\perp}_{\beta}
  \textrm{ or } (Y,Z,U,N)\in 		
					\mathcal{S}^{2}_{\beta}	\times	\Htwos{,\circ} \times
					\Htwos{,\natural} 	\times	\mathcal{H}^{2,\perp}_{\beta},
\end{align*}
for some $\beta\leq \hat{\beta}$ such that, $\Pm-a.s.$, for any $t\in\llbracket0,T\rrbracket$, 
	\[Y_t = \xi 	+ \int_t^T f(s,Y_s,Z_s,U_s(\cdot)) \dC{s} 
				- \int_t^T \ Z_s \ud X^\circ_s 
				-\int_t^T\int_{\R{n}}U_s(x) \widetilde{\mu}^\natural(\ud s,\ud x)
				- \int_t^T \ud N_s.
\numberthis\label{BSDE}
\]
\end{definition}

\vspace{0.5em}
\begin{remark}
We emphasize that in \eqref{BSDE}, the stochastic integrals are well defined since $(Z,U,N)\in\Htwos{,\circ} \times\Htwos{,\natural}\times\mathcal{H}^{2,\perp}_{\beta}$. Let us verify that the integral
$$ \int_0^\cdot f(s,Y_s,Z_s,U_s(\cdot)) \dC{s}$$
is also well--defined. First of all, we know by definition that for any $(y,z,u)\in\R{d}\times\R{d\times m}\times\mathfrak H$, there exists a $\ud C\otimes\ud\mathbb P-$null set $\mathcal N^{y,z,u}$ such that for any $(t,\omega)\notin\mathcal N^{y,z,u}$
$$f(t,\omega,y,z,u_t(\omega;\cdot))\ \text{is well defined and }u_t(\omega;\cdot)\in\mathfrak H_{t,\omega}.$$
Moreover, by Lemma \ref{lem:polish}, we know also that for some $\ud C\otimes\ud\mathbb P-$null set $\widetilde {\mathcal N}$, we have for every $(t,\omega)\notin\widetilde{\mathcal N}$, that $\mathfrak H_{t,\omega}$ is Polish for the norm $\tnorm{\cdot}_t(\omega)$, so that it admits a countable dense subset which we denote by $H_{t,\omega}$. Let us then define
$$H:=\left\{u\in\mathfrak H,\ u_t(\omega;\cdot)\in H_{t,\omega},\,�\forall(t,\omega)\notin\widetilde{\mathcal N}\right\},\ \mathcal N:=\bigcup\left\{\mathcal N^{y,z,u},\ (y,z,u)\in\mathbb Q^d\times\mathbb Q^{d\times m}\times H \right\},$$
where $\mathbb Q$ and $\mathbb Q^{d\times m}$ are the subsets of $\mathbb R$ and $\mathbb R^{d\times m}$ with rational components.

\vspace{0.5em}
Then, since $H$ is countable, $\mathcal N$ is still a $\ud C\otimes\ud\mathbb P-$null set. 
Then, it suffices to use \ref{Fthree} to realize that for any $(t,\omega)\notin\mathcal N\cup\widetilde{\mathcal N}$, $f$ is continuous in $(y,z,u)$, and conclude that we can actually define $f(t,\omega,y,z,u_t(\omega;\cdot))$ outside a universal $\ud C\otimes\ud\mathbb P-$null set. 
This implies in particular that for any $(Y,Z,U)\in\mathbb H^2_\beta\times \Htwos{,\circ} \times\Htwos{,\natural}$
$$f(t,\omega, Y_t(\omega),Z_t(\omega),U_t(\omega;\cdot)) \ \text{is defined for $\ud C\otimes\ud\mathbb P-a.e.$ $(t,\omega)\in\llbracket0,T\rrbracket\times\Omega$}.$$
Finally, it suffices to use \ref{Fthree} and \ref{Ffive} to conclude that 
$$\int_0^T|f(t,\omega, Y_t(\omega),Z_t(\omega),U_t(\omega;\cdot))|\ud C_t(\omega) \ \text{is finite for $\ud C\otimes\ud\mathbb P-a.e.$ $(t,\omega)\in\llbracket0,T\rrbracket\times\Omega$}.$$
\end{remark}

%% file: Existence_and_uniqueness_statement.tex
We devote this subsection to the statement of our main theorem. 
Before that, we need some preliminary results of a purely analytical nature, whose proofs are relegated to Appendix \ref{AppendixB}.

\begin{lemma}\label{lem:const}
Fix $\beta, \Psi>0$ and consider the set $\mathcal C_{\beta}:=\{(\gamma,\delta)\in(0,\beta]^2,\ \gamma < \delta\}$. 
We define the following quantity
\begin{align*}
\Pi^{\Psi}(\gamma,\delta):=\frac{9}{\delta}+(2+9\delta)\frac{\e{(\delta-\gamma)\Psi}}{\gamma(\delta-\gamma)}.
\end{align*}
Then, the infimum of \! $\Pi^\Psi$ is given by
\begin{align*}
M^\Psi(\beta)
	&:=\inf_{(\gamma,\delta)\in\mathcal C_{\beta}}\Pi^\Psi(\gamma,\delta) =\frac{9}{\beta} + \frac{\Psi^2(2+9\beta)}{\sqrt{\beta^2\phi^2+4}-2}\exp\Big(\frac{\beta\Psi+2-\sqrt{\beta^2\Psi^2+4}}{2}\Big),
\end{align*}
and is attained at the point $\big(\overline{\gamma}^\Psi(\beta),\beta\big)$ where
\begin{align*}
\overline{\gamma}^\Psi(\beta)
:=\frac{\beta\Psi - 2 + \sqrt{4+\beta^2\Psi^2}}{2\Psi}.
\end{align*}
In addition, if we define
\begin{align*}
\Pi_\star^\Psi(\gamma,\delta):=\frac8\gamma+\frac{9}{\delta}+{9\delta}\frac{\e{(\delta-\gamma)\Psi}}{\gamma(\delta-\gamma)},
\end{align*}
then the infimum of  $\Pi^\Psi_\star$ is given by $
M^\Psi_\star(\beta):=\inf_{(\gamma,\delta)\in\mathcal C_{\beta}}\Pi^\Psi_\star(\gamma,\delta)=\Pi^\Psi_\star(\overline\gamma^\Psi_\star(\beta),\beta)$, where $\overline{\gamma}_\star^\Psi(\beta)$ is the unique solution in $(\overline\gamma^\Psi_\star(\beta),\beta)$ of the equation with unknown $x$
\begin{align*}
8(\beta-x)^2 - 9 \beta\e{(\beta-x)\Psi}\big(\Psi x^2 - (\beta\Psi-2)x - \beta\big)=0.
\end{align*}
Moreover, it holds
\begin{align*}
\lim_{\beta\to\infty} M^\Psi(\beta)=\lim_{\beta\to\infty}M^\Psi_\star(\beta)= 9{\rm e}\Psi.
\end{align*}
\end{lemma}

\begin{theorem}\label{BSDEMainTheorem}
Let $(X,\mathbb G, T,\xi,f,C)$ be standard data under $\hat{\beta}$.
If $M^\Phi(\hat\beta)<\frac12$ $($resp. $M^\Phi_\star(\hat\beta)<\frac12)$, then there exists a unique quadruple $(Y,Z,U,N)$ which satisfies \eqref{BSDE} 
and with $\norm{(Y,Z,U,N)}_{\hat\beta}<\infty$ $($resp. $\norm{(Y,Z,U,N)}_{\star,\hat\beta}<\infty)$.
\end{theorem}

\begin{corollary}
Let $(X,\mathbb G, T,\xi,f,C)$ be standard data under $\hat{\beta}$, for $\hat{\beta}$ sufficiently large.
If for the constant $\Phi$ defined in \eqref{Ffouruniformjumps} holds
\begin{equation}\label{eq:condphi}
\Phi<\frac{1}{18\e{}},
\end{equation}
then the BSDE \eqref{BSDE} has a unique solution. 
\end{corollary}
\begin{proof}
Using the results of Lemma \ref{lem:const} and Theorem \ref{BSDEMainTheorem}, it is immediate that as soon as $\Phi<1/(18\text{e})$, then there always exists a unique solution of the BSDE for $\hat{\beta}$ large enough. 
\end{proof}

%% file: Comparison_literature.tex

\subsubsection{Some counterexamples}\label{subsec:some-counterexamples}

As mentioned already in the introduction, Confortola, Fuhrman and Jacod \cite[Section 4.3]{confortola2014backward} provided a counterexample to the existence or the uniqueness of the solution of a BSDE in case the integrator $C$ is not a continuous process.
We would like to shed some more light on their counterexample here, and discuss various situations in which a solution may or may not exist.

\vspace{.5em}
Let us first rewrite their counterexample using our notation. 
Let $T>0$, $\ell\in(0,T]$, $X$ be a piecewise constant process with potentially a single jump at time $\ell$, that is $\mathbb P(\Delta X_{\ell}\neq 0)=p\in(0,1)$ and $\mathbb P(\{\Delta X_t=0 \text{ for every }t\in(0,\infty)\})=1-p$.
Let $\Pi:=\{\omega\in\Omega,\; \Delta X_{\ell}(\omega)\neq 0\}$ and $\Pi^c$ be its complement.
Moreover, let $\mathbb G$ be the natural filtration of $X$, $C_\cdot=p\mathds{1}_{[{\ell},\infty)}(\cdot)$, and fix some generator $f:[0,T]\times\R{}\times\R{}\longrightarrow\R{}$. 
Given the structure of the filtration $\mathbb G$, the terminal condition $\xi$ can always have a decomposition of the form
$$
\xi(\omega) 
	 =: \xi^\Pi \mathds{1}_\Pi(\omega) + \xi^{\Pi^c} \mathds{1}_{\Pi^c}(\omega),
\;	 (\xi^\Pi, \xi^{\Pi^c}) \in\mathbb R\times\mathbb R.
$$
Then, \cite{confortola2014backward} considers the following BSDE
\begin{align}\label{BSDE:counterexample}
Y_t + \int_t^T U_s \mu^{X}(\ud s) = \xi + \int_{(t,T]}f(s,Y_{s-},U_s) \ud C_s,
\end{align} 
and shows that, if the generator has the form $f(t,y,u)=\frac{1}{p}(y + g(u))$ for a deterministic function $g:\R{}\longrightarrow\R{}$, then the BSDE can admit either infinitely many solutions or none.

\vspace{.5em}
Once again because of the structure of $\mathbb G$, one can show that the possible solutions for the BSDE \eqref{BSDE:counterexample} necessarily  have the following form
\begin{align*}
Y_t(\omega) 
	& = Y_0\mathds{1}_{[0,{\ell})}(t) + \xi^\Pi \mathds{1}_{[{\ell},\infty)}(t)\mathds{1}_\Pi(\omega) + \xi^{\Pi^c} \mathds{1}_{[{\ell},\infty)}(t)\mathds{1}_{\Pi^c}(\omega),\\
U_t(\omega) 
	& = \upsilon(t) + \upsilon^\Pi(t)\mathds{1}_\Pi(\omega)\mathds1_{({\ell},T]}(t) + \upsilon^{\Pi^c}(t)\mathds{1}_{\Pi^c}(\omega)\mathds1_{({\ell},T]}(t),
\end{align*}
for some $Y_0\in \R{}$ and some deterministic functions $\upsilon,\upsilon^\Pi,\upsilon^{\Pi^c}:[0,T]\longrightarrow\R{}$. 
However, only the value $\upsilon({\ell})$ is actually involved.
By \cite[Theorem II.3.26]{jacod2003limit} we have that $C$ is the compensator of $X$, \emph{i.e.} the process $\widetilde X_\cdot := X_\cdot - C_\cdot$ is a $\mathbb G-$martingale.
Now we can distinguish between the following cases.

\vspace{0.5em}
\begin{enumerate}[label=\textbf{(C\arabic*)}, leftmargin=*, itemindent=*] 
\item Consider the BSDE \eqref{BSDE:counterexample}.
	 	Then, there exists a solution if and only if there exists a fixed point, called $Y^\star_0$, for the equation
	 	\begin{align*}
	 		\xi^\Pi + p f({\ell}, x, \xi^\Pi-\xi^{\Pi^c}) = x.
	 	\end{align*}
	 	The pair $(Y_0^\star,\xi^\Pi- \xi^{\Pi^c})$ is a solution of \eqref{BSDE:counterexample}.
	 	The solution is unique if and only if the fixed point is unique.
	 	In case $f$ is globally Lipschitz with respect to its second argument, \emph{i.e.}
	 	\begin{align*}
	 		|f(t,y_1,u)-f(t,y_2,u)|^2 \le r |y_1-y_2|^2,
	 	\end{align*}
	 	then, a sufficient condition for the existence and uniqueness of the solution is $r\Delta C_{\ell}^2 <1.$
		
		\vspace{0.5em}
\item\label{Ctwo} Consider the following BSDE instead, where the stochastic integral is taken with respect to the compensated jump process
	\begin{align}\label{BSDE:counterexample-2}
		Y_t + \int_t^T U_s \widetilde{\mu}^{X}(\ud s) = \xi + \int_{(t,T]}f(s,Y_{s-},U_s) \ud C_s.
	\end{align}
	 	Then, there exists a solution if and only if there exists a fixed point, called $Y^\star_0$, for the equation
	 	\begin{align*}
	 		\xi^\Pi + p f({\ell}, x, \xi^\Pi-\xi^{\Pi^c}) + p(\xi^\Pi - \xi^{\Pi^c})= x.
	 	\end{align*}
	 	The pair $(Y_0^\star,\xi^\Pi- \xi^{\Pi^c})$ is a solution of \eqref{BSDE:counterexample-2}.
	 	The solution is unique if and only if the fixed point is unique.
	 	In case $f$ is globally Lipschitz with respect to the second argument as above, a sufficient condition for the existence and uniqueness of the solution is again $r\Delta C_{\ell}^2 <1.$
		
		\vspace{0.5em}
\item Consider now a BSDE similar to \eqref{BSDE:counterexample}, where the integrand of the Lebesgue--Stieltjes integral depends on $Y$ instead of $Y_-$, \textit{i.e.}
	\begin{align}\label{BSDE:counterexample-3}
		Y_t + \int_t^T U_s \mu^{X}(\ud s) = \xi + \int_{(t,T]}f(s,Y_s,U_s) \ud C_s.
	\end{align}
	 	Then, there exists a solution if and only if there exists a fixed point, called $\upsilon^\star({\ell})$, for the equation
	 	\begin{align*}
	 		\xi^\Pi - \xi^{\Pi^c} -  p f\big({\ell r},\xi^{\Pi^c}, x\big) + p f({\ell },\xi^\Pi, x) = x.
	 	\end{align*}
	 	The pair $(\xi^{\Pi^c} +  p f\big({\ell },\xi^{\Pi^c} + p \upsilon^\star({\ell }), \upsilon^\star({\ell})\big),\upsilon^\star({\ell }))$ is a solution of \eqref{BSDE:counterexample-3}.
	 	The solution is unique if and only if the fixed point is unique.
	 	In case $f$ is globally Lipschitz with respect to its third argument, \emph{i.e.}
	 	\begin{align*}
	 		|f(t,y,u_1)-f(t,y,u_2)|^2 \le \theta^\sharp |u_1-u_2|^2,
	 	\end{align*}
	 	then, a sufficient condition for the existence and uniqueness of the solution is $4\theta^\sharp\Delta C_{\ell}^2 <1.$
	 	This condition is not necessary: let $f'(t,y,u)=\frac{1}{p} (g(y) + u)$, where $g$ is a deterministic function, then it holds that $\theta^\sharp\Delta C_{\ell }^2=1$; however, \eqref{BSDE:counterexample-3} admits a unique solution, which is given by the pair 
	 	$$\big(\xi^\Pi + g(\xi^{\Pi}), \xi^\Pi - \xi^{\Pi^c} + g(\xi^\Pi) - g(\xi^{\Pi^c})\big).$$

\vspace{0.5em}
\item Finally, consider a BSDE similar to \eqref{BSDE:counterexample-3} where the stochastic integral is taken with respect to the compensated jump process
	\begin{align}\label{BSDE:counterexample-4}
		Y_t + \int_t^T U_s \widetilde{\mu}^{X}(\ud s) = \xi + \int_{(t,T]}f(s,Y_s,U_s) \ud C_s.
	\end{align}
	 	Then, there exists a solution if and only if there exists a fixed point, called $\upsilon^\star({\ell })$, for the equation
	 	\begin{align*}
	 		\xi^\Pi - \xi^{\Pi^c} -  p f\big({\ell },\xi^{\Pi^c}, x\big) + p f({\ell },\xi^\Pi, x) = x.
	 	\end{align*}
	 	The pair $(\xi^{\Pi^c} +  p f\big({\ell },\xi^{\Pi^c}, \upsilon^\star({\ell })\big),\upsilon^\star({\ell }))$ is a solution of \eqref{BSDE:counterexample-4}.
	 	The solution is unique if and only if the fixed point is unique.
	 	In case $f$ is globally Lipschitz with respect to its third argument as above, a sufficient condition for the existence and uniqueness of the solution is again $4\theta^\sharp\Delta C_{\ell }^2 <1.$
	 	Once again this condition in not necessary; indeed, for $f'$ as in (C3), $\theta^\sharp\Delta C_{\ell }^2=1$, while the unique solution of the BSDE \eqref{BSDE:counterexample-4} is the pair 
	 	$$\big((1-p) [\xi^{\Pi} + g(\xi^\Pi)] + p [\xi^{\Pi^c}  +g(\xi^{\Pi^c})], \xi^\Pi - \xi^{\Pi^c} + g(\xi^\Pi) - g(\xi^{\Pi^c})\big).$$
\end{enumerate}

\vspace{.5em}
Now, returning to the original counterexample of \cite{confortola2014backward}, we can observe that the sufficient condition $r \Delta C_{\ell }^2 < 1$ is violated there, which explains why wellposedness issues can arise.
However, an important observation here is that the structure of the generator plays a crucial role as well.
Indeed, if we consider the same BSDE with the following generator $f(t,y,u) = m (y+g(u))$ with $m\neq\frac1p$, then the BSDE admits a unique solution. 

\vspace{.5em}
Let us finally argue why condition \eqref{eq:condphi} rules out this counterexample from our setting. 
The generator $f(t,y,u) = \frac1p (y+g(u))$ needs to be Lipschitz so that it fits in our framework, and to satisfy \eqref{fstochLip}. 
Let us further assume that the function $g$ is also Lipschitz, say with associated constant $L^g.$
Then, using Young's Inequality, we can obtain 
\begin{align*}
|f(t,y,u) -f(t,y',u')|^2 \le \frac{1+\varepsilon}{p^2} |y-y'|^2 + \frac{1}{p^2}\bigg(1+\frac{(L^g)^2}{\varepsilon}\bigg)|u-u'|^2, \text{ for every }\varepsilon >0.
\end{align*} 
Before we proceed recall \eqref{eq:A} and \eqref{Ffouruniformjumps}, \emph{i.e.} $A_\cdot =\int_0^{\cdot} \alpha_s^2 \dC{s}$ and $\Delta A_t(\omega) \leq \Phi$, for $\ud C \otimes \ud \Pm-a.e.$ 
Therefore, we have that $\alpha^2_{\mathpzc r}=\max\Big\{ \sqrt{1+\varepsilon}/p,  \Big(1+\frac{(L^g)^2}{\varepsilon}\Big)/p^2\Big\}$, and 
\begin{align*}
\alpha^2_{\mathpzc r} \Delta C_{\mathpzc r} = \max\bigg\{ \sqrt{1+\varepsilon},  \frac1p+\frac{(L^g)^2}{p\varepsilon}\bigg\}\ge \sqrt{1+\varepsilon}>\frac{1}{18{\rm e}} \text{ for every }\varepsilon>0.
\end{align*}

\begin{remark}
Coming back to Remark \ref{rem:Ys-or-Ys-}, we observe that the dependence of the integrand on $Y$ or $Y_-$ is not always that innocuous.
Indeed, the same BSDE might have a solution in the one formulation but not in the other.
Observe furthermore that in the first situation the Lipschitz constant $r$ appears in the condition for the existence and uniqueness of a solution, while in the second case the Lipschitz constant $\theta^\sharp$ appears.
As stated already before, in our framework we can treat both cases simultaneously, hence naturally both Lipschitz constants appear in our condition, through the definition of $\alpha^2$ as the maximum of all the Lipschitz constants.
\end{remark}

\subsubsection{Related literature}

Let us now compare our work with the papers by Bandini {\rm \cite{bandini2015existence}} and Cohen and Elliott {\rm\cite{cohen2012existence}} who also consider BSDEs in stochastically discontinuous filtrations.
The setting in \cite{cohen2012existence} is rather different from ours. 
Indeed, in our case a driving martingale $X$ is given right from the start, and as a consequence the process $C$ with respect to which the generator $f$ is integrated is linked to the predictable bracket of $X$. 
However, the authors of \cite{cohen2012existence} do not choose any $X$ from the start, but consider instead a general martingale representation theorem involving countably many orthogonal martingales, which only requires the space of square integrable random variables on $(\Omega,\mathcal F,\mathbb P)$ to be separable to hold. 
Furthermore, their process $C$ can, unlike our case, be chosen arbitrary (in the sense that it does not have to be related to the driving martingales), but with the restriction that it has to be deterministic. 
Moreover, it has to assign positive measure on every interval, see Definition 5.1 therein, hence $C$ cannot be piecewise constant; the latter would naturally arise from a discrete-time martingale with independent increments, which is exactly the situation one encounters when devising numerical schemes for BSDEs.
Therefore, their setting cannot be embedded into our framework, and \textit{vice versa}.

\vspace{0.5em}
On the other hand, in \cite{bandini2015existence}, the author considers a BSDE driven by a pure--jump martingale without an orthogonal component, which is a special case of \eqref{BSDEintegralform}.
The martingale in this setting should actually have jumps of finite activity, hence many of the interesting models for applications in mathematical finance, such as the generalized hyperbolic, CGMY, and Meixner processes, are excluded. Such a restriction is not present at all in our framework.
Otherwise, the assumptions and the conclusion in \cite{bandini2015existence} are analogous to the present work.
A direct comparison is however not possible, \textit{i.e.} we cannot deduce the existence and uniqueness results in her work from our setting, since the assumptions are not exactly comparable. In particular, the integrability condition ${\rm (iii)}$ on page 3 in \cite{bandini2015existence} is not compatible with \ref{Ffive}.

\vspace{0.5em}
Let us also compare our result with the literature on BSDEs with random terminal time.
Royer \cite{royer2004bsdes}, for instance, considers a BSDE driven by Brownian motion, where the terminal time is a $[0,\infty)-$valued stopping time. 
Hence, her setting can be embedded in ours, by assuming the absence of jumps and of the orthogonal component, and further requiring that $C$ is a continuous process.
She shows existence and uniqueness of a solution under the assumptions that the generator is uniformly Lipschitz in $z$ and continuous in $y$, and the terminal condition is bounded. 
Moreover, she requires that either the generator is strictly monotone in $y$ and $f(t,0,0)$ is bounded (for all $t$) or that the generator is monotone in $y$ and $f(t,0,0)=0$ (for all $t$).
These conditions are not directly covered by our Assumptions \ref{Fone}--\ref{Ffive}, however if we consider her conditions and assume in addition that the generator is Lipschitz in $y$, then we can recover the existence and uniqueness result from our main theorem.
Let us point out that BSDEs with constant terminal time are related to semi--linear parabolic PDEs, while BSDEs with random terminal time are associated to semi--linear elliptic PDEs.

\vspace{.5em}
We would also like to comment briefly on the choice of the norms we consider here. 
They are mostly inspired from the ones defined in the seminal work of El Karoui and Huang \cite{elkaroui1997general}, and are equivalent to the usual norms found in the literature when the process $A$ and time $T$ are both bounded. 
Bandini \cite{bandini2015existence} uses different spaces, where the norm is defined using the Dol\'eans--Dade stochastic exponential instead of the natural exponential. 
In our setting where $A$ is allowed to be unbounded, we can only say that our norm dominates hers. 
This means that we require stronger integrability conditions, but as a result we will also obtain a solution of the BSDE with stronger integrability properties.
In any case, our method could be adapted to this choice of the norm, albeit with modified computations in our estimates. 
We refer the reader to Remark \ref{remark:norms} below for a more detailed discussion about the definition of the norms.

\vspace{0.5em}
Let us conclude this section by commenting on the condition \eqref{eq:condphi}. 
We start with the observation that the analysis of the counterexample of Confortola \emph{et al.} \cite{confortola2014backward} made in Sub--sub--section \ref{subsec:some-counterexamples} does not allow for a general statement of wellposedness of the BSDE when $\Phi\ge 1$.
In this light, the result of Cohen and Elliott \cite[Theorem 6.1]{cohen2012existence}, which implies that the condition $\Phi<1$ ensures the wellposedness of the BSDE, lies in the optimal range for $\Phi$.  
Analogously in the case of Bandini \cite{bandini2015existence}, once her results are translated using the Lipschitz assumption in \ref{Fthree}, $\Phi<1$ also ensures the wellposedness of the BSDE.
On the contrary, condition \eqref{eq:condphi} which reads as $\Phi<1/(18\e{})$, may seem much more restrictive. 
The first immediate remark we can make is that the stochasticity of the integrator $C$ considerably deteriorates the condition on $\Phi.$ 
In \cite{cohen2012existence} the integrator is deterministic, while in \cite{bandini2015existence} and in our case the integrator is stochastic.
However, we would like to remind the reader, that, as explained above, the level of generality we are working with is substantially higher than in these two references. 
We also want to emphasize the fact that our condition is clearly not the sharpest one possible, but we believe it is the sharpest that can be obtained using our method of proof. 
The main possibilities for improvement are, in our view, twofold:
\begin{enumerate}[label=\raisebox{0.08em}{\textbullet}]
\item First of all, in specific situations (e.g. $T$ bounded, $f$ Lipschitz, less general driving processes, \dots) one should be able to improve the {\it a priori} estimates of Lemma \ref{lemmaest} by refining several of the inequalities. This exactly what we will do in Section \ref{subsec:AlternEstim} below, by using an approach reminiscent of the one usually used in the BSDE literature.

\item Second, as highlighted in Remark \ref{remark:norms}, we actually have a degree of freedom in choosing the norms we are interested in. In this paper, we used exponentials, while Bandini \cite{bandini2015existence} used stochastic exponentials, but other choices, leading to potentially better estimates, could also be considered.
\end{enumerate}
We leave this interesting problem of finding the optimal $\Phi$ open for future research.

%% file: A_priori_estimates.tex
The method of proof we will use follows and extends the one of El Karoui and Huang \cite{elkaroui1997general}. 
In \cite{elkaroui1997general}, as well as in Pardoux and Peng \cite{pardoux1990adapted}, the result is obtained using fixed--point arguments and the so-called {\it a priori} estimates.  
However, we would like to underline that the proof of such estimates in our case is significantly harder, due to the fact that the process $C$ is not necessarily continuous.

\vspace{0.5em}
The following result can be seen as the {\it a priori} estimates for a BSDE whose generator does not depend on the solution. In order to keep notation as simple as possible, as well as to make the link with the data of the problem we consider clearer, we will reuse part of the notation of \ref{Fone}--\ref{Ffive}, namely $\xi, T, f,C,\alpha$ and $A$, only for the next two lemmata. 

\begin{lemma}\label{lemmaest}
	Let $y$ be a $d-$dimensional $\Fil-$semimartingale of the form 
	\begin{equation}\label{BSDE2}
		y_t=\xi + \int_t^T f_s \dC{s} - \int_t^T \ud\eta_s,
	\end{equation}
	where $T$ is a $\Fil-$stopping time, $\xi\in\mathbb{L}^2\pair{\G_{T};\R{d}}$, $f$ is a $d-$dimensional optional process,	$C$ is an increasing, predictable and \cadlag~ process and $\eta\in\Hcal^2$.
	
	\vspace{0.5em}
Let $A:=\alpha^2 \cdot C$ for some predictable process $\alpha$. Assume that there exists $\Phi>0$ such that property \eqref{Ffouruniformjumps} holds for $A$.
	Suppose there exists $\beta\in\Rp$ such that
	\[
		\E\left[ \e{\beta A_T} \abs{\xi}^2 \right]+ \E\left[\int_0^T\e{\beta A_t} \frac{\abs{f_t}^2}{\alpha_t^2} \dC{t}\right]<\infty.
	\]
	Then we have for any $(\gamma,\delta)\in(0,\beta]^2$, with $\gamma \ne \delta$,
$$
		\norm{\alpha y}_{\mathbb{H}^2_{\delta}}^2 \leq 
		\frac{2\e{\delta \Phi}}{\delta} \xinorm{\delta} + 
		2\Lambda^{\gamma,\delta,\Phi}\fnorm{\gamma\vee\delta},\ 
		\norm{y}_{\mathcal{S}^2_\delta}^2 
		\leq 8 \xinorm{\delta} + 
		\frac{8}{\gamma}\fnorm{\gamma\vee\delta},
$$
\vspace{-0.5em}
$$	
\norm{\eta}_{\mathcal{H}^2_{\delta}}^2  \leq 
	9\left(1 + \e{\delta\Phi}\right)	\xinorm{\delta}+
	9\left(\frac{1}{\gamma\vee\delta}+\delta \Lambda^{\gamma,\delta,\Phi}	\right)			\fnorm{\gamma\vee\delta},
$$
where we have defined
\[	
	\Lambda^{\gamma,\delta,\Phi}:=\frac{1\vee\e{(\delta-\gamma)\Phi}}{\gamma\abs{\delta-\gamma}}.
\]

\vspace{0.5em}
As a consequence, we have	\begin{align}\label{sumnormest}
		\norm{\alpha y}_{\mathbb{H}^2_{\delta}}^2 +
		\norm{\eta}_{\Hcal^2_{\delta}}^2 &\leq
		\widetilde{\Pi}^{\delta,\Phi} \xinorm{\delta} +
		\Pi^{\Phi}(\gamma,\delta)\fnorm{\gamma\vee\delta},\\
		\label{sumnormstarest}
		\norm{y}_{\mathcal{S}^2_{\delta}}^2 +
		\norm{\eta}_{\Hcal^2_{\delta}}^2 &\leq
		\widetilde{\Pi}_{\star}^{\delta,\Phi} \xinorm{\delta} +
		\Pi_{\star}^{\Phi}(\gamma,\delta)\fnorm{\gamma\vee\delta},
	\end{align}
where
$$\widetilde{\Pi}^{\delta,\Phi}	:=9+\pair{9+\frac{2}{\delta}}\e{\delta\Phi} \ \text{ and } \ \widetilde{\Pi}_{\star}^{\delta,\Phi}:=17+9\e{\delta\Phi}.$$
\end{lemma}
\begin{proof}
Recall the identity
\begin{equation}\label{identity1}
	y_t=\xi + \int_t^Tf_s \dC{s} - \int_t^T\ud\eta_s = \E\left[ \left.\xi + \int_t^T f_s \dC{s} \right| \mathcal{G}_t \right],
\end{equation}
and introduce the anticipating function
\begin{equation}
	F(t)=\int_t^Tf_s \dC{s}.
\end{equation}

For $\gamma\in \mathbb{R}_+ $, we have by the Cauchy--Schwarz inequality,
\begin{align*}
	\abs{F(t)}^2 &\leq \int_t^T \e{-\gamma A_s} \dA{s} 
	\int_t^T\e{\gamma A_s} \frac{\abs{f_s}^2}{\alpha_s^2} \dC{s}
	\leq \int_{A_t}^{A_T} \e{-\gamma A_{L_s}} \ds \int_t^T\e{\gamma A_s}\frac{\abs{f_s}^2}{\alpha_s^2} \dC{s}\\
	&\leq \int_{A_t}^{A_T} \e{-\gamma s} \ds \int_t^T\e{\gamma A_s}\frac{\abs{f_s}^2}{\alpha_s^2} \dC{s}
	\leq \frac{1}{\gamma} \e{-\gamma A_t} \int_t^T\e{\gamma A_s}\frac{\abs{f_s}^2}{\alpha_s^2} \dC{s},
	\numberthis \label{ineq1}
\end{align*}
where for the third inequality we used Lemma \ref{GeneralisedInverseLemma}.\ref{GILemmaSeven}. For $t=0$, since we assumed that 
$$\E\left[\int_0^T\e{\beta A_t} \frac{\abs{f_t}^2}{\alpha_t^2} \dC{t}\right]<\infty,$$ we have that the following holds for $0<\gamma<\beta$ 
\[
\E\left[\abs{F(0)}^2\right]<\infty.
\]

For $\delta\in\Rp$ and by integrating \eqref{ineq1} w.r.t. $\e{\delta A_t} \ud A_t$ it follows
\begin{align*}
	\int_0^T\e{\delta A_t}\abs{F(t)}^2 \ud A_t &\stackrel{\eqref{ineq1}}{\leq} 
	\frac{1}{\gamma}\int_0^T\e{(\delta-\gamma) A_t}
	\int_t^T\e{\gamma A_s} \frac{\abs{f_s}^2}{\alpha_s^2}\dC{s} \dA t \\
	&\stackrel{\quad \quad}{=}\frac{1}{\gamma}\int_0^T\e{\gamma A_s}\frac{\abs{f_s}^2}{\alpha_s^2} 
	\int_0^{s-}\e{(\delta-\gamma) A_t} \dA t \dC{s}\\
	&\stackrel{\quad \quad}{\le}\frac{1}{\gamma}\int_0^T\e{\gamma A_s}\frac{\abs{f_s}^2}{\alpha_s^2} 
	\int_0^{s}\e{(\delta-\gamma) A_t} \dA t \dC{s}, \numberthis \label{ineq2}
\end{align*}
where we used Tonelli's Theorem in the equality. We can now distinguish between two cases:

\vspace{0.5em}
$\bullet$ For $\delta>\gamma$, we apply Corollary \ref{maincorollary} for $g(x)=\e{(\delta-\gamma)x}$, and inequality \eqref{ineq2} becomes
\begin{align*}
	\int_0^T\e{\delta A_t}\abs{F(t)}^2 \dA t 
	&\leq \frac{\e{(\delta-\gamma)\Phi}}{\gamma}\int_0^T\e{\gamma A_s}\frac{\abs{f_s}^2}{\alpha_s^2}
	\int_{A_0}^{A_s}\e{(\delta-\gamma)t} \dt \dC{s}\leq \frac{\e{(\delta-\gamma)\Phi}}{\gamma(\delta-\gamma)}\int_0^T\e{\delta A_s}\frac{\abs{f_s}^2}{\alpha^2_s}\dC{s}, \numberthis \label{ineq3}
\end{align*}
which is integrable if $\delta\leq \beta$.

\vspace{0.5em}
$\bullet$ For $\delta<\gamma$, inequality \eqref{ineq2} can be rewritten as follows
\begin{align*}
	\int_0^T\e{\delta A_t}\abs{F(t)}^2 \dA t 
	&\stackrel{\quad \quad \quad \quad \ \ \ }{\leq} \frac{1}{\gamma}\int_0^T\e{\gamma A_s}\frac{\abs{f_s}^2}{\alpha_s^2} 
	\int_{A_0}^{A_s}\e{(\delta-\gamma) A_{L_t}} \dt\dC{s}\\
	&\stackrel{\textrm{Lem. }\ref{GeneralisedInverseLemma}.\ref{GILemmaSeven}}{\leq} 
	\frac{1}{\gamma }  
	\int_0^T\e{\gamma A_s}\frac{\abs{f_s}^2}{\alpha_s^2}\int_{A_0}^{A_s}\e{(\delta-\gamma) t} \dt \dC{s}\\
	&\stackrel{\quad \quad \quad \quad \ \ \ }{\leq}\frac{1}{\gamma \abs{\delta-\gamma}}  
	\int_0^T\e{\gamma A_s}\frac{\abs{f_s}^2}{\alpha_s^2} 
	\left(\e{(\delta-\gamma) A_0}-\e{(\delta-\gamma) A_s}\right) \dC{s}\\
	&\stackrel{\quad \quad \quad \quad \ \ \ }{\leq} \frac{1}{\gamma \abs{\delta-\gamma}}\int_0^T\e{\gamma A_s}\frac{\abs{f_s}^2}{\alpha_s^2}\dC{s}, 
	\numberthis \label{ineq4}
\end{align*}
which is integrable if $\gamma\leq \beta$. To sum up, for $\gamma,\delta\in(0,\beta]$, $\gamma\neq \delta$, we have
\begin{equation}\label{ineqA}
	\Expect{\int_0^T\e{\delta A_t}\abs{F(t)}^2 \dA t} \leq \Lambda^{\gamma,\delta,\Phi} \norm{\frac{f}{\alpha}}^2_{\mathbb{H}^2_{\gamma\vee\delta}}.
\end{equation}
For the estimate of $\norm{\alpha y}_{\mathbb{H}_{\delta}^2}$ we first use the fact that
\begin{align*}
	\norm{\alpha y}^2_{\mathbb{H}_{\delta}^2}=\Expect{\int_0^T\e{\delta A_t} \abs{y_t}^2 \dA t}&\leq 
	2 \Expect{\int_0^T\mathbb E\left[\left. \e{\delta A_t}\abs{\xi}^2 + \e{\delta A_t}\abs{F(t)}^2 \right|\mathcal G_t\right]\dA t }\\
	&= 
	2 \Expect{ \int_0^{\infty}\mathbb E\left[\left.\e{\delta A_t}\abs{\xi}^2 + \e{\delta A_t}\abs{F(t)}^2 \right|\mathcal G_t\right] \dA t^T }\\
	&\stackrel{\text{Cor. \ref{OptFubini}}}{=} 
	2 \Expect{ \int_0^{\infty}\e{\delta A_t}\abs{\xi}^2 + \e{\delta A_t}\abs{F(t)}^2 \dA t^T }\\
	&= 
	2 \Expect{ \int_0^T \e{\delta A_t}\abs{\xi}^2 + \e{\delta A_t}\abs{F(t)}^2 \dA t }\\
	&\overset{\textrm{Cor. }\ref{maincorollary}}{\underset{\eqref{ineqA}}{\leq} }
	\frac{2\e{\delta \Phi}}{\delta} \norm{\xi}^2_{\mathbb{L}^2_{\delta}} + 
	2 \Lambda^{\gamma,\delta,\Phi} \norm{\frac{f}{\alpha}}^2_{\mathbb{H}^2_{\gamma\vee\delta}}.
	\numberthis \label{ineqB}
\end{align*}
In the second equality we have used that the processes $\abs{\xi}^2\mathds{1}_{\Omega\times[0,\infty]}(\cdot)$ and $\abs{F(\cdot)}^2$ are uniformly integrable, hence their optional projections are well defined. 
Indeed, using \eqref{ineq1} and remembering that $\mathbb E[\abs{F(0)}^2]<\infty$, we can conclude the uniform integrability of $\abs{F(\cdot)}^2$. 
Then, by \cite[Theorem 5.4]{he1992semimartingale} it holds that
\begin{align*}
^o\!\left( \e{\delta A_{\cdot}}\abs{\xi}^2 + \e{\delta A_{\cdot}}\abs{F(\cdot)}^2\right)_t
&= 	\e{\delta A_t}~ ^o\!\left(\abs{\xi}^2 \mathds{1}_{\Omega\times[0,\infty]}(\cdot)\right)_t+ 
	\e{\delta A_t}~ ^o\!\left(\abs{F(\cdot)}^2\right)_t \\
&=	\e{\delta A_t} \mathbb E\left[\left.\abs{\xi}^2 \right|\mathcal G_t\right]+ 
	\e{\delta A_t} \mathbb E\left[\left.\abs{F(t)}^2 \right|\mathcal G_t\right]=	\mathbb E\left[\left.\e{\delta A_t}\abs{\xi}^2+ 
	\e{\delta A_t} \abs{F(t)}^2 \right|\mathcal G_t\right],
\end{align*}
which justifies the use of Corollary \ref{OptFubini}. For the estimate of $\norm{y}_{\mathcal{S}_{\delta}^2}$ we have
\begin{align*}
	\norm{y}_{\mathcal{S}_{\delta}^2} 
	&= \Expect{\sup_{0\leq t\leq T}\left(\e{\frac{\delta}{2}A_t} \abs{y_t} \right)^2}\\
	&\leq 2\Expect{\sup_{0\leq t\leq T} 
	\mathbb E\left[\left.\sqrt{\e{\delta A_t} \abs{\xi}^2 +
	\e{\delta A_t}\abs{F(t)}^2}\right|\mathcal G_t\right]^2} \\
	&\leq2\Expect{\sup_{0\leq t\leq T} 
	\mathbb E\left[\left.\sqrt{\e{\delta A_T} \abs{\xi}^2 +
	\frac{1}{\gamma}\e{(\delta-\gamma) A_t} \int_t^T \e{\gamma A_s}\frac{\abs{f_s}^2}{\alpha_s^2} \dC{t}} \right|\mathcal G_t\right]^2}\\
	&\leq 2\Expect{\sup_{0\leq t\leq T} 
	\mathbb E\left[\left.\sqrt{\e{\delta A_T} \abs{\xi}^2 +
	\frac{1}{\gamma} \int_0^T \e{(\gamma+(\delta-\gamma)^+ )A_s}\frac{\abs{f_s}^2}{\alpha_s^2} \dC{s}} \right|\mathcal G_t\right]^2}\\
	&\leq 8\Expect{
	\e{\delta A_T} \abs{\xi}^2 +
	\frac{1}{\gamma} \int_0^T \e{(\gamma\vee\delta) A_s}\frac{\abs{f_s}^2}{\alpha_s^2} \dC{s}}
	\leq 8\xinorm{\delta} + 
	\frac{8}{\gamma} \fnorm{\gamma\vee\delta}
	\numberthis \label{ineqC}
\end{align*}
for $\gamma\vee\delta\leq \beta$, where in the second and third inequalities we used the inequality $a+b\leq \sqrt{2(a^2+b^2)}$ and \eqref{ineq1} respectively.

\vspace{0.5em}
What remains is to control $\norm{\eta}_{\mathcal{H}^2_{\delta}}$. We remind once more the reader that $
	\int_t^T \ud\eta_s =\xi - y_t + F(t),
$
hence
\begin{equation}\label{ineq5}
	\mathbb E\left[\left.\abs{\xi-y_t-F(t)}^2\right|\mathcal G_t\right] = 
	\mathbb E\left[\left. \int_t^T\ud{\rm Tr}[\,\pqc{\eta}_s]\right|\mathcal G_t\right].
\end{equation}
In addition, we have
\begin{align*}
	\int_0^T\e{\delta A_s} \ud\tr\pqc{\eta}_s 
	&\stackrel{\quad \quad \quad \quad \ }{=}
	\int_0^T\int_{A_0}^{A_s} \delta \e{\delta t} \dt \ud{\rm Tr}[\,\pqc{\eta}_s] +{\rm Tr}[\,\pqc{\eta}_T] \\
	&\stackrel{\textrm{Lem. }\ref{GeneralisedInverseLemma}.\ref{GILemmaSeven}}{\leq} 
	\int_0^T\int_{A_0}^{A_s} \delta \e{\delta A_{L_t}} \dt \ud{\rm Tr}[\,\pqc{\eta}_s] + {\rm Tr}[\,\pqc{\eta}_T]\\  
	&\stackrel{\quad \quad \quad \quad \ }{=}
	\delta \int_0^T\int_{0}^{s} \e{\delta A_{t}} \dA t\ud{\rm Tr}[\,\pqc{\eta}_s] + {\rm Tr}[\,\pqc{\eta}_T]\\
	&\stackrel{\quad \quad \quad \quad \ }{\leq}
	\delta \int_{0}^{T} \e{\delta A_{t}} \int_t^T  \ud{\rm Tr}[\,\pqc{\eta}_s]\dA t +{\rm Tr}[\,\pqc{\eta}_T]  ,
\end{align*}
so that
\begin{align}\label{ineqD}
	\norm{\eta}_{\mathcal{H}^2_{\delta}}
	&\leq \delta \Expect{\int_{0}^{T} \e{\delta A_{t}} \int_t^T  \ud{\rm Tr}[\,\pqc{\eta}_s] \dA t} + 
	\Expect{{\rm Tr}[\,\pqc{\eta}_T]}.
\end{align}

For the first summand on the right-hand-side of \eqref{ineqD}, we have
\begin{align*}
	\Expect{\int_0^T \e{\delta A_t} \int_t^T\ud{\rm Tr}[\,\pqc{\eta}_s] \dA t}  
	&\stackrel{\text{Cor. \ref{OptFubini}}}{=}
	\Expect{\int_0^T \e{\delta A_t} \mathbb E\left[\left. \int_t^T \ud{\rm Tr}[\,\pqc{\eta}_s]\right|\mathcal G_t\right] \dA t}\\
	&\stackrel{(\ref{ineq5})}{=}
	\Expect{\int_0^T  \e{\delta A_t} \mathbb E\left[\left. \abs{\xi-y_t+F(t)}^2 \right|\mathcal G_t\right] \dA t}\\
	&\leq
	3 \Expect{\int_0^T \e{\delta A_t} \mathbb E\left[\left.\abs{\xi}^2+ \abs{y_t}^2+\abs{F(t)}^2\right|\mathcal G_t\right] \dA t}\\
	&\stackrel{\eqref{identity1}}{\leq}
			3 \Expect{ \int_0^T \e{\delta A_t} \abs{\xi}^2\ \dA t} + 3 \Expect{ \int_0^T \e{\delta A_t} \abs{F(t)}^2\ \dA t }
			\\
& \hspace{1.9em} +6 \Expect{ \int_0^T \e{\delta A_t}\mathbb E\left[\left. \abs{\xi}^2 + \abs{F(t)}^2 \right|\mathcal G_t\right] \dA t}\\
	&\stackrel{\text{Cor. \ref{OptFubini}}}{=}
	9 \Expect{ \int_0^T \e{\delta A_t} \abs{\xi}^2\ \dA t} +
	9 \Expect{ \int_0^T \e{\delta A_t} \abs{F(t)}^2\ \dA t }\\
	&\overset{\textrm{Cor. } \ref{maincorollary}}{\underset{\eqref{ineqA}}{\leq}}	
	\frac{9\e{\delta \Phi}}{\delta} \xinorm{\delta}
	+9 \Lambda^{\gamma,\delta,\Phi} \fnorm{\gamma\vee\delta}.
\end{align*}

We now need an estimate for $\Expect{\int_0^T\ud{\rm Tr}[\,\pqc{\eta}_s]}$, \emph{i.e.} the second summand of \eqref{ineqD}, which is given by 
\begin{align*}
	\Expect{{\rm Tr}[\,\pqc{\eta}_T]} &= \Expect{\abs{\xi-y_0+F(0)}^2 } 
	\leq 3 \Expect{\abs{\xi}^2 + \abs{y_0}^2+\abs{F(0)}^2 }\\
	&\overset{\eqref{identity1}}{\leq} 9 \Expect{\abs{\xi}^2} + 9\Expect{ \abs{F(0)}^2} 
	\overset{\eqref{ineq1}}{\leq} 9 \xinorm{\delta} + \frac{9}{\gamma\vee\delta}\fnorm{\gamma\vee\delta},
\end{align*}
where we used the fact that $\Expect{\abs{y_0}^2}\leq 2 \Expect{\abs{\xi}^2+\abs{F(0)}^2}$. 

\vspace{0.5em}
Then \eqref{ineqD} yields
\begin{equation}\label{ineqE}
		\norm{\eta}_{\mathcal{H}^2_{\delta}}^2  \leq 
	9\left(1 + \e{\delta\Phi}\right)									\xinorm{\delta} +
	9\left(\frac{1}{{\gamma\vee\delta}}+\delta \Lambda^{\gamma,\delta,\Phi}	\right)			\fnorm{\gamma\vee\delta}. \qedhere
\end{equation}
\end{proof}

\begin{remark}\label{remark:norms}
An alternative framework can be provided if we define the norms in Subsection \ref{subsec:SuitableSpaces} using another positive and increasing function $h$ instead of the exponential function. 
In order to obtain the required \emph{a priori} estimates, we need to assume that $h$ is sub--multiplicative\footnote{In the proof of \cite[Proposition 25.4]{sato1999levy} we can find a convenient tool for constructing sub--multiplicative functions.} 
and that it shares some common properties with the exponential function.
The following provides the \emph{a priori} estimates of the semi--martingale decomposition \eqref{BSDE2} in the case $h:\mathbb R\to [1,\infty)$, with $h(x)=(1+x)^\zeta$, for $\zeta\ge1,$ with the additional assumption that the process $A$ defined in \ref{Ffour} is $\mathbb P-$a.s. bounded by $\Psi.$
It holds for $\frac 1 \zeta < \gamma<\delta<\hat \beta$
\begin{align*}
\Vert\alpha y\Vert^2_{\mathbb{H}_{\delta}^2} + \Vert\eta\Vert_{\mathbb{H}_{\delta}^2}^2 
	&\le \left(2h(\Psi) h^\delta(\Phi) + 9 + \frac{9[h(\Psi)]^{1-\frac 1 \zeta} [h(\Phi)]^{\delta - \frac 1 \zeta}}{\delta-\frac 1 \zeta +1}	 \right)
		\Vert \xi\Vert^2_{\mathbb L^2_\delta} \\
	&\hspace{0.5em} + \left(		\frac{2 [h(\Psi)]^{1+\frac 1 \zeta}[h(\Phi)]^{\delta -\gamma +\frac 1 \zeta}}{\delta- \gamma + \frac 1 \zeta +1}
				+ 	\frac{9h(\Psi)[h(\Phi)]^{\delta-\gamma}}{\delta-\gamma+1}   + \frac{9}{\delta\zeta-1}	\right)		
		\left\Vert\frac{f}{\alpha}\right\Vert_{\mathbb H^2_{\delta}},\\
\Vert y\Vert^2_{\mathcal {S}_{\delta}^2} + \Vert\eta\Vert_{\mathbb{H}_{\delta}^2}^2 
	&\le \left(8 +2h(\Psi) h^\delta(\Phi)\right)
		\Vert \xi\Vert^2_{\mathbb H^2_\delta} + \left( \frac{[h(\Psi)]^{\frac 1 \zeta}}{\gamma\zeta - 1}  + \frac{2 [h(\Psi)]^{1+\frac 1 \zeta}[h(\Phi)]^{\delta -\gamma + \frac 1 \zeta }}{\delta - \gamma + \frac 1 \zeta +1}\right)
		\left\Vert \frac{f}{\alpha}\right\Vert^2_{\mathbb H^2_{\delta}}.
\end{align*}

\end{remark}

Let us also provide the following pathwise estimates.

\begin{lemma}\label{PathwiseAPriori}
Let $\xi, T, C,\alpha$, $A$ and $\Phi$ as in Lemma \ref{lemmaest}.
Assume that the $d-$dimensional $\mathbb G-$semimartingales $y_t^1$ and $y_t^2$ can be decomposed as follows
\begin{equation}\label{identity-yi-path}
		y_t^i=\xi + \int_t^T f_s^i\,\textup{d}C_{s} - \int_{t}^{T}\,\textup{d}\eta_s^i, \text{ for }i=1,2
\end{equation}
where $f^1,$ $f^2$ are $d-$dimensional $\mathbb G-$optional processes such that 
\begin{equation*}
	\mathbb E\left[\int_0^T\textup{e}^{\beta A_t} \frac{\vert f_t^i\vert^2}{\alpha_t^2}\,\textup{d}C_{t}\right]<\infty,
\end{equation*}
for $i=1,2$ and for some $\beta\in\mathbb R_+$, and $\eta^1,$ $\eta^2\in\mathcal H^2(\mathbb R^d).$
Then, for $\gamma, \delta\in(0,\beta]$, with $\gamma\neq \delta$
\begin{gather}
\begin{multlined}[c][0.9\displaywidth]
\int_0^T\textup{e}^{\beta A_t}\vert y_t^i\vert^2\,\textup{d} A_t 
\le
\frac{2}{\beta}\textup{e}^{\beta \Phi}\textup{e}^{\beta A_T}\sup_{t\in[0,T]}\mathbb E[ \vert \xi\vert^2|\mathcal G_t]
+\frac{2}{\gamma (\beta - \gamma)}
\textup{e}^{(\beta-\gamma) \Phi}\textup{e}^{(\beta-\gamma)A_T}
	\sup_{t\in[0,T]}\mathbb E\biggl[\int_0^T\textup{e}^{\beta A_t}\frac{\bigl\vert f^i_s \big\vert^2}{\alpha^2_s}\,\textup{d}C_s\bigg|\mathcal G_t\biggr],
\label{PathwiseAPriori-yi}
\end{multlined}
\shortintertext{and}
\int_0^T\textup{e}^{\beta A_t}\vert y_t^1-y_t^2\vert^2\,\textup{d} A_t 
\le	\frac{\textup{e}^{(\beta-\gamma)\Phi}}{\gamma(\beta-\gamma)}\textup{e}^{(\beta-\gamma)A_T}
	\sup_{t\in[0,T]}\mathbb E\biggl[\int_0^T\textup{e}^{\beta A_t}\frac{\bigl\vert f^1_s - f^2_s \big\vert^2}{\alpha^2_s}\,\textup{d}C_s\bigg|\mathcal G_t\biggr].
\label{PathwiseAPriori-Deltay}
\end{gather}
Moreover, for the martingale parts $\eta^1,\eta^2\in\mathcal H^2(\mathbb G;\mathbb R^n)$ of the aforementioned decompositions, we have
\begin{gather}
\sup_{t\in[0,T]}\big\vert \eta^1_t - \eta^1_0\big\vert^2
\le
6\sup_{t\in[0,T]} \mathbb E\bigg[ \big\vert\xi\big\vert^2 + \frac{1}{\beta}\int_0^T\textup{e}^{\beta A_s}\frac{\big\vert f^1_s\big\vert^2}{\alpha^2_s} \,\textup{d}C_s\bigg| \mathcal G_t \bigg]
 + 3\int_0^T\textup{e}^{\beta A_t} \frac{\big\vert f^1_s\big\vert^2}{\alpha^2_s}\,\textup{d} C_s,
 \label{PathwiseAPriori-etai}
\shortintertext{and}
\begin{multlined}[c][0.9\displaywidth]
\sup_{t\in[0,T]}\big\vert (\eta_t^1 - \eta_t^2)- (\eta_0^1-\eta_0^2 )\big\vert^2
\le
\frac{6}{\beta} \int_0^T\textup{e}^{\beta A_t} \frac{\big\vert f^1_s - f^2_s\big\vert^2}{\alpha^2_s}\,\textup{d} C_s 
+ \frac{3}{\beta} \sup_{t\in[0,T]}\mathbb E\bigg[ \int_0^T \textup{e}^{\beta A_s}\frac{\big\vert f^1_s -f^2_s\big\vert^2}{\alpha^2_s} \,\textup{d}C_s\bigg| \mathcal G_t\bigg].
\end{multlined}\label{PathwiseAPriori-Deltaeta}
\end{gather}
\end{lemma}

\begin{proof}
For the following assume $\gamma, \delta\in(0,\beta]$ with $\gamma\neq \delta.$
\begin{enumerate}[label=\raisebox{1.5pt}{\tiny$\bullet$}, itemindent=0.7cm]
	\item We will prove Inequality \eqref{PathwiseAPriori-yi} for $i=1$ by following analogous to Lemma \ref{lemmaest} calculations.
The sole difference will be that we are going to apply the conditional form of the Cauchy--Schwartz Inequality.
Moreover, by Identity \eqref{identity-yi-path}, we have
\begin{gather}
\big\vert y_t^1\big\vert^2 = \bigg\vert\mathbb E\bigg[ \xi^1 + \int_t^T f_s^1\,\textup{d}C_{s}\bigg|\mathcal G_t\bigg] \bigg\vert^2.
\label{identity-yi}
\end{gather}
In view of these comments, we have
\begin{align*}
&\int_0^T\textup{e}^{\beta A_t}\vert y_t^1\vert^2\,\textup{d} A_t 
\int_0^T\textup{e}^{\beta A_t}\bigg\vert\mathbb E\bigg[ \xi  + \int_t^Tf_s^1\,\textup{d}C_{s}\bigg|\mathcal G_t\bigg]\bigg\vert^2\,\textup{d} A_t\\
&\hspace{1em}\overset{\phantom{\text{C-S Ineq.}}}{\le}
2\int_0^T\textup{e}^{\beta A_t}\big\vert\mathbb E\big[ \xi\big|\mathcal G_t\big]\big\vert^2\,\textup{d} A_t  +
2\int_0^T\textup{e}^{\beta A_t}\bigg\vert\mathbb E\bigg[ \int_t^Tf_s^1\,\textup{d}C_{s}\bigg|\mathcal G_t\bigg]\bigg\vert^2\,\textup{d} A_t\\
&\hspace{1em}
\overset{\text{C-S Ineq.}}{\le}
2\int_0^T\textup{e}^{\beta A_t}\big\vert\mathbb E\big[ \xi\big|\mathcal G_t\big]\big\vert^2\,\textup{d} A_t
	 + 2\int_0^T \textup{e}^{\beta A_t} \mathbb E\Bigl[\frac{1}{\gamma}\textup{e}^{-\gamma A_t}\Big|\mathcal G_t\Bigr]
			\mathbb E\biggl[\int_t^T\textup{e}^{\gamma A_t}\frac{\bigl\vert f^1_s \big\vert^2}{\alpha^2_s}\,\textup{d}C_s\bigg|\mathcal G_t\biggr]\,\textup{d}A_t\\
&\hspace{1em}
\overset{\phantom{\text{C-S Ineq.}}}{\le}
2\int_0^T\textup{e}^{\beta A_t}\mathbb E\big[ \vert\xi\vert^2\big|\mathcal G_t\big]\,\textup{d} A_t
+\frac{2}{\gamma}\int_0^T \textup{e}^{(\beta-\gamma) A_t} 
		\mathbb E\biggl[\int_t^T\textup{e}^{\gamma A_t}\frac{\bigl\vert f^1_s\big\vert^2}{\alpha^2_s}\,\textup{d}C_s\bigg|\mathcal G_t\biggr]\,\textup{d}A_t\\
&\hspace{1em}
\overset{\phantom{\text{C-S Ineq.}}}{\le}
2\sup_{t\in[0,T]}\mathbb E\big[\vert\xi\vert^2\big|\mathcal G_t\big]
\int_0^T\textup{e}^{\beta A_t}\,\textup{d} A_t 
+\frac{2}{\gamma}
\int_0^T \textup{e}^{(\beta-\gamma) A_t} 
		\mathbb E\biggl[\int_0^T\textup{e}^{\beta A_t}\frac{\bigl\vert f^1_s\big\vert^2}{\alpha^2_s}\,\textup{d}C_s\bigg|\mathcal G_t\biggr]\,\textup{d}A_t\\
&\begin{multlined}[c][0.9\displaywidth]
\hspace{0.8em}
\underset{\phantom{\text{C-S Ineq.}}}{\overset{\text{Cor. \ref{maincorollary}}}{\le}}
2\sup_{t\in[0,T]}\mathbb E\big[ \vert\xi\vert^2\big|\mathcal G_t\big]
\int_0^T\textup{e}^{\beta A_t}\,\textup{d} A_t
+\frac{2}{\gamma}
	\sup_{t\in[0,T]}\mathbb E\biggl[\int_0^T\textup{e}^{\beta A_t}\frac{\bigl\vert f^1_s \big\vert^2}{\alpha^2_s}\,\textup{d}C_s\bigg|\mathcal G_t\biggr]
	\int_0^T\textup{e}^{(\beta-\gamma) A_t}\,\textup{d} A_t
\end{multlined}\\
&\begin{multlined}[c][0.9\displaywidth]
\hspace{0.8em}
\underset{\phantom{\text{C-S Ineq.}}}{\overset{\text{Cor. \ref{maincorollary}}}{\le}}
\frac{2}{\beta}\textup{e}^{\beta \Phi}\textup{e}^{\beta A_T}\sup_{t\in[0,T]}\mathbb E\big[ \vert\xi\vert^2\big|\mathcal G_t\big]
+\frac{2}{\gamma (\beta - \gamma)}
\textup{e}^{(\beta-\gamma) \Phi}\textup{e}^{(\beta-\gamma)A_T}
	\sup_{t\in[0,T]}\mathbb E\biggl[\int_0^T\textup{e}^{\beta A_t}\frac{\bigl\vert f^1_s \big\vert^2}{\alpha^2_s}\,\textup{d}C_s\bigg|\mathcal G_t\biggr].
\end{multlined}
\end{align*}

\item We will prove Inequality \eqref{PathwiseAPriori-Deltay}.
We will follow analogous arguments as in the previous case, but we are going to use  instead of \eqref{identity-yi} the identity 
\begin{gather}
\vert y^1_t- y^2_t\vert^2 = \bigg\vert\mathbb E\bigg[\int_t^Tf_s^1 - f_s^2\,\textup{d}C_{s}\bigg|\mathcal G_t\bigg]\bigg\vert^2.
\label{identity-Deltay}
\end{gather}
Now we have
\begin{align*}
&\int_0^T\textup{e}^{\beta A_t}\vert y_t^1-y_t^2\vert^2\,\textup{d} A_t 
	\overset{\eqref{identity-Deltay}}{=}\int_0^T\textup{e}^{\beta A_t}\bigg\vert\mathbb E\bigg[\int_t^Tf_s^1 - f_s^2\,\textup{d}C_{s}\bigg|\mathcal G_t\bigg]\bigg\vert^2\,\textup{d} A_t  \\
&\hspace{1em}
\overset{\text{C-S Ineq.}}{\le}
	\int_0^T \textup{e}^{\beta A_t} \mathbb E\Bigl[\frac{1}{\gamma}\textup{e}^{-\gamma A_t}\Big|\mathcal G_t\Bigr]
										\mathbb E\biggl[\int_t^T\textup{e}^{\gamma A_t}\frac{\bigl\vert f^1_s - f^2_s \big\vert^2}{\alpha^2_s}\,\textup{d}C_s\bigg|\mathcal G_t\biggr]\,\textup{d}A_t\\
&\hspace{1em}
\overset{\phantom{\text{C-S Ineq.}}}{\le}
	\frac{1}{\gamma}\int_0^T \textup{e}^{(\beta-\gamma) A_t} 
		\mathbb E\biggl[\int_0^T\textup{e}^{\beta A_t}\frac{\bigl\vert f^1_s - f^2_s \big\vert^2}{\alpha^2_s}\,\textup{d}C_s\bigg|\mathcal G_t\biggr]\,\textup{d}A_t\\
&\hspace{1em}
\underset{\phantom{\text{C-S Ineq.}}}{\overset{\text{Cor. \ref{maincorollary}}}{\le}}
	\frac{\textup{e}^{(\beta-\gamma)\Phi}}{\gamma(\beta-\gamma)}\textup{e}^{(\beta-\gamma)A_T}
	\sup_{t\in[0,T]}\mathbb E\biggl[\int_0^T\textup{e}^{\beta A_t}\frac{\bigl\vert f^1_s - f^2_s \big\vert^2}{\alpha^2_s}\,\textup{d}C_s\bigg|\mathcal G_t\biggr].
\end{align*}

\item Now we are going to prove \eqref{PathwiseAPriori-etai} for $i=1$.
We use initially the analogous to Inequality \eqref{ineq1} in order to obtain 
\begin{equation}\label{ineq1-eta1}
	\bigg\vert\int_0^T f^1_s\,\textup{d} C_s\bigg\vert^2 \le \frac{1}{\beta} \int_0^T\textup{e}^{\beta A_t} \frac{\big\vert f^1_s\big\vert^2}{\alpha^2_s}\,\textup{d} C_s.
\end{equation}
Moreover, by Identity \eqref{identity-yi} we obtain
\begin{align*}
\big\vert y^1_t\big\vert^2 
&\overset{\eqref{identity-yi}}{\le}
\bigg\vert \mathbb E\bigg[ \xi + \int_t^Tf^1_s \,\textup{d}C_s\bigg| \mathcal G_t \bigg]\bigg\vert^2
\le
2 \mathbb E\bigg[ \big\vert\xi\big\vert^2 + \bigg\vert\int_t^Tf^1_s \,\textup{d}C_s\bigg\vert^2\bigg| \mathcal G_t \bigg]\\
&\overset{\eqref{ineq1-eta1}}{\le}
2 \mathbb E\bigg[ \vert\xi\vert^2 + \frac{1}{\beta}\int_0^T\textup{e}^{\beta A_s}\frac{\big\vert f^1_s\big\vert^2}{\alpha^2_s} \,\textup{d}C_s\bigg| \mathcal G_t \bigg],
\end{align*}
and consequently
\begin{equation}
\sup_{t\in[0,T]}\big\vert y^1_t\big\vert^2 
\le
2\sup_{t\in[0,T]} \mathbb E\bigg[ \vert\xi\vert^2 + \frac{1}{\beta}\int_0^T\textup{e}^{\beta A_s}\frac{\big\vert f^1_s\big\vert^2}{\alpha^2_s} \,\textup{d}C_s\bigg| \mathcal G_t \bigg].
\label{ineq2-eta1}
\end{equation}
Now, by Identity \eqref{identity-yi-path} we have that
\begin{align*}
&\sup_{t\in[0,T]}\big\vert\eta^1_t - \eta^1_0\big\vert^2 = \sup_{t\in[0,T]}\bigg\vert y^1_t - y^1_0 + \int_0^Tf^1_s\,\textup{d}C_s\bigg\vert^2
\le 6\sup_{t\in[0,T]}\big\vert y^1_t\big\vert^2 
+ 3\bigg\vert\int_0^T f^1_s\,\textup{d} C_s\bigg\vert^2\\
&\hspace{2em}
\underset{\eqref{ineq2-eta1}}{\overset{\eqref{ineq1-eta1}}{\le}}
6\sup_{t\in[0,T]} \mathbb E\bigg[ \vert\xi\vert^2 + \frac{1}{\beta}\int_0^T\textup{e}^{\beta A_s}\frac{\big\vert f^1_s\big\vert^2}{\alpha^2_s} \,\textup{d}C_s\bigg| \mathcal G_t \bigg]
 + 3\int_0^T\textup{e}^{\beta A_t} \frac{\big\vert f^1_s\big\vert^2}{\alpha^2_s}\,\textup{d} C_s.
\end{align*}

\item We are going to prove, now, the Inequality \eqref{PathwiseAPriori-Deltaeta}.
We use initially the analogous to Inequality \eqref{ineq1} in order to obtain 
\begin{equation}\label{ineq1-Delta-eta-1}
\bigg\vert\int_0^T (f^1_s- f^2_s)\,\textup{d} C_s\bigg\vert^2 
	\le \frac{1}{\beta} \int_0^T\textup{e}^{\beta A_t} \frac{\big\vert f^1_s - f^2_s\big\vert^2}{\alpha^2_s}\,\textup{d} C_s.
\end{equation}
Moreover, by Identity \eqref{identity-Deltay} we have by Conditional Cauchy-Schwartz Inequality (analogously to the second case)
\begin{equation}\label{ineq1-Delta-eta-2}
\sup_{t\in[0,T]} \big\vert y^1_t - y^2_t\big\vert^2
\le \frac{1}{\beta} \sup_{t\in[0,T]}\mathbb E\bigg[ \int_0^T \textup{e}^{\beta A_s}\frac{\big\vert f^1_s -f^2_s\big\vert^2}{\alpha^2_s} \,\textup{d}C_s\bigg| \mathcal G_t\bigg].
\end{equation}
By Identity \eqref{identity-yi-path}, we have
\begin{equation*}
 (\eta_t^1 - \eta_t^2)- (\eta_0^1-\eta_0^2 )= (y_t^1 - y_t^2) - (y_0^1 - y_0^2) + \int_0^t (f_s^1 - f_s^2)\,\textup{d}C_{s}.
\end{equation*}
Finally, we have
\begin{align*}
&\sup_{t\in[0,T]}\big\vert (\eta_t^1 - \eta_t^2)- (\eta_0^1-\eta_0^2 )\big\vert
\le \sup_{t\in[0,T]}\bigg\vert(y_t^1 - y_t^2) - (y_0^1 - y_0^2) + \int_0^t (f_s^1 - f_s^2)\,\textup{d}C_{s}\bigg\vert^2\\
&\hspace{1em}\overset{\phantom{\eqref{ineq1-Delta-eta-2}}}{\le}
6\sup_{t\in[0,T]} \big\vert y^1_t - y^2_t\big\vert^2 
+\frac{3}{\beta}  \int_0^T\textup{e}^{\beta A_s} \frac{\big\vert f^1_s - f^2_s\big\vert^2}{\alpha_s^2}\,\textup{d}C_s\\
&\hspace{1em}\overset{\eqref{ineq1-Delta-eta-1}}{\underset{\eqref{ineq1-Delta-eta-2}}{\le}}
\frac{6}{\beta} \int_0^T\textup{e}^{\beta A_t} \frac{\big\vert f^1_s - f^2_s\big\vert^2}{\alpha^2_s}\,\textup{d} C_s 
+ \frac{3}{\beta} \sup_{t\in[0,T]}\mathbb E\bigg[ \int_0^T \textup{e}^{\beta A_s}\frac{\big\vert f^1_s -f^2_s\big\vert^2}{\alpha^2_s} \,\textup{d}C_s\bigg| \mathcal G_t\bigg].
\end{align*}

\end{enumerate}
\end{proof}

\begin{remark}
Viewing \eqref{BSDE2} as a BSDE whose generator does not depend on $y$ and $\eta$, then this BSDE has a solution, which can be uniquely determined by the pair $(y,\eta)$. 
Indeed, consider the data $(\FilG, T, \xi, f, C)$ and the processes $\alpha$ and $ A$, which all satisfy the respective assumptions of Lemma \ref{lemmaest} for some $\hat{\beta}>0$.
Then the semimartingale 
$$y_t=\E\bigg[\xi+\int_t^Tf_s \ud C_s\bigg|\G_t\bigg] = \E\bigg[\xi+\int_0^T f_s \ud C_s\bigg|\G_t\bigg] - \int_0^t f_s\ \ud C_s, \ t\in\Rp$$
satisfies $y_T=\xi$ and for $\eta_\cdot := \E\big[\xi+\int_0^T f_s\ \ud C_s\big|\G_{\cdot}\big]$
\begin{align*}
y_t-y_T	&= \E\bigg[\xi+\int_t^Tf_s \ud C_s\bigg|\G_t\bigg] - \xi
		 = \E\bigg[\xi+\int_0^Tf_s \ud C_s\bigg|\G_t\bigg] - \int_0^tf_s \ud C_s - \xi\\
		&= \E\bigg[\xi+\int_0^Tf_s \ud C_s\bigg|\G_t\bigg] + \int_t^T f_s \ud C_s - \int_0^T f_s \ud C_s - \xi= \eta_t + \int_t^T f_s \ud C_s -\eta_T.
\end{align*}
Now, one possible choice of a square--integrable $\FilG-$martingale $X$ such that $(X,\FilG, T, \xi, f, C)$ become standard data for any arbitrarily chosen integrator $C$, is the zero martingale.
Hence, given the standard data $(0,\FilG, T, \xi, f, C)$, the quadruple $(y,Z,U,\eta)$ satisfies the BSDE
$$y_t = \xi + \int_t^T f_s \ud C_s -\int_t^T\ud\eta_s,\ t\in\llbracket 0,T \rrbracket,$$
for any pair $(Z,U)$. 
Assume now that there exists a quadruple $(\widetilde{y},\widetilde{Z},\widetilde{U},\widetilde{\eta})$ which satisfies  
$$\widetilde{y}_t = \xi + \int_t^T f_s \ud C_s -\int_t^T\ud\widetilde{\eta}_s,\ t\in\llbracket 0,T \rrbracket.$$
Then, the pair $(y-\widetilde{y},\eta-\widetilde{\eta})$ satisfies
$$y-\widetilde{y}_t = -\int_t^T\ud\big(\eta-\widetilde{\eta}\big)_s,\ t\in\llbracket 0,T \rrbracket,$$
and by Lemma \ref{lemmaest}, for $\xi=0$ and $f=0$, we conclude that $\norm{y-\widetilde{y}}_{\mathcal{S}^2}=\norm{\eta-\widetilde{\eta}}_{\mathcal{H}^2}=0.$ Therefore $y$ and $\widetilde{y}$, resp. $\eta$ and $\widetilde{\eta}$, are indistinguishable, which implies our initial statement that every solution can be uniquely determined by the pair $(y,\eta)$. 
\end{remark}

In order to obtain the {\it a priori} estimates for the BSDE \eqref{BSDEintegralform}, we will have to consider solutions $(Y^i,Z^i,U^i,N^i),$ $i=1,2$, associated with the data $(X,\mathbb G,T,\xi^i,f^i,C), \ i=1,2$ under $\hat{\beta}$, where we also assume that $f^1,f^2$ have common $r,\vartheta$ bounds. Denote the difference between the two solutions by $(\delta Y, \delta Z,\delta U, \delta N),$ as well as $\delta\xi:=\xi^1-\xi^2$ and
$$ \delta_2 f_t := (f^1-f^2)(t,Y_t^2,Z_t^2,U_t^2(\cdot)),\ \psi_t:= f^1(t,Y^1_t,Z_t^1,U_t^1(\cdot)) - f^2(t,Y^2_t,Z_t^2,U_t^2(\cdot)).$$ 
We have the identity
\begin{align}\label{BSDE3}
	\delta Y_t =&\ \delta \xi + \int_t^T\psi_s \dC s- \int_t^T\delta Z_s \ud X^\circ_s - \int_t^T\int_{\R n}\delta U_s(x)\widetilde{\mu}^\natural(\ud s,\ud x) - \int_t^T \ud \delta N_s.
\end{align}
For the wellposedness of this last BSDE we need the following lemma.

\begin{lemma}
The processes
\begin{align*}
\int_0^\cdot\delta Z_s \ud X^\circ_s  
	\quad \textrm{ and } \quad 
\int_0^\cdot\int_{\R n}\delta U_s(x)\compd{\natural}
\end{align*}
are square-integrable martingales with  finite associated $\norm{\cdot}_{\hat\beta}-$norms.
\end{lemma}

\begin{proof}
The square-integrability is obvious. The inequalities 
\begin{align*}
\Expect{\tr [\,\pqc{\delta Z\cdot X^\circ]}}&\leq 2\Expect{\tr[\langle Z^1\cdot X^\circ\rangle]} +2\Expect{\tr[\langle Z^2\cdot X^\circ\rangle]},\\
\Expect{\tr [\,\pqc{\delta U\star\widetilde{\mu}^\natural}]}&\leq 2\Expect{\tr[\langle U^1\star\widetilde{\mu}^\natural\rangle]} +2\Expect{\tr[\langle U^2\star\widetilde{\mu}^\natural\rangle]},
\end{align*}
together with Lemma \ref{lemma1} guarantee that
\begin{equation*}
	\Expect{\int_0^T \e{\hat{\beta} A_t}\abs{c_t \delta Z_t}^2 \dC t} + 
	\Expect{\int_0^T \e{\hat{\beta} A_t}\tnorm{\delta U}^2_t \dC t} <\infty. \qedhere
\end{equation*}
\end{proof} 

Therefore, by defining
\begin{equation}\label{sumofmart}
	H_t := \int_0^t\delta Z_s \ud X^\circ_s + \int_0^t\int_{\R n}\delta U_s(x)\compd{\natural} + \int_0^t \ud \delta N_s,
\end{equation}
we can treat the BSDE \eqref{BSDE3} exactly as the BSDE \eqref{BSDE2}, where the martingale $H$ will play the role of the martingale $\eta$.

\begin{proposition}[\textit{A priori estimates for the {\rm BSDE} \eqref{BSDEintegralform}}]\label{aprioriBSDE}
Let $(X,\mathbb G,T,\xi^i,f^i,C),$ be standard data under $\hat{\beta}$ for $i=1,2$. 
Then $\psi/\alpha\in \mathbb{H}^2_{\hat{\beta}}$ and, if $M^\Phi(\hat\beta)<1/2$,
	 the following estimates hold 	
\begin{align*}
\norm{(\alpha\delta Y,\delta Z, \delta U, \delta N)}^2_{\hat\beta} 
	&\leq \widetilde{\Sigma}^{\Phi}(\hat\beta)  \norm{\delta \xi}^2_{\mathbb{L}^2_{\hat\beta}} 
		+ \Sigma^\Phi(\hat\beta) \norm{\frac{\delta_2 f}{\alpha}}_{\mathbb{H}^2_{\hat\beta}}^2 ,\\
\norm{(\delta Y,\delta Z, \delta U, \delta N)}^2_{\star,\hat\beta} 
	&\leq \widetilde{\Sigma}_\star^{\Phi}(\hat\beta)  \norm{\delta \xi}^2_{\mathbb{L}^2_{\hat\beta}} 
		+ \Sigma^\Phi_\star(\hat\beta) \norm{\frac{\delta_2 f}{\alpha}}_{\mathbb{H}^2_{\hat\beta}}^2,
		\end{align*}
where 
\begin{align*}
&\widetilde{\Sigma}^{\Phi}(\hat\beta) 
	:=\frac{\widetilde{\Pi}^{\hat\beta,\Phi}}{1-2M^\Phi(\hat\beta)},\; \widetilde{\Sigma}^{\Phi}_\star(\hat\beta) 
	:=\min\Big\{\widetilde{\Pi}_{\star}^{\hat\beta,\Phi} + 2M^\Phi_{\star}(\hat\beta), \widetilde{\Sigma}^{\Phi}(\hat\beta),\ 8 + \frac{16}{\hat\beta} \widetilde{\Sigma}^{\Phi}(\hat\beta) \Big\},\\
	&\Sigma^\Phi(\hat\beta) :=\frac{2M^\Phi(\hat\beta)}{1-2M^\Phi(\hat\beta)},\; \Sigma^\Phi_\star(\hat\beta) 
	:=\min\Big\{2 M^\Phi_{\star}(\hat\beta) (1 + \Sigma^\Phi(\hat\beta) ),\ \frac{16}{\hat\beta}( 1 + \Sigma^\Phi(\hat\beta)) \Big\}.
	\end{align*}
\end{proposition}
\begin{proof}
For the integrability of $\psi$, using the Lipschitz property \ref{Fthree} of $f^1,f^2$, we get
\begin{align*}
	|\psi_t|^2 \leq 	2r_t \abs{\delta Y_t}^2 + 	2\theta^\circ_t \norm{c_t \delta Z_t}^2 + 	2\theta^\natural_t \tnorm{\delta U}_t^2 +	2|\delta_2 f_t|^2.
\end{align*}
Hence by the definition of $\alpha$, which implies that $\frac{r}{\alpha^2}\leq \alpha^2$ and the obvious $\frac{\theta^\circ}{\alpha^2}, \frac{\theta^\natural}{\alpha^2}\leq 1$, we get
\begin{align*}
	\frac{|\psi_t|^2}{\alpha_t^2} 
	&\leq 2 \left(		\alpha_t^2 \abs{\delta Y_t}^2 + 
	\norm{c_t \delta Z_t}^2 + 		\tnorm{\delta U_t(\cdot)}_t^2 +
	\frac{\abs{\delta_2 f}^2}{\alpha^2}			\right)
	\numberthis\label{ineq10}\\
&\leq 2 \alpha_t^2 \abs{\delta Y_t}^2 + 
	2 \norm{c_t \delta Z_t}^2 + 
	2 \tnorm{\delta U_t(\cdot)}_t^2\\
&\hspace{0.9em}+ \frac{4}{\alpha^2}\left(\abs{f^1(s,0,0,\mathbf 0)}^2 + 
	r_t \abs{Y_t^2}^2 + 
	\theta^\circ_t \norm{c_t Z_t^2}^2 + 
	\theta^\natural_t \tnorm{\delta U^2_t(\cdot)}^2_t\right)\\
&\hspace{0.9em}+ \frac{4}{\alpha^2} \left(\abs{f^2(s,0,0,\mathbf 0)}^2 +
	r_t \abs{\delta Y_t}^2 + 
	\theta^\circ_t \norm{c_t \delta Z^2_t}^2 + 
	\theta^\natural_t \tnorm{\delta U_t(\cdot)}^2_t\right)\\
&\leq 6\big( \alpha_t^2 \abs{\delta Y_t}^2 + 
	\norm{c_t \delta Z_t}^2 + 
	 \tnorm{\delta U_t(\cdot)}_t^2\big)+ \frac{4}{\alpha^2}\left(\abs{f^1(s,0,0,\mathbf 0)}^2 + 
	\abs{f^2(s,0,0,\mathbf 0)}^2\right),
\end{align*}
where, having used once more that 
	$\frac{r}{\alpha^2}\leq \alpha^2$ and $\frac{\theta^\circ}{\alpha^2}, \frac{\theta^\natural}{\alpha^2}\leq 1$, it follows that 
	$\frac{\psi}{\alpha}\in \mathbb{H}^2_{\hat{\beta}} $. Next, for the $\norm{\cdot}_{\hat\beta}-$norm, we have 
	\begin{align*}
		\norm{\left(\delta Y, \delta Z, \delta U, \delta N \right)}^2_{\hat\beta} 
		&=\norm{\alpha \delta Y}^2_{\mathbb{H}^2_{\hat\beta}} +
		\norm{\delta Z}^2_{\mathbb{H}^{2,\circ}_{\hat\beta}} +
		\norm{\delta U}^2_{\mathbb{H}^{2,\natural}_{\hat\beta}} +
		\norm{\delta N}^2_{\mathcal{H}^{2,\perp}_{\hat\beta}}\\
		&\stackrel{(\ref{sumofmart})}{=}\norm{\alpha \delta Y}^2_{\mathbb{H}^2_{\hat\beta}} + 
		\norm{H}^2_{\mathcal{H}^{2}_{\hat\beta}}\\
		&\stackrel{(\ref{sumnormest})}{\leq} 
		\widetilde{\Pi}^{\hat\beta,\Phi} \norm{\delta \xi}^2_{\mathbb{L}^2_{\hat\beta}} +
		M^\Phi(\hat\beta)\norm{\frac{\psi}{\alpha}}^2_{\mathbb{H}^2_{\hat\beta}}\\
		&\leq \widetilde{\Pi}^{\hat\beta,\Phi} \norm{\delta \xi}^2_{\mathbb{L}^2_{\hat\beta}} +
		2 M^\Phi(\hat\beta) \left( \norm{\alpha \delta Y}^2_{\mathbb{H}^2_{\hat\beta}} +
		\norm{\delta Z}^2_{\mathbb{H}^{2,\circ}_{\hat\beta}} +
		\norm{\delta U}^2_{\mathbb{H}^{2,\natural}_{\hat\beta}} \right)	+	2M^\Phi(\hat\beta)		
		\norm{\frac{\delta_2f}{\alpha}}^2_{\mathbb{H}^2_{\hat\beta}}\\
		&\leq 
		\widetilde{\Pi}^{\hat\beta,\Phi} \norm{\delta \xi}^2_{\mathbb{L}^2_{\hat\beta}} +
		  2 M^\Phi(\hat\beta) \left( \norm{\alpha \delta Y}^2_{\mathbb{H}^2_{\hat\beta}} + \norm{H}^2_{\mathcal{H}^{2}_{\hat\beta}} \right) 
		+ 2 M^\Phi(\hat\beta) \norm{\frac{\delta_2 f}{\alpha}}^2_{\mathbb{H}^2_{\hat\beta}}.
	\end{align*}
Therefore, this implies
\begin{equation}\label{SigmaEstimates}
		\norm{(\alpha\delta Y,\delta Z, \delta U, \delta N)}^2_{\hat\beta} 
	\leq \widetilde{\Sigma}^{\Phi}(\hat\beta)  \norm{\delta \xi}^2_{\mathbb{L}^2_{\hat\beta}} 
		+ \Sigma^\Phi(\hat\beta) \norm{\frac{\delta_2 f}{\alpha}}_{\mathbb{H}^2_{\hat\beta}}^2 .
\end{equation}
We can obtain \emph{a priori} estimates for the $\norm{\cdot}_{\star,\hat\beta}-$norm by arguing in two different ways:

\vspace{0.5em}
$\bullet$ The identity \eqref{BSDE3} gives
\begin{align*}
\norm{\left(\delta Y, \delta Z, \delta U, \delta N \right)}^2_{\star,\hat\beta}
&=\norm{\delta Y}^2_{\mathcal{S}^2_{\hat\beta}} +
	\norm{\delta Z}^2_{\mathbb{H}^{2,\circ}_{\hat\beta}} +
	\norm{\delta U}^2_{\mathbb{H}^{2,\natural}_{\hat\beta}} +
	\norm{\delta N}^2_{\mathcal{H}^{2,\perp}_{\hat\beta}}\\
&
\overset{\eqref{sumofmart}}{=}\norm{\delta Y}^2_{\mathcal{S}^2_{\hat\beta}} + 
	\norm{H}^2_{\mathcal{H}^{2}_{\hat\beta}} 
\overset{\eqref{sumnormstarest}}{\leq} 
	\widetilde{\Pi}_{\star}^{\hat\beta,\Phi} \norm{\delta \xi}^2_{\mathbb{L}^2_{\hat\beta}} +
	M^\Phi_{\star}(\hat\beta) \norm{\frac{\psi}{\alpha}}^2_{\mathbb{H}^2_{\hat\beta}}\\
&
\stackrel{(\ref{ineq10})}{\leq}
	\widetilde{\Pi}_{\star}^{\hat\beta,\Phi} \norm{\delta \xi}^2_{\mathbb{L}^2_{\hat\beta}} +
	  2M^\Phi_\star(\hat\beta)\norm{\alpha\delta Y}^2_{\mathbb{H}^{2}_{\hat\beta}} +
	  2M^\Phi_\star(\hat\beta)\norm{H}^2_{\mathcal{H}^{2}_{\hat\beta}} 
	+ 2M^\Phi_\star(\hat\beta) \norm{\frac{\delta_2 f}{\alpha}}^2_{\mathbb{H}^2_{\hat\beta}}\\
&
\stackrel{(\ref{SigmaEstimates})}{\leq}
	\widetilde{\Pi}_{\star}^{\hat\beta,\Phi}\norm{\delta \xi}^2_{\mathbb{L}^2_{\hat\beta}} +
	  2M^\Phi_\star(\hat\beta)	\norm{\frac{\delta_2 f}{\alpha}}^2_{\mathbb{H}^2_{\hat\beta}}+ 2M^\Phi_\star(\hat\beta) 	
		\left(\widetilde{\Sigma}^{\Phi}(\hat\beta) \norm{\delta\xi}^2_{\mathbb{L}^2_{\hat\beta}} +	\Sigma^\Phi(\hat\beta)   \norm{\frac{\delta_2 f}{\alpha}}^2_{\mathbb{H}^2_{\hat\beta}}\right)\\
&
\hspace{0.2cm}=\left( \widetilde{\Pi}_{\star}^{\hat\beta,\Phi} + 2M^\Phi_\star(\hat\beta) \widetilde{\Sigma}^{\Phi}(\hat\beta)\right)
	\norm{\delta \xi}^2_{\mathbb{L}^2_{\hat\beta}} 	+	2M^\Phi_\star(\hat\beta) \left(1 + \Sigma^\Phi(\hat\beta) \right) 
	\norm{\frac{\delta_2 f}{\alpha}}^2_{\mathbb{H}^2_{\hat\beta}}.
\end{align*}

\vspace{0.5em}
$\bullet$ The identity \eqref{identity1} gives
\begin{align*}
	\norm{\delta Y}^2_{\mathcal{S}_{\hat\beta}^2}
&=	\Expect{\sup_{0\leq t\leq T}\left(\e{\frac{\hat\beta}{2}A_t} \abs{\delta Y_t} \right)^2}
	\stackrel{\eqref{identity1}}{\leq}
	\Expect{\sup_{0\leq t\leq T}\mathbb E\left[\left.\e{\frac{\hat\beta}{2}A_t} \abs{\delta\xi} 	+	\e{\frac{\hat\beta}{2}A_t} \abs{\int_t^T\psi_s \dC s}\right|\mathcal G_t\right]^2 }\\
&\hspace{-0.2cm}\stackrel{\eqref{ineq1}}{\leq} 
	2\Expect{\sup_{0\leq t\leq T} 
	\mathbb E\left[\left.\sqrt{\e{\hat\beta A_t} \abs{\delta \xi}^2 	+	\frac{1}{\hat\beta}\int_t^T\e{\hat\beta A_s} \frac{\abs{\psi_s}^2}{\alpha_s^2} \dC s}\right|\mathcal G_t\right]^2}\\
&\leq 	2\Expect{\sup_{0\leq t\leq T}
	\mathbb E\left[\left.\sqrt{\e{\hat\beta A_T} \abs{\delta\xi}^2 + 	\frac{1}{\hat\beta}\int_0^T\e{\hat\beta A_s} \frac{\abs{\psi_s}^2}{\alpha_s^2} \dC s	}\right|\mathcal G_t\right]^2}\\
&\leq 	8\Expect{	\e{\hat\beta A_T} \abs{\delta\xi}^2 +	\frac{1}{\hat\beta} \int_0^T \e{\hat\beta A_s}\frac{\abs{\psi_s}^2}{\alpha_s^2} \dC s}\\
&\hspace{-0.2cm}\stackrel{\eqref{ineq10}}{\leq} 	
	8 \Vert\delta\xi\Vert^2_{\mathbb L^2_{\hat\beta}} 
	+ \frac{16}{\hat\beta}\left( \norm{\frac{\delta_2 f}{\alpha}}^2_{\mathbb{H}^2_{\hat\beta}} 
	+ \norm{\alpha \delta Y}^2_{\mathbb{H}^2_{\hat\beta}}
	+ \norm{\delta Z}^2_{\mathbb{H}^{2,\circ}_{\hat\beta}} 
	+  \norm{\delta U}^2_{\mathbb{H}^{2,\natural}_{\hat\beta}}\right),
	\numberthis \label{ineq6}
\end{align*}
where, in the second and fifth inequality we used the inequality $a+b\leq \sqrt{2(a^2+b^2)}$ and Doob's inequality respectively. 
Then we can derive the required estimate
	\begin{align*}
		\norm{\left(\delta Y, \delta Z, \delta U, \delta N \right)}^2_{\star,\hat\beta} 
			&\stackrel{\hspace{0.8cm}}{=}\norm{\delta Y}^2_{\mathcal{S}^2_{\hat\beta}} +
		\norm{\delta Z}^2_{\mathbb{H}^{2,\circ}_{\hat\beta}} +
		\norm{\delta U}^2_{\mathbb{H}^{2,\natural}_{\hat\beta}} +
		\norm{\delta N}^2_{\mathcal{H}^{2,\perp}_{\hat\beta}}\\
			&\stackrel{\eqref{sumofmart}}{=}
		\norm{\delta Y}^2_{\mathcal{S}^2_{\hat\beta}} + \norm{H}^2_{\mathcal{H}^{2}_{\hat\beta}}\\
			&\underset{\eqref{sumofmart}}{\overset{\eqref{ineq6}}{\leq}}  
						8  \norm{\delta \xi}^2_{\mathbb{L}^2_{\hat\beta}} +
		\frac{16}{\hat\beta} \norm{\frac{\delta_2 f}{\alpha}}^2_{\mathbb{H}^2_{\hat\beta}} +
		\frac{16}{\hat\beta} \norm{\alpha \delta Y}^2_{\mathbb{H}^2_{\hat\beta}}+	
		\frac{16}{\hat\beta} \norm{ H }^2_{\mathcal{H}^{2}_{\hat\beta}}\\
			&\overset{\eqref{SigmaEstimates}}{\leq}
		\left( 8 + \frac{16}{\hat\beta} \widetilde{\Sigma}^{\Phi}(\hat\beta)\right) \norm{\delta \xi}^2_{\mathbb{L}^2_{\hat\beta}} +
		\frac{16}{\hat\beta}\left( 1 + \Sigma^\Phi(\hat\beta)\right) 
		\norm{\frac{\delta_2 f}{\alpha}}^2_{\mathbb{H}^2_{\hat\beta}}. \qedhere
	\end{align*}
\end{proof}

%% file: Proof_of_the_main_theorem.tex
We will use now the previous estimates to obtain the existence of a unique solution using a fixed point argument.

\begin{proof}[Proof of Theorem \ref{BSDEMainTheorem}.]\label{MainTheoremProof}
Let \  $\pair{y, z, u, n}$ be such that 
$(\alpha y, z, u, n)\in\Htwoshat{}\times \Htwoshat{,\circ}\times \Htwoshat{,\natural}\times\Hcal^{2,\perp}$.
Then the process $M$ defined by 
$$M_\cdot:=\mathbb E\left[\left.\xi + \int_0^T f(s,y_s,z_s,u_s(\cdot)) \ud C{s}\right|\mathcal G_\cdot\right]+n_\cdot\in\Hcal^2,$$ and by Proposition \ref{prop:OrthogDecomp}  it has a unique, up to indistinguishability, orthogonal decomposition
\[
M _\cdot	= M_0 	+	\int_0^{\cdot} Z_s  \ud X^\circ_s
			+	\int_0^{\cdot}\int_{\R{n}} U_s(x) \widetilde{\mu}^\natural(\ud s,\ud x)
			+ 	L_\cdot,
\]
where 
$\pair{Z,U,L}\in \mathbb{H}^{2,\circ}\times \mathbb{H}^{2,\natural}\times\Hcal^{2,\perp}$.
In view of the identity
\[
M_T - M_t 	=	\int_t^{T} Z_s \ud X^\circ_s
			+	\int_t^{T}\int_{\R{n}} U_s(x) \widetilde{\mu}^\natural(\ud s,\ud x)
			+ 	\int_t^T  \ud L_s,\ 0\leq t\leq T,
\]
we obtain
\begin{align*}
\mathbb E\left[\left.\xi + \int_t^T f(s,y_s,z_s,u_s(\cdot)) \ud C{s}\right|\mathcal G_t\right] =&\ 	\xi + 	\int_t^{T} f(s,y_s,z_s,u_s(\cdot)) \ud C{s}-	\int_t^{T} Z_s  \ud X^\circ_s\\
														&
														-	\int_t^{T}\int_{\R{n}} U_s(x) \widetilde{\mu}^\natural(\ud s,\ud x)
														- 	\int_t^{T}\ud N_s,
\end{align*}
where $N:=L-n$. Define 
$$Y_t:=\mathbb E\left[\left.\xi + \int_t^T f(s,y_s,z_s,u_s(\cdot)) \ud C{s}\right|\mathcal G_t\right].$$ 
In order to construct a contraction using Lemma \ref{lemmaest}, we need to choose $\delta>\gamma$.            
Then by Lemma \ref{lem:const} we can choose $\gamma^\star\in(0,\hat\beta]$ such that $\inf_{(\gamma,\delta) \in\mathcal C_{\hat\beta}}\Pi^\Phi(\gamma,\delta) = \Pi^\Phi(\gamma^\star(\hat\beta),\hat\beta)$.
Now we get that 
$(\alpha Y, Z\cdot X^\circ + U\star\mutilde + N)\in \mathbb H^2_{\hat\beta}\times\Hcal^2_{\hat\beta}$, and due to the orthogonality of the martingales we conclude that $(\alpha Y,Z,U,N)\in
\mathbb H^2_{\hat\beta}\times \mathbb H^{2,\circ}_{\hat\beta}\times \mathbb H^{2,\natural}_{\hat\beta}\times\mathcal H^{2,\perp}_{\hat\beta}$.
Hence, the operator
\[
S:	\mathbb H^2_{\hat\beta}\times \mathbb H^{2,\circ}_{\hat\beta}\times \mathbb H^{2,\natural}_{\hat\beta}\times\mathcal H^{2,\perp}_{\hat\beta} \longrightarrow 
	\mathbb H^2_{\hat\beta}\times \mathbb H^{2,\circ}_{\hat\beta}\times \mathbb H^{2,\natural}_{\hat\beta}\times\mathcal H^{2,\perp}_{\hat\beta},
\]
with the associated norms, that maps the processes $(\alpha y, z, u, n)$ to the processes $(\alpha Y, Z, U, N)$ defined above, is indeed well-defined. 

\vspace{0.5em}
Let $(\alpha y^i, z^i,u^i,n^i)\in \mathbb H^2_{\hat\beta}\times \mathbb H^{2,\circ}_{\hat\beta}\times \mathbb H^{2,\natural}_{\hat\beta}\times\mathcal H^{2,\perp}_{\hat\beta}$ for $i=1,2$, with 
\[
	S\left( \alpha y^i, z^i, u^i, n^i \right) = (\alpha Y^i, Z^i, U^i, N^i), \ \textrm{ for }\  i=1,2.
\]
Denote, as usual, $\delta y, \delta z, \delta u, \delta n$ the difference of the processes and $\psi_t:=f\left(t,y_t^1,z_t^1, u^1_t(\cdot)\right)-f\left(t,y_t^2,z_t^2,u_t^2(\cdot)\right).$
It is immediate that $\frac{\psi}{\alpha} \in \mathbb H^2_{\hat\beta}$ and that
\begin{align*}
	\norm{S\left( \alpha y^1, z^1, u^1, n^1\right) -S\left( \alpha y^2, z^2, u^2, n^2\right)}^2_{\hat\beta} 
&=	\norm{\alpha \delta Y}^2_{\mathbb H^{2}_{\hat\beta}} 
	+ 	\norm{\delta Z}^2_{\mathbb H^{2,\circ}_{\hat\beta}} 
	+ 	\norm{\delta U}^2_{\mathbb H^{2,\natural}_{\hat\beta}}
	+	\norm{\delta N}^2_{\Hcal^{2,\perp}_{\hat\beta}}\\
&\hspace{-0.35cm}\overset{\delta \xi=0}{\underset{\textrm{Lem. \ref{lemmaest}}}{\leq}}
	M^\Phi(\hat\beta) \norm{\frac{\psi}{\alpha}}^2_{\mathbb H^2_{\hat\beta}}\\
&\hspace{-0.2cm}\overset{\eqref{ineq10}}{\leq} 2M^\Phi(\hat\beta) \left(
		\norm{\alpha \delta y}^2_{\mathbb H^2_{\hat\beta}} 
	+ 	\norm{\delta z}^2_{\mathbb H^{2,\circ}_{\hat\beta}} 
	+ 	\norm{\delta u}^2_{\mathbb H^{2,\natural}_{\hat\beta}}
	\right) \\
&\leq	2M^\Phi(\hat\beta)
	\norm{\left( \alpha y^1, z^1, u^1, n^1\right) -\left( \alpha y^2, z^2, u^2, n^2\right)}^2_{\hat\beta}.
\end{align*}

Hence, for 
$M^\Phi(\hat\beta)<1/2$, we can apply Banach's fixed point theorem to obtain the existence of a unique fixed point $(\widetilde Y,Z, U,N)$. 
To obtain a solution in the desirable spaces we substitute $\widetilde{Y}$ in the quadruple with $Y$, the corresponding c\`adl\`ag version; indeed, $\mathbb G$ satisfies the usual conditions and $\widetilde Y$ is a semimartingale. The exact same reasoning using the $\norm{\cdot}_{\mathcal S^2_{\hat\beta}}-$norm for $Y$ leads to a contraction when $M^\Phi_\star(\hat\beta) < 1/2$.
\end{proof}

\begin{remark}\label{rem:proofMainTheorem}
Let us have a closer look at the proof of Theorem \ref{BSDEMainTheorem}.
In the following we adopt the notation introduced there.
Let us fix an initial point $(\alpha y^0,z^0,u^0,n^0)\in\mathbb H^2_{\hat{\beta}}\times\mathbb H^{2,\circ}_{\hat{\beta}}\times\mathbb H^{2,\natural}_{\hat{\beta}}\times\mathcal H^{2,\perp}_{\hat{\beta}}$ and define
$(\alpha y^k,z^k,u^k,n^k)\in\mathbb H^2_{\hat{\beta}}\times\mathbb H^{2,\circ}_{\hat{\beta}}\times\mathbb H^{2,\natural}_{\hat{\beta}}\times\mathcal H^{2,\perp}_{\hat{\beta}}$ by
\begin{align*}
(\alpha y^{k},z^{k},u^{k},n^{k}):=S(\alpha y^{k-1},z^{k-1},u^{k-1},n^{k-1}) \hspace{1em}\text{ for every }k\in\mathbb N.
\end{align*}
Let, moreover, $(Y,Z,U,N)$ denote the fixed-point.
Then, by Corollary \ref{cor:SpaceDecomposition} we can verify that
\begin{align*}
	z^k\cdot X^\circ \xrightarrow{\hspace{1em}\mathcal H^2\hspace{1em}}Z\cdot X^\circ,\hspace{1em}
	u^k\star \widetilde{\mu}^\natural\xrightarrow{\hspace{1em}\mathcal H^2\hspace{1em}} U\star\widetilde{\mu}^\natural
	\hspace{1em}\text{ and }\hspace{1em}
	n^k\xrightarrow{\hspace{1em}\mathcal H^2\hspace{1em}} N.
\end{align*} 
\end{remark}

\begin{corollary}[Picard approximation]\label{PicardCorollary}
	Assume that $M^\Phi(\hat\beta)<1/2$ $($resp. $M^\Phi_\star(\hat\beta)<1/2)$ 
	and define a sequence $(\Upsilon^{(p)})_{p\in\mathbb{N}}$
	on $\mathbb H^2_{\hat\beta}\times \mathbb H^{2,\circ}_{\hat\beta}\times \mathbb H^{2,\natural}_{\hat\beta}\times\Hcal^{2,\perp}_{\hat\beta}$
	$($resp. on $\mathcal S^2_{\hat\beta}\times \mathbb H^{2,\circ}_{\hat\beta}\times \mathbb H^{2,\natural}_{\hat\beta}\times\Hcal^{2,\perp}_{\hat\beta})$ 
	such that $\Upsilon^{(0)}$ is the zero element of the product space and $\Upsilon^{(p+1)}$ is the solution of 
	\begin{align*}
		Y^{(p+1)}_t =&\ \xi + \int_t^Tf(s,Y^{(p)}_s,Z^{(p)}_s,U^{(p)}_s(\cdot)) \ud C_s 
		- \int_t^T Z^{(p+1)}_s dX^\circ_s 
		- \int_t^T dN^{(p+1)}_s\\
		&-\int_t^T\int_{\R{n}}U^{(p+1)}_s(x) \widetilde{\mu}^\natural(\ud s,\ud x)
	\end{align*}
	Then 
\begin{enumerate}[label={\rm(\roman*)}, itemindent=0.8cm, leftmargin=0cm]
\item \label{Picardi} The sequence $(\Upsilon^{(p)})_{p\in\mathbb{N}}$ converges
	 in $\mathbb H^2_{\hat\beta}\times \mathbb H^{2,\circ}_{\hat\beta}\times \mathbb H^{2,\natural}_{\hat\beta}\times\Hcal^{2,\perp}_{\hat\beta}$
	$($resp. in $ \mathcal S^2_{\hat\beta}\times \mathbb H^{2,\circ}_{\hat\beta}\times \mathbb H^{2,\natural}_{\hat\beta}\times\Hcal^{2,\perp}_{\hat\beta})$
	to the solution of the BSDE \eqref{BSDE}.
\item \label{SeparateMartConv}
	The following convergence holds
	\begin{align*}
	\big(Z^{(p)}, U_s^{(p)}, N^{(p)} \big)
		\xrightarrow[p\to\infty]{\hspace{1cm}}
	(Z, U, N),\ \text{in }\ \mathbb H^{2,\circ}_{\hat{\beta}}\times\mathbb H^{2,\natural}_{\hat{\beta}}\times\mathcal H^{2,\perp}_{\hat{\beta}}.
	\end{align*}
	\item\label{Picardii} There exists a subsequence $(\Upsilon^{(p_m)})_{m\in\mathbb{N}}$ which converges $e^{\hat\beta A}\ud C \otimes \ud\Pm-a.e.$
\end{enumerate}

\end{corollary}

\begin{proof}
As in the proof of Theorem \ref{BSDEMainTheorem}, we obtain
\begin{equation}\label{CauchySequence}
	\norm{\Upsilon^{(p+1)}-\Upsilon^{(p)}}^2_{\hat\beta} 
		\leq \left( 2M^\Phi(\hat\beta)\right)^p \norm{\Upsilon^{(1)}}^2_{\hat\beta} \ 
	\left(\text{resp. }\norm{\Upsilon^{(p+1)}-\Upsilon^{(p)}}^2_{\star,\hat\beta} 
		\leq \left( 2M^\Phi_\star(\hat\beta)\right)^p \norm{\Upsilon^{(1)}}^2_{\star,\hat\beta}\right),
\end{equation}
	and consequently, since $\sum_{p\in\mathbb{N}}\norm{\Upsilon^{(p+1)} -\Upsilon^{(p)}}^2_{\hat\beta} <\infty$ (resp. $\sum_{p\in\mathbb{N}}\norm{\Upsilon^{(p+1)} -\Upsilon^{(p)}}^2_{\star,\hat\beta} <\infty$), 
	the sequence $(\Upsilon^{(p)})_{p\in\mathbb{N}}$ is Cauchy 
	in $\mathbb H^2_{\hat\beta}\times \mathbb H^{2,\circ}_{\hat\beta}\times \mathbb H^{2,\natural}_{\hat\beta}\times\Hcal^{2,\perp}_{\hat\beta}$
	(resp. in $ \mathcal S^2_{\hat\beta}\times \mathbb H^{2,\circ}_{\hat\beta}\times \mathbb H^{2,\natural}_{\hat\beta}\times\Hcal^{2,\perp}_{\hat\beta}$).
	Denote by $\Upsilon$ the unique limit on the product space. Then, it coincides with the unique fixed point 
	for the contraction $S$ (see the proof of Theorem \ref{BSDEMainTheorem} above) due to the construction of
	$(\Upsilon^{(p)})_{p\in\mathbb{N}}$, which proves \ref{Picardi}. 

\vspace{0.5em}
	For \ref{SeparateMartConv}, the result is immediate by the Cauchy property of the sequence $(\Upsilon^{(p)})_{p\in\mathbb{N}}$ and Corollary \ref{cor:SpaceDecomposition}\footnote{The reader may recall the Remark \ref{rem:proofMainTheorem}.}.

\vspace{0.5em}
	Finally, for \ref{Picardii}, by the $\norm{\cdot}_{\hat\beta}-$convergence, we can extract a subsequence $\set{p_m}_{m\in\mathbb{N}}$ such that
	\begin{equation}\label{SuitableRate}
		\norm{\Upsilon^{(p_{m+1})} -\Upsilon^{(p_{m})}}_{\hat\beta}\leq 2^{-2m},\ \text{for every $m\geq 0$}.
	\end{equation}
	Define, for any $\varepsilon\geq 0$, $N^{p,\varepsilon}:=
\set{\left(\omega,t\right)\in\Omega\times\stochint{0}{T},\ |Y^{(p)}_t(\omega)-Y_t(\omega) |>\varepsilon }$. Then we have
\begin{align*}	
	e^{\hat\beta A}\ud C\otimes \ud \Pm \left(\limsup_{m\to \infty} N^{p_m,\varepsilon}	\right) 
	&=\lim_{m\to \infty} e^{\hat\beta A}\ud C\otimes \ud\Pm
	\left(\bigcup_{\ell=m}^{\infty} \left[ \abs{Y^{(p_\ell)}-Y}	> 	\varepsilon	\right]	\right) \\
	&\leq \lim_{m\to \infty} \frac{1}{\varepsilon^2}\sum_{\ell=m}^{\infty} 
				\Expect{\int_0^T \e{\hat\beta A_t} \abs{Y^{(p_\ell)}_t-Y_t}^2 \ud C{t}}\\
	&\leq \lim_{m\to \infty} \frac{1}{\varepsilon^2}\sum_{\ell=m}^{\infty} \norm{Y^{(p_\ell)}-Y}^2_{\hat\beta}\\
		&\leq \lim_{m\to \infty} \frac{1}{\varepsilon^2}\sum_{\ell=m}^{\infty} 
	\left( \sum_{n=1}^{\infty} 2^n\norm{Y^{(p_{\ell+n+1})}-Y^{(p_{\ell+n})}}^2_{\mathbb H^2_{\hat\beta}} \right)\\
	&\leq \lim_{m\to \infty} \frac{1}{\varepsilon^2}\sum_{m=\ell}^{\infty} 
	\left( \sum_{n=1}^{\infty} 2^n\norm{\Upsilon^{(p_{\ell+n+1})}-\Upsilon^{(p_{\ell+n})}}^2_{\hat\beta} \right)\\
	&\stackrel{\eqref{SuitableRate}}{\leq} \lim_{m\to \infty} \frac{1}{\varepsilon^2}\sum_{\ell=m}^{\infty} 
	\left( \sum_{n=1}^{\infty} 2^n 2^{-2(\ell+n)} \right)=0,\ \text{for any $\varepsilon >0$}.
\end{align*}

Hence
\begin{align*}
	e^{\hat\beta A}\ud C\otimes \ud\Pm\left(\limsup_{m\to \infty} N^{p_m,0}	\right) 
	\leq \sum_{n\in\mathbb{N}}e^{\hat\beta A}\ud C\otimes \ud\Pm
	\left(\limsup_{m\to \infty} N^{p_m,1/n}	\right) = 0.
\end{align*}

Following the same arguments, we have the almost sure convergence of $Z^{p_m}, U^{p_m}, N^{p_m} $ to the corresponding processes of the $\norm{\cdot}_{\beta}-$solution of the BSDE \eqref{BSDE}.
Moreover, using the same steps, we can obtain the analogous result for the $\norm{\cdot}_{\star,\hat{\beta}}-$norm.
\end{proof}

%% file: New_estimates.tex

In this subsection we derive the \emph{a priori} estimates for the BSDE \eqref{BSDEintegralform} by means of an alternative method. 
It is essentially the classical one used to obtain estimates in a BSDE setting, namely apply It\=o's formula to an appropriately weighted $\mathbb L^2-$type norm of the $Y$ part of the solution, and then take conditional expectations. 
We will see that even though this approach still works in this setting (albeit with significant complications) and leads to sufficient conditions for wellposedness which are very similar to the ones obtained in \cite{bandini2015existence,cohen2012existence}, it also requires an additional assumption, which is completely inherent to the approach, and turns out to be slightly restrictive in terms of applications, see Remark \ref{rem:H6} for more details.

\vspace{0.5em}
Let us, initially, introduce some auxiliary processes. 
Let $\varepsilon$ be a $\mathbb G-$predictable process such that $\varepsilon_s(\omega)\ge \Delta C_s(\omega)$, for $\ud C\otimes\ud \mathbb P-a.e$ $(s,\omega)\in\mathbb R_+\times\Omega$. 
Fix, moreover, a non-negative, $\mathbb G-$predictable process $\gamma$ and define the increasing, $\mathbb G-$predictable and \cadlag process
\begin{equation}\label{eq:v}
	v_\cdot :=\int_0^{\cdot} \gamma_s \dC{s}.
\end{equation}

Let $\mathcal E$ denote the stochastic exponential operator.
The following assumptions will be in force throughout this subsection\footnote{The reader may recall the notation introduced at the beginning of the section.}.
\vspace{1em}
\begin{enumerate}[label=\bf{(H\arabic*)},leftmargin=*,itemindent=0.0cm]
	\item\label{Hone} The martingale $\oX$ belongs to $\mathcal H^2(\mathbb R^m)\times \mathcal H^2(\mathbb R^n)$ and $(\oX,C)$ satisfies Assumption \ref{assumptionC}.
	\item\label{Htwo} The terminal condition $\xi$ satisfies $\mathbb E\big[\mathcal E(v)_T\xi^2\big]<\infty$. 
	\item\label{Hthree} The generator of the equation $f:\domain $ is such that for any 
				$(y,z,u)\in\R d\times \R {d\times m}\times \mathfrak H$, 
				the map 
				$$(t,\omega)\longmapsto f(t,\omega,y,z,u_t(\omega;\cdot))\ \text{is $\mathbb G-$predictable.}$$ 
				Moreover, $f$ satisfies a stochastic Lipschitz condition\footnote{This is exactly the same as \ref{Fthree}.}, that is to say there exist 
				\[
					r:\left(\Omega\times\mathbb{R}_+,\Pred\right)\longrightarrow(\Rp,\borel{\Rp})\ \textrm{ and }\
					\vartheta =(\theta^\circ,\theta^\natural):\left(\Omega\times\mathbb{R}_+,\Pred\right)\longrightarrow(\Rp^2,\borel{\Rp^2}),
				\] such that, for $\ud C \otimes \ud \Pm-a.e.$ $(t,\omega)\in\mathbb R_+\times\Omega$
				\begin{align}\label{fstochLip-2}
				\begin{split}
					& \quad \left|f(t,\omega,y,z,u_t(\omega;\cdot))-f(t,\omega,y',z',u_t'(\omega;\cdot))\right|^2\\
							&\hspace{2em}\leq 	r_t(\omega) |y-y'|^2 	+ \theta^\circ_t(\omega)\Vert c_t(\omega) (z-z')\Vert^2 + \theta^\natural_t(\omega) \left(\tnorm{u_t(\omega;\cdot)-u'_t(\omega;\cdot)}_t(\omega)\right)^2.		
				\end{split}
				\end{align}
	\item\label{Hfour} We have 
 		 	\begin{equation*}
 				\E\left[\int_0^T\mathcal E(v)_{s-}(1+\gamma_s\Delta C_s)(\varepsilon_s-\Delta C_s)\abs{f(s,0,0,\mathbf 0)}^2 \dC{s}\right] <\infty,
 			\end{equation*}
 	where $\mathbf 0$ denotes the null application from $\mathbb R^n$ to $\mathbb R$.
	\item\label{Hfive} For every pair $Y^1,Y^2\in\mathbb H^2$, the measure $\ud C\otimes \ud \mathbb P$ is such that $Y^1,Y^2$ are equal $\ud C\otimes \ud \mathbb P-a.e.$ if and only if $Y^1_{-},Y^2_{-}$ are equal $\ud C\otimes \ud \mathbb P-a.e.$
\item\label{Hsix}
 		If $X^{\circ,c}\ne 0$ then one of the following is true $\ud C\otimes\ud \mathbb P-a.e.$
		
		\vspace{0.5em}
			\begin{enumerate}[label=(\roman*)]
				\item  $\displaystyle C_{s-}(\theta^\circ_s\vee\theta^{\natural}_s)<1$ and 
				$\displaystyle r_s<\min\bigg\{\frac{(\theta^\circ_s\vee\theta^{\natural}_s)(1-C_{s-}(\theta^\circ_s\vee\theta^{\natural}_s))}{C_s(1+(\theta^\circ_s\vee\theta^{\natural}_s)\Delta C_s)},\frac{(\sqrt{C_s}-\sqrt{C_{s-}})^2}{C_s(\Delta C_s)^2}\bigg\}$,
				
				\vspace{0.5em}
				\item  $\displaystyle C_{s-}\Delta C_s\theta^\circ_s<C_s$ and $\displaystyle (\Delta C_s)^2 r_s<\min\bigg\{\frac{\Delta C_s+C_{s-}(1-\theta^\circ_s\Delta C_s)}{C_s},\frac{(\sqrt{C_s}-\sqrt{C_{s-}})^2}{C_s}\bigg\}$,
				
				\vspace{0.5em}
				\item  $ \displaystyle (\theta^\circ_s\vee\theta^{\natural}_s) C_{s-}<1$ and $ \displaystyle r_s<\min\bigg\{ \frac{(\sqrt{C_s}-\sqrt{C_{s-}})^2}{C_s(\Delta C_s)^2},\frac{1-\theta^\circ_sC_{s-}}{C_s\Delta C_s}\bigg\}$.
			\end{enumerate}

			\vspace{0.8em}
\noindent If $X^{\circ,c}=0$ then one of the following holds true $\ud C\otimes\ud \mathbb P-a.e.$
 		\begin{enumerate}[label=(\roman*)]
				\item $\displaystyle r_s<\min\bigg\{\frac{(\sqrt{C_s}-\sqrt{C_{s-}})^2}{C_s(\Delta C_s)^2},\frac{\theta^\circ_s\vee\theta^{\natural}_s(1-C_{s-}(\theta^\circ_s\vee\theta^{\natural}_s))}{C_s(1+(\theta^\circ_s\vee\theta^{\natural}_s)\Delta C_s)}\bigg\}$,
				\vspace{0.5em}
				\item $\displaystyle C_{s-}\Delta C_s\theta^\circ_s<C_s$ and $\displaystyle (\Delta C_s)^2 r_s<\min\bigg\{\frac{\Delta C_s+C_{s-}(1-\theta^\circ_s\Delta C_s)}{C_s},\frac{(\sqrt{C_s}-\sqrt{C_{s-}})^2}{C_s}\bigg\}$.
		\end{enumerate}
\end{enumerate}

For this subsection we understand the term \emph{standard data} as follows: we will say that the sextuple $(\Fil,\oX,T,\xi,f,C)$ are \emph{standard data}, whenever its elements satisfy Assumptions \ref{Hone}--\ref{Hsix}.
Therefore, we also modify the definition of a solution of the BSDE \eqref{BSDEintegralform} given the standard data $(\Fil,\oX,T,\xi,f,C)$.

\begin{definition}\label{solutiondef-newapriori}
A \emph{solution of the {\rm BSDE} \eqref{BSDEintegralform} with standard data $(\oX,\mathbb G, T,\xi,f,C)$} is a quadruple of processes 
\begin{align*}
(Y,Z,U,N)&\in 	\mathbb H^2 \times	\mathbb H^{2,\circ} \times	\mathbb H^{2,\natural} \times	\mathcal H^{2,\perp}
  \textrm{ or } (Y,Z,U,N)\in 		
					\mathcal{S}^2 \times	\mathbb H^{2,\circ} \times	\mathbb H^{2,\natural} \times	\mathcal H^{2,\perp}
\end{align*}
such that, $\Pm-a.s.$, for any $t\in\llbracket0,T\rrbracket$, 
	\[Y_t = \xi 	+ \int_t^T f(s,Y_{s-},Z_s,U_s(\cdot)) \dC{s} 
				- \int_t^T \ Z_s \ud X^\circ_s 
				-\int_t^T\int_{\R{n}}U_s(x) \widetilde\mu^\natural(\ud s,\ud x)
				- \int_t^T \ud N_s.
\]
\end{definition}

\begin{remark}
In order to obtain the \emph{a priori} estimates by means of the method described in the next sub--sub--section, we will need to distinguish between the cases $X^{\circ,c}\ne 0$ and $X^{\circ,c}=0$.
After doing so, the analysis that we are going to make will lead us to specific conditions that the processes $C, r$ and $\vartheta$ should satisfy.
These conditions are described in the sub--parts of Assumption \ref{Hsix}.
\end{remark}

\begin{remark}
The attentive reader may have observed that the conditions appearing in \ref{Hsix} never impose a lower bound on the process $r$. Notice however that the analysis that will be carried out in Appendix \ref{App:AuxAnalysis} could, in principle, lead to additional sufficient conditions imposing such lower bounds. We have decided to ignore them because we wanted to be able to always include in our framework the case where the generator of the BSDE does not depend on the solution $Y$.
\end{remark}

\begin{remark}\label{rem:H6}
Assumption \ref{Hfive} is the tricky one here, and already appears in the work of Cohen and Elliott \cite{cohen2012existence}. The main point is that the fixed point argument which will be used here only allows to define uniquely the process of left--limits of $Y$, $\ud C \otimes \ud \Pm-a.e.$ Without Assumption \ref{Hfive}, we cannot define $Y$ itself from its left--limits alone. We emphasize that we did not require this condition with our first approach, and that it is inherent to the current approach and cannot be avoided with this method. This is the main advantage of our approach. We will now present two situations where it is actually satisfied.
\begin{enumerate}[label=(\roman*)]
\item 
	If $C$ is a deterministic process such that $C$ assigns positive measure to every non--empty open subinterval of $\mathbb R_+$, then Condition \ref{Hfive} is satisfied; see \cite[Lemma 5.1]{cohen2012existence}. 
		
	\vspace{0.5em}
\item\label{rem:C-ii}
The previous case is somehow of limited interest, as it excludes the case where $C$ is a piecewise--constant, increasing integrator. This corresponds to a so--called backward stochastic difference equation (BS$\Delta$E), and would be the object of interest in numerical schemes where one would approximate the martingale driving the BSDE by, for instance, appropriate random walks. The following describes how we can allow for some discrete--time processes $C$.

\vspace{0.5em}
\noindent Let us initially describe the properties the standard data should have in order to embed a BS$\Delta$E into a continuous--time framework.
Let $\pi:=\{0=t_0<t_1<\dots<t_n<\dots\}$ be a partition of $\mathbb R_+$ and $(\mathbb G,\oX,T,\xi,f,C)$ be standard data with the following properties.
\begin{enumerate}[label=\textbullet,itemindent=0.7cm]
\item The filtration $\mathbb G:=(\mathcal G_t)_{t\in\mathbb R_+}$ is such that $\mathcal G_t=\mathcal G_{t_n}$ for every $t\in[t_n,t_{n+1})$ and for every $n\in\mathbb N\cup\{0\}$.
\item The martingale $\oX$ is such that $\oX_t=\oX_{t_n}$ $\mathbb P-a.s.$ for every $t\in[t_n,t_{n+1})$ and for every $n\in\mathbb N\cup\{0\}$. 
\item The generator $f$ is such that $f(t,\omega,y,z,u(\omega;\cdot))=f(t_n,\omega,y,z,u(\omega;\cdot))$, $\mathbb P-a.s.$ for every $t\in[t_n,t_{n+1})$ and for every $n\in\mathbb N\cup\{0\}$, $y\in\mathbb R^d$, 
		$z\in\mathbb R^{m\times m}$ and 
		$u(\omega;\cdot):(\mathbb R^n,\mathcal B(\mathbb R^n))\longrightarrow (\mathbb R^d,\mathcal B(\mathbb R^d))$.
\item The integrator $C$ is of the form
\begin{align*}
C_{\cdot} = C(0) + \sum_{n\in\mathbb N}C(n) \mathds{1}_{[t_n,t_{n+1})}(\cdot),
\end{align*}
where $\big(C(n)\big)_{n\in\mathbb{N}\cup\{0\}}$ is a sequence of non-negative random variables such that
\begin{enumerate}[label={\tiny$\blacksquare$}]
\item the random variable $C(0)$ is $\mathcal G_0-$measurable and the random variable $C(n)$ is $\mathcal G_{t_{n-1}}-$measurable for every $n\in\mathbb N$,
\item for every $n\in\mathbb N\cup\{0\}$ holds $0\le C(n)\le C({n+1})$ $\mathbb P-a.s.$
\end{enumerate}
\end{enumerate}

\vspace{0.5em}
\noindent Let now $\tau$ be a stopping time for the (discretely--indexed) filtration $(\mathcal G_{t_n})_{n\in\mathbb N\cup\{0\}},$ \emph{i.e.} $\tau\in\pi$ $\mathbb P-a.s.,$ and $[\tau=t_n]\in\mathcal G_{t_n}$, for every $n\in\mathbb N\cup\{0\}.$
Let us assume, moreover, that there are $\bar{t}_1<\bar{t}_2\in\pi\backslash\{0\}$ such that $\mathbb P([\tau=\bar{t}_i])>0$ for every $i\in\{1,2\}.$
Additionally, let us assume that $\mathbb P([\Delta C_{\bar{t}_2}> 0]\cap[\tau=\bar{t}_1])>0$.
We can assume without loss of generality that there exists a $\delta>0$ such that  $\mathbb P([\Delta C_{\bar{t}_2}> \delta]\cap[\tau=\bar{t}_1])>0$.
Define, now, the stopping times $\sigma_1,\sigma_2$ as follows
\begin{gather*}
\sigma_1,\sigma_2\phantom:=\tau \text{ on } \Omega\backslash \bigcup_{t_n\ne \bar{t}_1}[\tau=t_n],\; \sigma_1:=\bar{t}_2 \text{ and }
\sigma_2:=s\in(\bar{t}_1, \bar{t}_2) \text{ on }[\tau=\bar{t}_1].
\end{gather*}
Then, $\ud \mathbb P\otimes \ud C\big(\mathds{1}_{\llbracket\sigma_1, \infty\llbracket}\neq\mathds{1}_{\llbracket\sigma_2, \infty\llbracket} \big)=0$, however
\begin{align*}
\ud \mathbb P\otimes \ud C\big(\mathds{1}_{\rrbracket\sigma_1, \infty\llbracket}\neq\mathds{1}_{\rrbracket\sigma_2, \infty\llbracket} \big)
	=\ud \mathbb P\otimes \ud C\big(\{\bar{t}_2\}\times[\tau=\bar{t}_1]\big)
	=\mathbb E[\Delta C_{\bar{t}_2} \mathds{1}_{[\tau=\bar{t}_1]}\big]
	\ge \delta \mathbb P([\Delta C_{\bar{t}_2}\ge \delta]\cap[\tau=\bar{t}_1])>0.
\end{align*}
If, however, we restrict ourselves in the subspace 
	\begin{align*}
		\mathbb H^2_\pi:=\big\{Y\in\mathbb H^2, (Y_{t_n})_{n\in\mathbb N\cup\{0\}} \text{ is adapted to } (\mathcal G_{t_n})_{n\in\mathbb N\cup\{0\}} \text{ and } Y_t=Y_{t_n} \text{ on } [t_n,t_{n+1}) \text{ for every } n\in\mathbb N\cup\{0\}\big\},
	\end{align*}
then, under the additional assumption that $\mathbb P(\Delta C_{t_n}>0)=1$ for every $n\in\mathbb N\cup\{0\}$, we have that $Y^1,Y^2\in\mathbb H^2_\pi$ are equal $\ud \mathbb P\otimes\ud C-a.e.$, if and only if $Y^1_{-},Y^2_{-}$ are equal $\ud \mathbb P\otimes\ud C-a.e.$, where $Y_{0-}:=Y_0$ for every $Y\in\mathbb H^2_\pi$.
In this case, we can also conclude that $Y^1$ and $Y^2$ are indistinguishable.
However, the reader may observe that for $Y^1,Y^2\in\mathbb H^2$ such that $Y^1,Y^2$ are equal $\ud \mathbb P\otimes \ud C-a.e.$ we cannot conclude that they are indistinguishable, as we concluded above. 
\end{enumerate}
\end{remark}

\subsubsection{New estimates}\label{sebsub:NewEst}
As we have already mentioned, we are going to derive the \emph{a priori} estimates for BSDE \eqref{BSDEintegralform}. 
To this end, let us fix a $d-$dimensional, $\mathbb G-$predictable process $h$ and consider the BSDE 
\begin{equation}\label{BSDEintegralform2}
	Y_t = \xi 	+ \int_t^T h_s \dC s 
				- \int_t^T Z_s \ud X^\circ_s 
				-\int_t^T\int_{\R{n}}U_s(x) \widetilde\mu^{\natural}(\ud s,\ud x)
				- \int_t^T \ud N_s,
\end{equation}
for some $(Z,U,N)\in\mathbb H^{2,\circ}\times\mathbb H^{2,\natural}\times\mathcal H^{2,\perp}$.
Moreover, we abuse notation (see Footnote \ref{foot:semimart-d}) and for the finite variation process $C$ we define $C^d_\cdot := \sum _{s\le \cdot} \Delta C_s.$
Our first result is the following estimates, which in conjunction with Theorem \ref{th:estimates2} can be seen as the analogous of Lemma \ref{lemmaest}.
\begin{theorem}\label{th:estimates}
For any positive $\mathbb G-$predictable process $(\varepsilon_t)_{t\geq 0}$, if $(Y,Z,U,N)\in 	\mathbb H^{2}\times \mathbb H^{2,\circ}\times  \mathbb H^{2,\natural}\times \mathcal H^{2,\perp}$
 solves {\rm BSDE} \eqref{BSDEintegralform2}, we have the estimate 
\begin{align*}
\mathcal E(v)_t\vert Y_t\vert^2&+\mathbb E\bigg[\int_t^T\mathcal E(v)_{s-}\big(\gamma_s-(1+\gamma_s\Delta C_s)\varepsilon_s^{-1}\big)\vert Y_{s-}\vert^2\dC s\bigg|\mathcal G_t\bigg]\\
&+\mathbb E\bigg[\int_t^T\big(\mathcal E(v)_{s}-\Delta\mathcal E(v)_s\mathds{1}_{\{X^{\circ,c}\neq 0\}} \big)\ud \mathrm{Tr}[\langle Z\cdot X^\circ\rangle]_s\bigg|\mathcal G_t\bigg]
+\mathbb E\bigg[\int_t^T{\mathcal E(v)_s}\ud\mathrm{Tr}[\langle U\star\widetilde{\mu}^\natural\rangle]_s\bigg|\mathcal G_t\bigg]\\
&+\mathbb E\bigg[\int_t^T\mathcal E(v)_{s-}\ud \mathrm{Tr}[\langle N\rangle]_s
+\int_t^T\Delta\mathcal E(v)_s  \mathds{1}_{\{X^{\circ,c}=0\}\cup\{X^{\circ,d}=0\}}\ud \mathrm{Tr}[\langle N^d\rangle]_s\bigg|\mathcal G_t\bigg]\\
&\hspace{10em}\leq\mathbb E\bigg[\mathcal E(v)_T\vert\xi\vert^2+\int_t^T\mathcal E(v)_{s-}(1+\gamma_s\Delta C_s)(\varepsilon_s-\Delta C_s)\vert h_s\vert^2\dC s\bigg|\mathcal G_t\bigg].
\end{align*}
\end{theorem}

\begin{proof}
In the following we are going to use the identities $[C] = [C^d]$ and $$\sum_{0\le s\le \cdot}\Delta C_s \Delta L_s=[C^d,L^d]_\cdot = [C,L^d]_\cdot = [C,L]_\cdot$$ for every semimartingale $L$\footnote{\label{foot:semimart-d}Recall that a semimartingale $L$ can be written in the form $L=L_0 + M + A$, where $M$ is a local martingale and $A$ is a finite variation process. We will denote by $L^d$ the process $M^d+A$.}, which are true due to the fact that $C$ is of finite variation; see \cite[Theorem I.4.52]{jacod2003limit}. 
We also use the fact that since $v$ is predictable and has finite variation, so is $\mathcal E(v)$. Moreover,  $\Delta \mathcal E(v)_t(\omega)\neq 0$ if and only if $\Delta C_t(\omega)\neq 0$.
This allows us to write $(\Delta\mathcal E(v))\cdot C^d$ as $(\Delta \mathcal E(v))\cdot C$.

\vspace{0.5em}
Let us fix an $i=1,\ldots,d$.
We apply It\=o's product rule in order to calculate the differential of the process $\mathcal E(v)(Y^i)^2$ and obtain
\begin{align*}
\mathcal E(v)_\cdot (Y_\cdot^i)^2
&=
					 \mathcal E(v)_0(Y_0^i)^2
					-2\int_0^\cdot\mathcal E(v)_{s-}Y^i_{s-}h^i_s\ud C_s 
 					+2\underbrace{\int_0^\cdot\mathcal E(v)_{s-}Y^i_{s-}\ud(Z\cdot X^\circ)^i_s }_{\text{martingale}} 
					+2\underbrace{\int_0^\cdot\mathcal E(v)_{s-}Y^i_{s-}\ud(U\star \widetilde{\mu}^\natural)^i_s}_{\text{martingale}} 					\\
&\phantom{=}
 					+2\underbrace{\int_0^\cdot\mathcal E(v)_{s-}Y^i_{s-}\ud N^i_s}_{\text{martingale}}
 					+\int_0^\cdot\mathcal E(v)_{s-}(h^i_s)^2\ud [C]_s
 					-2\underbrace{\int_0^\cdot\mathcal E(v)_{s-}h^i_s\ud[C,(Z\cdot X^\circ)^i]_s}_{\text{martingale; see \cite[Proposition I.4.49]{jacod2003limit}}}+\int_0^\cdot\mathcal E(v)_{s-}(Y^i_{s-})^2\gamma_s\ud C_s
\\
&\phantom{=}
 					-2\underbrace{\int_0^\cdot\mathcal E(v)_{s-}h^i_s\ud[C,U^i\star\widetilde{\mu}^\natural]_s}_{\text{martingale}}
 					-2\underbrace{\int_0^\cdot\mathcal E(v)_{s-}h^i_s\ud[C,N^i]_s}_{\text{martingale}}
 					+\int_0^\cdot\mathcal E(v)_{s-}\ud[(Z\cdot X^\circ)^i]_s\\
&\phantom{=}
 					+2\underbrace{\int_0^\cdot\mathcal E(v)_{s-}\ud[(Z\cdot X^\circ)^i,U^i\star\widetilde{\mu}^\natural]_s}_{\text{martingale, since }M_{\mu}[\Delta X^\circ|\widetilde{\mathcal P}]=0}
 					+2\underbrace{\int_0^\cdot\mathcal E(v)_{s-}h_s\ud[(Z\cdot X^\circ)^i,N^i]_s}_{\text{martingale, since }\langle X^\circ,N\rangle=0}+\int_0^\cdot\mathcal E(v)_{s-}\ud[U\star\widetilde{\mu}^\natural]_s\\
&\phantom{=}
 					+2\underbrace{\int_0^\cdot\mathcal E(v)_{s-}\ud[U^i\star\widetilde{\mu}^\natural,N^i]_s}_{\text{martingale, since } M_{\mu}[\Delta N|\widetilde{\mathcal P}]=0}
 					+\int_0^\cdot\mathcal E(v)_{s-}\ud[N^i]_s\\
&\phantom{=}		
					+\int_0^\cdot \Delta\mathcal E(v)_s \ud \big(\big[\big((Z\cdot X^\circ)^i\big)^d\big] + [U^i\star\widetilde{\mu}^\natural] + [(N^i)^d]\big)_s
					+\int_0^\cdot \Delta\mathcal E(v)_sh^i_s(h^i_s \Delta C_s-2Y^i_{s-}) \ud C_s\\
&\phantom{=}		+\underbrace{\big[\big(2\mathcal E(v)_{s-}\gamma_sY^i_{s-}-2\Delta\mathcal E(v)_sh^i_s\big)\cdot C, \big((Z\cdot X^\circ)^i\big)^d + U^i\star\widetilde{\mu}^\natural + (N^i)^d \big]_\cdot}_{\text{martingale; see \cite[Proposition I.4.49]{jacod2003limit}}}\\
&\phantom{=} 		+\int_0^\cdot\Delta \mathcal E(v)_{s}\big\{2 \ud \underbrace{\big[\big((Z\cdot X^\circ)^i\big)^d, U^i\star\widetilde{\mu}^\natural]_s}_{\text{martingale}} +2 \ud [\big((Z\cdot X^\circ)^i\big)^d,(N^i)^d]_s +2 \ud \underbrace{[ U^i\star\widetilde{\mu}^\natural,(N^i)^d]_s}_{\text{martingale} }  \big\}.
				\numberthis\label{Stop1}
\end{align*}
By writing Identity \eqref{Stop1} in its integral form on the interval $[t,T]$ and by taking conditional expectation with respect to the $\sigma-$algebra $\mathcal G_t$, we deduce
\begin{align*}
\mathbb E\big[\mathcal E(v)_T (\xi^i)^2\big|\mathcal G_t\big] - \mathcal E(v)_t(Y^i_t)^2 
&=					 
	\mathbb E\bigg[\int_t^T\mathcal E(v)_{s-}(h^i_s)^2\ud [C]_s-2\int_t^T\mathcal E(v)_{s-}Y^i_{s-}h^i_s\ud C_s+\int_t^T\mathcal E(v)_{s-}\ud[(Z\cdot X^\circ)^i]_s \bigg|\mathcal G_t\bigg]
	\\
&\phantom{=}
 	+\mathbb E\bigg[\int_t^T\mathcal E(v)_{s-}\ud[U^i\star\widetilde{\mu}^\natural]_s+\int_t^T\mathcal E(v)_{s-}\ud[N^i]_s+\int_t^T\mathcal E(v)_{s-}(Y^i_{s-})^2\gamma_s\ud C_s\bigg|\mathcal G_t\bigg]
 	\\
&\phantom{=}
	+\mathbb E\bigg[\int_t^T\Delta\mathcal E(v)_sh^i_s(h^i_s \Delta C_s-2Y^i_{s-}) \ud C_s+\int_t^T\Delta\mathcal E(v)_s \ud \big[\big((Z\cdot X^\circ)^i\big)^d\big]_s\bigg|\mathcal G_t\bigg]\\
&\phantom{=} 	
	+\mathbb E\bigg[\int_t^T\Delta\mathcal E(v)_s\ud [U^i\star\widetilde{\mu}^\natural]_s+2\int_t^T\Delta \mathcal E(v)_s \ud \big[\big((Z\cdot X^\circ)^i\big)^d,(N^i)^d\big]_s\bigg|\mathcal G_t\bigg]\\
&\phantom{=}
	+\mathbb E\bigg[\int_t^T\Delta\mathcal E(v)_s \ud [(N^i)^d]_s\bigg|\mathcal G_t\bigg].
\end{align*}
Reordering the terms in the above equality we obtain
\begin{align*}
0\le &\ \mathcal E(v)_t(Y_t^i)^2 
+\mathbb E\bigg[\int_t^T\mathcal E(v)_{s-}\ud[(Z\cdot X^\circ)^i]_s
+\int_t^T\Delta\mathcal E(v)_s \ud \big[\big((Z\cdot X^\circ)^i\big)^d\big]_s\bigg|\mathcal G_t\bigg]\\
&
+\mathbb E\bigg[\int_t^T\mathcal E(v)_{s-}\ud[U^i\star\widetilde{\mu}^\natural]_s\bigg|\mathcal G_t\bigg]
+\mathbb E\bigg[\int_t^T\Delta\mathcal E(v)_s\ud [U^i\star\widetilde{\mu}^\natural]_s\Big|\mathcal G_t\bigg]\\
&+\mathbb E\bigg[\int_t^T\mathcal E(v)_{s-}\ud[N^i]_s
+\int_t^T\Delta\mathcal E(v)_s \ud [(N^i)^d]_s\bigg|\mathcal G_t\bigg]
+\mathbb E\bigg[\int_t^T\mathcal E(v)_{s-}(Y^i_{s-})^2\gamma_s\ud C_s\Big|\mathcal G_t\bigg]\\
=&	\				 
\mathbb E\big[\mathcal E(v)_T (\xi^i)^2\big|\mathcal G_t\big]
	+2\mathbb E\bigg[\int_t^T\mathcal E(v)_{s-}Y^i_{s-}h^i_s\ud C_s \Big|\mathcal G_t\bigg]
	-\mathbb E\bigg[\int_t^T\mathcal E(v)_{s-}(h^i_s)^2\ud [C]_s\Big|\mathcal G_t\bigg]\\
&	-\mathbb E\bigg[\int_t^T\Delta\mathcal E(v)_sh^i_s(h^i_s \Delta C_s-2Y^i_{s-}) \ud C_s\Big|\mathcal G_t\bigg]
	-2\mathbb E\bigg[\int_t^T\Delta \mathcal E(v)_s \ud \big[\big((Z\cdot X^\circ)^i\big)^d,(N^i)^d\big]_s\Big|\mathcal G_t\bigg]\\
\le&\					 
\mathbb E\big[\mathcal E(v)_T (\xi^i)^2\big|\mathcal G_t\big]
	+2\mathbb E\bigg[\int_t^T\mathcal E(v)_{s-}Y^i_{s-}h^i_s\ud C_s \Big|\mathcal G_t\bigg]
	-\mathbb E\bigg[\int_t^T\mathcal E(v)_{s-}(h^i_s)^2\Delta C_s \ud C_s\Big|\mathcal G_t\bigg]\\
&	-\mathbb E\bigg[\int_t^T\Delta\mathcal E(v)_sh^i_s(h^i_s \Delta C_s-2Y^i_{s-}) \ud C_s\Big|\mathcal G_t\bigg]\\	
&	+\mathbb E\bigg[\int_t^T\Delta \mathcal E(v)_s \ud \big[\big((Z\cdot X^\circ)^i)^d] +\int_t^T\Delta \mathcal E(v)_s \ud [(N^i)^d]_s\bigg|\mathcal G_t\bigg],
	\numberthis\label{Stop3}
\end{align*}
where we obtained Inequality \eqref{Stop3} by using the Kunita--Watanabe inequality and then Young's inequality for the summand 
$2\mathbb E\big[\int_t^T\Delta \mathcal E(v)_s \ud \big[\big((Z\cdot X^\circ)^i\big)^d,(N^i)^d\big]_s\big|\mathcal G_t\big]$.
More precisely, we have 
\begin{align*}
-2\int_t^T\Delta (\mathcal E(v))_s\ud \big[\big((Z\cdot X^\circ)^i\big)^d,(N^i)^d\big]_s
&\le
		2\int_t^T\Delta (\mathcal E(v))_s\ud \text{Var}\big(\big[\big((Z\cdot X^\circ)^i\big)^d,(N^i)^d\big]\big)_s\\
&\le 
		2\Big(\int_t^T\Delta (\mathcal E(v))_s\ud \big[\big((Z\cdot X^\circ)^i\big)^d\big]_s\Big)^\frac12
		 \Big(\int_t^T\Delta (\mathcal E(v))_s\ud [(N^i)^d]_s\Big)^\frac12\\
&\le 
		\int_t^T \Delta (\mathcal E(v))_s\ud \big[\big((Z\cdot X^\circ)^i\big)^d]_s 
	+ 	\int_t^T \Delta (\mathcal E(v))_s\ud [(N^i)^d]_s.
\end{align*}

Therefore, by Inequality \eqref{Stop3} we obtain
\begin{align*}
0\le &\ \mathcal E(v)_t(Y^i_t)^2 
+\mathbb E\bigg[\int_t^T\mathcal E(v)_{s-}(Y^i_{s-})^2\gamma_s\ud C_s\Big|\mathcal G_t\bigg]
+\mathbb E\bigg[\int_t^T\mathcal E(v)_{s-}\ud[(Z\cdot X^\circ)^i]_s\Big|\mathcal G_t\bigg]\\
&+\mathbb E\bigg[\int_t^T\underbrace{(\mathcal E(v)_{s-} + \Delta\mathcal E(v)_s)}_{\mathcal E(v)_s}\ud[U^i\star\widetilde{\mu}^\natural]_s\Big|\mathcal G_t\bigg]
+\mathbb E\bigg[\int_t^T\mathcal E(v)_{s-}\ud[N^i]_s\Big|\mathcal G_t\bigg]\\
\le& \ 					 
\mathbb E\big[\mathcal E(v)_T (\xi^i)^2\big|\mathcal G_t\big]
	+\mathbb E\bigg[\int_t^T\mathcal E(v)_{s-}h^i_s(1+\gamma_s\Delta C_s)(2Y^i_{s-}-h^i_s\Delta C_s)\ud C_s \Big|\mathcal G_t\bigg].
	\numberthis\label{Stop4}
\end{align*}
Notice now that if $X^\circ\in\mathcal H^{2,c}$ or $X^\circ\in\mathcal H^{2,d}$, then the term $\mathbb E\big[ \int_t^T\Delta \mathcal E(v)_s \ud \big[\big((Z\cdot X^\circ)^i\big)^d, (N^i)^d\big]_s\big|\mathcal G_t\big]$ in the right--hand side of Identity \eqref{Stop5} vanishes. 
This is true since in the former case $\big[\big((Z\cdot X^\circ)^i\big)^d,(N^i)^d\big]=0$, while in the latter case the process $\big[\big((Z\cdot X^\circ)^i\big)^d,(N^i)^d\big]$ is a uniformly integrable martingale; recall that by the Galtchouk--Kunita--Watanabe decomposition we have that $\langle X^\circ,N\rangle=0$ and since $X^\circ=X^{\circ,d}$ we can easily conclude.
Therefore, we can incorporate the above special cases into Inequality \eqref{Stop4} as follows
\begin{align*}
0\le &\ \mathcal E(v)_t(Y^i_t)^2 
+\mathbb E\bigg[\int_t^T\mathcal E(v)_{s-}\ud[(Z\cdot X^\circ)^i]_s\bigg|\mathcal G_t\bigg]
+\mathbb E\bigg[\int_t^T\Delta\mathcal E(v)_s \mathds{1}_{\{X^\circ\in\mathcal H^{2,d}\}}\ud \big[\big((Z\cdot X^\circ)^i\big)^d\big]_s\bigg|\mathcal G_t\bigg]\\
&+\mathbb E\bigg[\int_t^T\mathcal E(v)_{s-}\ud[U^i\star\widetilde{\mu}^\natural]_s\bigg|\mathcal G_t\bigg]
+\mathbb E\bigg[\int_t^T\mathcal E(v)_{s-}(Y^i_{s-})^2\gamma_s\ud C_s\bigg|\mathcal G_t\bigg]
+\mathbb E\bigg[\int_t^T\Delta\mathcal E(v)_s\ud [U^i\star\widetilde{\mu}^\natural]_s\bigg|\mathcal G_t\bigg]\\
&+\mathbb E\bigg[\int_t^T\mathcal E(v)_{s-}\ud[N^i]_s\bigg|\mathcal G_t\bigg]
+\mathbb E\bigg[\int_t^T\Delta\mathcal E(v)_s  \mathds{1}_{\{X^\circ\in\mathcal H^{2,c}\}\cup\{X^\circ\in\mathcal H^{2,d}\}}\ud [(N^i)^d]_s\bigg|\mathcal G_t\bigg]\\
\leq& \					 
\mathbb E\big[\mathcal E(v)_T (\xi^i)^2\big|\mathcal G_t\big]
	+\mathbb E\bigg[\int_t^T\mathcal E(v)_{s}h^i_s(2Y^i_{s-}-h^i_s \Delta C_s) \ud C_s\bigg|\mathcal G_t\bigg].
\numberthis\label{Stop5}
\end{align*}
This rewrites equivalently as
\begin{align*}
0\le &\ \mathcal E(v)_t(Y^i_t)^2 +\mathbb E\bigg[\int_t^T\mathcal E(v)_{s-}(Y^i_{s-})^2\gamma_s\ud C_s\Big|\mathcal G_t\bigg]
+\mathbb E\bigg[\int_t^T(\mathcal E(v)_{s-}\mathds{1}_{\{X^{\circ,c}\neq 0\}} 
+ \mathcal E(v)_s \mathds{1}_{\{X^{\circ,c}=0\}})\ud[(Z\cdot X^\circ)^i]_s\bigg|\mathcal G_t\bigg]\\
&+\mathbb E\bigg[\int_t^T\mathcal E(v)_{s}\ud[U^i\star\widetilde{\mu}^\natural]_s\Big|\mathcal G_t\bigg]
+\mathbb E\bigg[\int_t^T\mathcal E(v)_{s-}\ud[N^i]_s
+\int_t^T\Delta\mathcal E(v)_s  \mathds{1}_{\{X^{\circ,c}=0\}\cup\{X^{\circ,d}=0\}}\ud [(N^i)^d]_s\bigg|\mathcal G_t\bigg]\\
\leq&	\				 
\mathbb E\big[\mathcal E(v)_T (\xi^i)^2\big|\mathcal G_t\big]
	+\mathbb E\bigg[\int_t^T\mathcal E(v)_sh^i_s(2Y^i_{s-}-h^i_s \Delta C_s) \ud C_s\Big|\mathcal G_t\bigg].
\end{align*}
Next, for any positive $\mathbb G-$predictable process $(\varepsilon_t)_{t\geq 0}$, we have the estimate
\begin{align*}
h^i_s\big(2Y^i_{s-}-h^i_s\Delta C_s\big)\leq \varepsilon_s^{-1}(Y^i_{s-})^2+(\varepsilon_s-\Delta C_s)(h^i_s)^2,
\end{align*}
so that we deduce the desired result using that $[L]-\langle L\rangle$ is a uniformly integrable martingale for any $L\in\mathcal H^2$. Finally, taking the sum for $i=1,\ldots, d$ we obtain the required estimates
\begin{align*}
&\hspace{-3em}\mathcal E(v)_t\vert Y_t\vert^2+\mathbb E\bigg[\int_t^T\mathcal E(v)_{s-}\big(\gamma_s-(1+\gamma_s\Delta C_s)\varepsilon_s^{-1}\big)\vert Y_{s-}\vert^2\dC s\bigg|\mathcal G_t\bigg]\\
&+\mathbb E\bigg[\int_t^T\big(\mathcal E(v)_{s}-\Delta\mathcal E(v)_s\mathds{1}_{\{X^{\circ,c}\neq 0\}} \big)\ud \mathrm{Tr}[\langle Z\cdot X^\circ\rangle]_s\bigg|\mathcal G_t\bigg]
+\mathbb E\bigg[\int_t^T{\mathcal E(v)_s}\ud\mathrm{Tr}[\langle U\star\widetilde{\mu}^\natural\rangle]_s\bigg|\mathcal G_t\bigg]\\
&+\mathbb E\bigg[\int_t^T\mathcal E(v)_{s-}\ud \mathrm{Tr}[\langle N\rangle]_s
+\int_t^T\Delta\mathcal E(v)_s  \mathds{1}_{\{X^{\circ,c}=0\}\cup\{X^{\circ,d}=0\}}\ud \mathrm{Tr}[\langle N^d\rangle]_s\bigg|\mathcal G_t\bigg]\\
&\hspace{10em}\leq\mathbb E\bigg[\mathcal E(v)_T\vert\xi\vert^2+\int_t^T\mathcal E(v)_{s-}(1+\gamma_s\Delta C_s)(\varepsilon_s-\Delta C_s)\vert h_s\vert^2\dC s\bigg|\mathcal G_t\bigg]. \qedhere
\end{align*}

\end{proof}
Our next result provides now estimates for the difference of the solutions of two BSDEs.
\begin{theorem}\label{th:estimates2}
Fix some processes $(y,\bar y, z,\bar z,u,\bar u)\in\big(\mathbb H^{2}(\oX)\big)^2\times
				\big(\mathbb H^{2}(X^\circ)\big)^2\times
				\big(\mathbb H^{2}(X^\natural)\big)^2$ and consider the following two {\rm BSDEs}
\begin{align*}
	Y_t &= \xi 	+ \int_t^T f_s(y_{s-},z,u)\dC s 
				- \int_t^T Z_s \ud X^\circ_s 
				-\int_t^T\int_{\R{n}}U_s(x) \widetilde\mu^\natural(\ud s,\ud x)
				- \int_t^T \ud N_s,\\
	\overline{Y}_t &= \xi 	+ \int_t^T f_s(\bar y_{s-},\bar z,\bar u)\dC s 
				- \int_t^T \overline{Z}_s \ud X^\circ_s 
				-\int_t^T\int_{\R{n}}\overline{U}_s(x) \widetilde\mu^\natural(\ud s,\ud x)
				- \int_t^T \ud \overline{N}_s,
\end{align*}
where $f$ is the process in \ref{Hthree}.
Denoting $\delta L:=L-\overline{L}$, for $L=y,Y,z,Z,u,U,N$, we have
\begin{align*}
&\mathbb E\bigg[\int_0^T\mathcal E(v)_{s-}\big(1+C_s\big(\gamma_s-(1+\gamma_s\Delta C_s)\varepsilon_s^{-1}\big)\big)\vert\delta Y_{s-}\vert^2\dC s\bigg]\\
&\hspace{0.9em}+\mathbb E\bigg[\int_0^TC_s\big(\mathcal E(v)_{s}-\Delta\mathcal E(v)_s\mathds{1}_{\{X^{\circ,c}\neq 0\}}\big)\ud \mathrm{Tr}[ \langle \delta Z\cdot X^\circ\rangle]_s\bigg]
+\mathbb E\bigg[\int_0^T\mathcal E(v)_{s}C_s\ud \mathrm{Tr}[\langle \delta U\star \widetilde{\mu}^\natural\rangle]_s\bigg]\\
&\hspace{0.9em}+\mathbb E\bigg[\int_0^T\mathcal E(v)_{s-}C_s\ud\mathrm{Tr}[\langle\delta N\rangle]_s\bigg] 
+ \mathbb E\bigg[\int_0^T\Delta \mathcal E(v)_s C_s\mathds{1}_{\{X^{\circ,c}=0\}\cup\{X^{\circ,d}=0\}}\ud \mathrm{Tr}[\langle \delta N^d\rangle]_s\bigg]\\
&\hspace{4em}\leq\mathbb E\bigg[\int_0^T\mathcal E(v)_{s-}(1+\gamma_s\Delta C_s)C_s(\varepsilon_s-\Delta C_s)\big(r_s |\delta y_{s-}|^2
+\theta^{\circ}_s\ud \mathrm{Tr}[\langle \delta z\cdot X^\circ\rangle]_s 
+ \theta^{\natural}_s\ud \mathrm{Tr}[\langle \delta u\star\widetilde{\mu}^\natural\rangle]_s\big)\bigg].
\numberthis\label{Contraction}
\end{align*}
\end{theorem}
\begin{proof}
First, we have
\begin{align*}
\delta Y_t =&\  \int_t^T f_s( y_{s-},z_s,u_s(\cdot))-f_s(\bar y_{s-},\bar z_s,\bar u_s(\cdot)) \dC s 
				- \int_t^T \delta Z_s \ud X^\circ_s 
				-\int_t^T\int_{\R{n}}\delta U_s(x) \widetilde\mu^\natural(\ud s,\ud x)- \int_t^T \ud \delta N_s.
\end{align*}
The estimate from Theorem \ref{th:estimates} and the Lipschitz property of $f$ ensure then that 
\begin{align*}
&\mathcal E(v)_t\vert\delta Y_t\vert^2
+\mathbb E\bigg[\int_t^T\mathcal E(v)_{s-}\big(\gamma_s-(1+\gamma_s\Delta C_s)\varepsilon^{-1}\big)\vert\delta Y_{s-}\vert^2\dC s\bigg|\mathcal G_t\bigg]+\mathbb E\bigg[\int_t^T{\mathcal E(v)_s}\ud\mathrm{Tr}[\langle \delta U\star\widetilde{\mu}^\natural\rangle]_s\bigg|\mathcal G_t\bigg]\\
&\hspace{0.9em}+\mathbb E\bigg[\int_t^T\big(\mathcal E(v)_{s}-\Delta\mathcal E(v)_s\mathds{1}_{\{X^{\circ,c}\neq 0\}} \big)\ud \mathrm{Tr}[\langle \delta Z\cdot X^\circ\rangle]_s\bigg|\mathcal G_t\bigg]\\
&\hspace{0.9em}+\mathbb E\bigg[\int_t^T\mathcal E(v)_{s-}\ud \mathrm{Tr}[\langle \delta N\rangle]_s
+\int_t^T\Delta\mathcal E(v)_s  \mathds{1}_{\{X^{\circ,c}=0\}\cup\{X^{\circ,d}=0\}}\ud \mathrm{Tr}[\langle \delta N^d\rangle]_s\bigg|\mathcal G_t\bigg]\\
&\leq\mathbb E\bigg[\int_t^T\mathcal E(v)_{s-}(1+\gamma_s\Delta C_s)(\varepsilon_s-\Delta C_s)\big(r_s|\delta y_{s-}|^2+\theta^{\circ}_s\norm{c_s\delta z_s}^2+\theta^{\natural}_s\tnorm{\delta u_s}^2\big)\dC s\bigg|\mathcal G_t\bigg].
\end{align*}
Writing the integrals above as difference of integrals on $(0,T]$ and $(0,t]$, we can take left--limits for $t\uparrow u$ (use \cite[Theorem 3.4.11]{he1992semimartingale} and apply L\'evy's upward theorem), then integrate between $0$ and $T$ with respect to the measure $\ud C$, which, using \cite[Theorem 8.2.5]{cohen2015stochastic} for the predictable projection and Fubini's theorem, leads us to the desired estimate.
\end{proof}
The previous result leads to naturally define new norms. 
More precisely, assume that $v$ and $\varepsilon$ additionally satisfy 
$$1+C_s\big(\gamma_s-(1+\gamma_s\Delta C_s)\varepsilon_s^{-1}> 0,\; {\rm d}\mathbb P\otimes\dC s-a.e.,$$
and define for any $(Y,Z,U,N)$ with appropriate dimensions and measurability the following norms, as well as the associated spaces
$$\norm{Y}_{\mathbb H^2(\oX,v)}^2:=\E\bigg[\int_0^T\mathcal E(v)_{s-}\big(1+C_s\big(\gamma_s-(1+\gamma_s\Delta C_s)\varepsilon_s^{-1}\big)\big)\vert Y_{s-}\vert^2\dC s\bigg],
\; \norm{U}_{\mathbb H^2(X^\natural,v)}^2:=\mathbb E\bigg[\int_0^T\mathcal E(v)_{s}C_s\ud \mathrm{Tr}[\langle U\star \widetilde{\mu}^\natural\rangle]_s\bigg],$$

\vspace{-0.5em}
$$ \norm{N}_{\mathbb H^2(\oX^\perp,v)}^2:=\mathbb E\bigg[\int_0^T\mathcal E(v)_{s-}C_s\ud\mathrm{Tr}[\langle N\rangle]_s+\int_0^T\Delta \mathcal E(v)_s C_s\mathds{1}_{\{X^{\circ,c}=0\}\cup\{X^{\circ,d}=0\}}\ud \mathrm{Tr}[\langle N^d\rangle]_s\bigg],$$

\vspace{-0.5em}
$$\norm{Z}_{\mathbb H^2(X^\circ,v)}^2:=\mathbb E\bigg[\int_0^TC_s\big(\mathcal E(v)_{s}-\Delta\mathcal E(v)_s\mathds{1}_{\{X^{\circ,c}\neq 0\}}\big)\ud \mathrm{Tr}[\langle Z\cdot X^\circ\rangle]_s\bigg].$$

Our goal is now to use the the results of Theorems \ref{th:estimates} and \ref{th:estimates2} to obtain sufficient conditions ensuring that the map associating the quadruplet $(y,z,u,n)$ to the quadruplet $(Y,Z,U,N)$ defined as the solution to the BSDE
$$Y_t = \xi 	+ \int_t^T f(s,y_{s-},z_s,u_s)\dC s 
				- \int_t^T Z_s \ud X^\circ_s 
				-\int_t^T\int_{\R{n}}U_s(x) \widetilde\mu(\ud s,\ud x)
				- \int_t^T \ud N_s,$$
				is a contraction in the Banach space $\mathbb H^2(\oX,v)\times\mathbb H^2(X^\circ,v)\times\mathbb H^2(X^\natural,v)\times\mathbb H^2(X^\perp,v)$. 
These spaces are defined completely analogously to the $\mathbb H^2_\beta(\cdot)$ spaces on page \pageref{def-spaces}, with the requirement that the respective norms, defined above, are finite.
Given that the norm for $Z$ depends a lot on whether $X^\circ$ is purely discontinuous or not, we will distinguish between these two cases. 
The detailed analysis will be relegated to Appendix \ref{App:AuxAnalysis}.

\begin{remark}
The above norms may seem curious at first sight. 
However, we believe they are the natural ones in the current framework for the following two reasons:
\begin{itemize}
\item[$(i)$] First of all, when $C$ is bounded, $T$ is finite, and the generator is actually Lipschitz, these norms are equivalent to the usual ones considered in the BSDE literature. This result therefore subsumes earlier and simpler ones in the literature.
\item[$(ii)$] Second, we believe that the natural spaces for solutions to BSDEs should somehow be dictated by the {\it a priori} estimates that can be obtained, and the method we used here to derive them is, by any means, a simple generalization of the classical one based on It\=o's formula and classical inequalities.
\end{itemize}
\end{remark}

After these remarks, we can state our main result of this subsection.

\begin{theorem}\label{th:main2}
Let Assumptions {\rm\ref{Hone}--\ref{Hsix}} hold true. Then we can find a non--decreasing process $v$ such that the {\rm BSDE} \eqref{BSDEintegralform} has a unique solution in $\mathbb H^2(\oX,v)\times\mathbb H^2(X^\circ,v)\times\mathbb H^2(X^\natural,v)\times\mathbb H^2(X^\perp,v)$. 
\end{theorem}

\subsubsection{Comparison with the literature}\label{subsubsec:ComparLiterature}
In Subsection \ref{subsec:Compare} we have already discussed the differences between the related literature and Theorem \ref{BSDEMainTheorem}. 
In this sub--sub--section we are going to make an analogous discussion regarding the conditions for the existence and uniqueness of the solution under the framework of Subsection \ref{subsec:AlternEstim}.

\vspace{0.5em}
\begin{enumerate}[label=\raisebox{0.08em}{\textbullet}]
\item Having in mind the counterexample provided in \cite{confortola2014backward}, we would like to see how our conditions translate if we consider the BSDE \eqref{BSDE:counterexample-2}. 
To this end, we assume $X^\circ=0$ and $\theta^\circ=0$. Then, our conditions are equivalent to $p^2r_s<1,$ which is weaker than the condition extracted in case \ref{Ctwo} of Subsection \ref{subsec:Compare}. 
The reader may observe that this condition is analogous to that of \cite[Theorem 6.1]{cohen2012existence}, under of course a different framework.
Comparing also with \cite[Theorem 4.1]{bandini2015existence}, we recall that the condition in this work read here $r_s (\Delta C_s)^2<1-\varepsilon$ $\ud\mathbb P\otimes\ud C-a.e.$, for some $\varepsilon\in(0,1).$
	
	\vspace{0.5em}
\item Having in mind the BSDE (3) of \cite{cohen2012existence}, we assume that $X^\natural=0$ and $\theta^\natural=0$. Then, our conditions translate to\footnote{When $X^{\circ,c}\neq 0$ and $\theta^\natural=0$, one can verify that condition \ref{Hsix}.(ii) is weaker than \ref{Hsix}.(i) as well as than \ref{Hsix}.(iii).}
\begin{align*}
	&\hspace{2em}C_{s-}\Delta C_s \theta^\circ_s <C_s \text{ and } (\Delta C_s)^2 r_s <\min\left\{\frac{(\sqrt{C_s}-\sqrt{C_{s-}})^2}{C_s}, \frac{\Delta C_s + C_{s-}(1-\theta^\circ_sC_{s-})}{C_s} \right\},\; \ud \mathbb P\otimes \ud C-a.e.	
\end{align*}
The second condition is reminiscent of the ones found in \cite[Theorem 6.1]{cohen2012existence}, as the upper bound is also upper bounded by $1$, and can be as close as $1$, depending on the properties of $C$. The different form and the additional constraint may appear because of the method we have followed here, which slightly differs from that of \cite{cohen2012existence}. Indeed, the approach of \cite{cohen2012existence} is to apply It\=o's formula to $Y^2$, and then use Gronwall's lemma to make the stochastic exponential appear. This can be done in this order because their process $C$ is deterministic, and thus can be taken out of conditional expectations. Since in our case $C$ is random, we apply It\=o's formula to the product of $Y^2$ and $\mathcal E(v)$ immediately. Since $C$ jumps, this creates additional cross--variation terms that need to be controlled as well, and which are the ones worsening the estimates. Obviously, in the case where $C$ is deterministic, the method of \cite{cohen2012existence} could be readily applied and would lead us to similar results.
\end{enumerate}

\subsubsection{A comparison theorem in dimension $1$}

The comparison theorem has always been recognized as a powerful tool in BSDEs analysis. 
In this sub--sub--section, we specialize the discussion to the one--dimensional case, that is to say $d=1$. 
	
\vspace{0.5em}	
We will need to work under the following assumptions.
\begin{enumerate}[label=\bf{\text{(Comp\arabic*)}},leftmargin=*,itemindent=0.0cm]
\item\label{Comp1} 
The martingale $X^\circ$ is continuous.
\item\label{Comp2}
The generator $f$ is such that for any $(s,y,z,u,u^\prime)\in\mathbb R_+\times \mathbb R\times\mathbb R^m\times\mathfrak H\times\mathfrak H$, there is some map $\rho\in\mathbb H^{2,\natural}$ 
with $\Delta \big( \rho\star\widetilde\mu^{\natural}\big)>-1$ on $\llbracket0 , T\rrbracket$, such that for $\ud \mathbb P \otimes \ud C-a.e.$ $(\omega,s)\in\mathbb R_+\times\Omega$, denoting $\delta u:=u-u^\prime$,
$$f\big(\omega,s,y,z, u(\cdot)\big)-f\big(\omega,s,y,z, u^\prime(\cdot)\big)
	\leq \widehat{K}_s\big( \delta u_s(\cdot) - \widehat{\delta u_s}\big)(\rho_s(\cdot)-\widehat{\rho}_s)\big)
	+ (1-\zeta_s)\Delta C_s\widehat{K}_s\big( \delta u_s(\cdot) - \widehat{\delta u_s}\big)\widehat{K}_s\big( \rho_s(\cdot)-\widehat{\rho_s}\big)\Big).$$
\item\label{Comp3}
The generator $f$ is such that it satisfies Assumption \ref{Hthree} and $r\Delta C^2<1,\; \ud \mathbb P\otimes \ud C-a.e.$
\end{enumerate}
\begin{remark}\label{rem:linearize-generator}
Let $f$ be a generator. 
Then, for every $(y,z,u), (y',z',u')\in\mathbb R\times\mathbb R^m\times \mathfrak H$ and $\ud \mathbb P\otimes\ud C-a.e.$ $(\omega,s)\in\Omega\times\mathbb R_+$ we can write 
\begin{align*}
&f\big(s,y,z,u_s(\omega;\cdot)\big) - f\big(s,y',z',u'_s(\omega;\cdot)\big)\\ 
	&\hspace{2em}= \lambda_s^{y,y',z,u}(\omega) (y-y') + \eta_s^{y',z,z',u,c}(\omega)c_s(\omega)(z-z')^\top + f\big(s,y',z',u_s(\omega;\cdot)\big) - f\big(s,y',z',u'_s(\omega;\cdot)\big),
\end{align*}
where
\begin{align*}
\lambda_s^{y,y',z,u}(\omega)
	&:= \begin{cases}
		\frac{f\big(s,y,z,u_s(\omega;\cdot)\big) - f\big(s,y',z,u_s(\omega;\cdot)\big)}{y-y'}, &\text{for }y-y'\ne 0\\
		0, &\text{otherwise}
		\end{cases}
\shortintertext{and}
\eta_s^{y',z,z',u,c}(\omega)
	&:=\begin{cases}
		\frac{f\big(s,y',z,u_s(\omega;\cdot)\big) - f\big(s,y',z',u_s(\omega;\cdot)\big)}{|(z-z') c_s(\omega)|^2}(z-z')c_s(\omega), &\text{for }(z-z')c_s(\omega)\neq0\\
		0,	&\text{otherwise}
	\end{cases}.
\end{align*}
Moreover, if $f$ satisfies Assumption \ref{Hthree}, then for every $(y,z,u), (y',z',u')\in\mathbb R\times\mathbb R^m\times \mathfrak H$ holds $|\lambda_s^{y,y',z,u}(\omega)|^2\le r_s(\omega)$ as well as $|\eta_s^{y',z,z',u,c}(\omega)|^2\le \theta^\circ_s(\omega)$ $\ud \mathbb P\otimes\ud C-a.e.$ on $\Omega\times \mathbb R_+$.
In the following, whenever no confusion may arise and in order to simplify the introduced notation, we will omit $y,y',z'z',u,u',c$ and we will simply write $\lambda$ and $\eta$, instead of $\lambda^{y,y',z,u}$ and $\eta^{y',z,z',u,c}$.
\end{remark}
\begin{theorem}
For $i=1,2$, let $(Y^i,Z^i,U^i,N^i)$ be solutions, in the sense of Definition \ref{solutiondef-newapriori} of the BSDEs with standard data $(\oX,\mathbb G, T,\xi^i,f^i,C)\footnotemark$.
\footnotetext{The stochastic Lipschitz bounds of the generator $f^1$, resp. $f^2$, will be denoted by $r^1,\theta^{1,\circ},\theta^{1,\natural}$, resp. $r^2,\theta^{2,\circ},\theta^{2,\natural}$.} 
Assume that Assumption \ref{Comp1} holds and that $f^1$ satisfies Assumptions \ref{Comp2}--\ref{Comp3}. 
If
\begin{itemize}
\item $\xi^1\leq \xi^2$, $\mathbb P-$a.s.
\item $f^1\big(s,Y_{s-}^2,Z_s^{2}, U_s^2(\cdot)\big)\leq f^2\big(s,Y_{s-}^2,Z_s^{2}, U_s^2(\cdot)\big),\; \ud \mathbb P\otimes \ud C-a.e.$,
\item the process $\mathcal E\big(\frac{\eta}{1-\lambda \Delta C}\cdot X^\circ + \rho\star\widetilde{\mu}^\natural\big)$ is a uniformly integrable martingale, where $\lambda, \eta$ are the processes associated to the generator $f^1$ by Remark \ref{rem:linearize-generator} for $(Y^1_-,Z^1,U^1(\cdot))$ and $(Y^2_-,Z^2,U^2(\cdot))$\footnote{More precisely, $\lambda=\lambda^{Y^1_{-},Y^2_{-}, Z^1, U^1}$ and $\eta=\eta^{Y^2_{-},Z^1, Z^2, U^1, c}$.} and $\rho$ comes from Assumption \ref{Comp2},
\end{itemize}
then, we have $Y_t^1\leq Y_t^2$, for any $t\in\llbracket 0,T\rrbracket$, $\mathbb P-$a.s.
\end{theorem}

\begin{proof}
Fix some non--negative predictable process $\gamma$ and define $v_\cdot:=\int_0^\cdot\gamma_s\ud C_s$. Define also
$$\delta Y_\cdot:=Y^1_\cdot-Y^2_\cdot,\; \delta Z_\cdot:=Z^1_\cdot-Z^2_\cdot,\; \delta U_\cdot:=U^1_\cdot-U^2_\cdot,\; \delta N_\cdot:=N^1_\cdot-N^2_\cdot,\; \delta \xi_\cdot:=\xi^1-\xi^2,$$
$$ 	\delta f^{1,2}_\cdot:= f^1\big(\cdot,Y_{\cdot-}^2,Z_\cdot^{2}, U_\cdot^2(\cdot)\big)-f^2\big(\cdot,Y_{\cdot-}^2,Z_\cdot^{2}, U_\cdot^2(\cdot)\big),\; 
	\delta f^{1}_\cdot:= f^1\big(\cdot,Y_{\cdot-}^1,Z_\cdot^{1}, U_\cdot^1(\cdot)\big)-f^1\big(\cdot,Y_{\cdot-}^2,Z_\cdot^{2}, U_\cdot^2(\cdot)\big).$$
Arguing similarly as in the proof of Theorem \ref{th:estimates}, we deduce from It\=o's formula\footnote{Recall that the notation $C^d$ was introduced before Theorem \ref{th:estimates}.}
\begin{align}\label{eq:comp}
\nonumber\mathcal E(v)_t\delta Y_t
=&\ \mathcal E(v)_T\delta\xi +\int_t^T\mathcal E(v)_{s-}\Big(\big(\delta f_s^1+\delta f_s^{1,2}\big)\big(1+\gamma_s\Delta C_s\big)-\gamma_s\delta Y_{s-}\Big)\ud C_s-\int_t^T\mathcal E(v)_{s-}\delta Z_s\ud X_s^\circ\\
&-\int_t^T\mathcal E(v)_{s-}\ud\delta N_s-\int_t^T\int_{\mathbb R^n}\mathcal E(v)_{s-}\delta U_s(x)\widetilde\mu^{\natural}(\ud s,\ud x)-\int_t^T \mathcal E(v)_{s-}\gamma_s\ud\big[ C^d, \delta U \star\widetilde{\mu}^{\natural} +  \delta N \big]_s.
\end{align}
Now, by Remark \ref{rem:linearize-generator} we can write
	$$\delta f_s^1=\lambda_s^{Y^1_{s-},Y^2_{s-}, Z^1_s, U^1_s} \delta Y_{s-}+\eta_s^{Y^2_{s-},Z_s^1, Z_s^2, U_s^1, c_s} c_s\delta Z_s^\top+f^1\big(s,Y_{s-}^2,Z_s^{2}, U_s^1(\cdot)\big)-f^1\big(s,Y_{s-}^2,Z_s^{2}, U_s^2(\cdot)\big),$$
where $\lambda$ is a $1-$dimensional predictable process such that $|\lambda|^2\leq  r^1$, $\ud \mathbb P\otimes \ud C-a.e.$, and $\eta$ is an $\mathbb R^m-$dimensional predictable process such that $\vert\eta\vert^2\leq  \theta^{1,\circ}$, $\ud \mathbb P\otimes \ud C-a.e.$\footnote{Observe that there exists a constant $D_m$, which depends only on the dimension of $\mathbb R^m$, such that $\eta^\top \frac{\ud \langle X^\circ\rangle}{\ud C} \eta \le D_m \theta^{1,\circ} \sum_{j=1}^m\sum_{i=1}^m \left(\frac{\ud \langle X^\circ\rangle}{\ud C}\right)^{ij}$ $\ud \mathbb P\otimes \ud C-a.e.$ 
Given the Assumption \ref{Comp3}, the process $\frac{\eta}{1-\lambda \Delta C}\cdot X^\circ$ is a well--defined (local) martingale.}

\vspace{0.5em}
At this point we choose $\gamma=\frac{\lambda}{1-\lambda \Delta C}$; our choice will be justified later. 
Define then the following measure $\mathbb Q$ with density
$$\frac{\ud \mathbb Q}{\ud \mathbb P}:=\mathcal E\bigg(\eta(1+\gamma\Delta C) \cdot X^\circ+\rho\star\widetilde\mu^{\natural}\bigg)_T=\mathcal E\bigg(\frac{\eta}{(1-\lambda\Delta C)} \cdot X^\circ+\rho\star\widetilde\mu^{\natural}\bigg)_T .$$
Since $\rho$ has been assumed to verify $\Delta \big( \rho\star\widetilde\mu^{\natural}\big)>-1$ on $\llbracket0 , T\rrbracket$, up to $\mathbb P-$indistinguishability, the stochastic exponential process remains (strictly) positive on $\llbracket 0, T\rrbracket$, as well as its c\`agl\`ad version; see \cite[Remark 15.3.1]{cohen2015stochastic}. 
In other words, the measure $\mathbb Q$ is equivalent to the measure $\mathbb P$.
Let us initially translate the stochastic integrals appearing in \eqref{eq:comp} into semimartingales under the measure $\mathbb Q.$
To this end, we will apply Girsanov's Theorem in its form \cite[Theorem 15.2.6]{cohen2015stochastic}.
For convenience define $M:=\eta(1+\gamma\Delta C) \cdot X^\circ+\rho\star\widetilde\mu^{\natural}.$ We have
\begin{align*}
(\mathcal E(v)_-\delta Z)\cdot X^\circ 
	&= \underbrace{\big((\mathcal E(v)_-\delta Z)\cdot X^\circ  - \mathcal E(M)_-^{-1}\cdot \langle (\mathcal E(v)_-\delta Z)\cdot X^\circ, \mathcal E(M)\rangle \big)}_{=:P^\circ \text{ which is a }\mathbb Q-\text{martingale}} 
		+ \mathcal E(M)_-^{-1}\cdot \langle (\mathcal E(v)_-\delta Z)\cdot X^\circ, \mathcal E(M)\rangle\\
	&= P^\circ	+ \mathcal E(M)_-^{-1}\mathcal E(v)_-\mathcal E(M)\cdot \langle \delta Z\cdot X^\circ, M\rangle = P^\circ	+ (\mathcal E(v)_-\delta Zc \eta^\top)\cdot C,
\end{align*}
and
\begin{align*}
&\mathcal E(v)_-\cdot \big((\delta U)\star \widetilde{\mu}^{\natural}\big)\\
	&= \underbrace{\big[\mathcal E(v)_-\cdot \big((\delta U)\star \widetilde{\mu}^{\natural}\big) - \mathcal E(M)_-^{-1}\cdot \langle \mathcal E(v)_-\cdot \big((\delta U)\star \widetilde{\mu}^{\natural}\big), \mathcal E(M)\rangle \big]}_{=:P^\natural \text{ which is a }\mathbb Q-\text{martingale}}+\mathcal E(M)_-^{-1}\cdot \langle \mathcal E(v)_-\cdot \big((\delta U)\star \widetilde{\mu}^{\natural}\big), \mathcal E(M)\rangle\\
	&= P^\natural
	+\mathcal E(v)_-\cdot \langle  (\delta U)\star \widetilde{\mu}^{\natural}, M\rangle\\
	&= P^\natural
	+\mathcal E(v)_-\cdot \langle  (\delta U)\star \widetilde{\mu}^{\natural}, \rho\star\widetilde{\mu}^\natural\rangle\\
	&\hspace{-0.3em}\overset{\eqref{repr:pqc-PDM}}{=} K^\natural
		+\int_0^\cdot \mathcal E(v)_{s-}\Big(\widehat{K}_s\Big( (\delta U_s(\cdot) - \widehat{\delta U_s})(\rho_s(\cdot)-\widehat{\rho}_s)\Big)+ (1-\zeta_s)\Delta C_s\widehat{K}_s\big( \delta U_s(\cdot) - \widehat{\delta U_s}\big)\widehat{K}_s\big( \rho_s(\cdot)-\widehat{\rho_s}\big) \Big)\ud C_s.
\end{align*}
For the term $\mathcal E(v)_-\cdot (\delta N)$ observe that it is a $\mathbb Q-$martingale, since 
\begin{align*}
\langle \delta N, M\rangle
	= \mathcal E(M)_-\cdot\big(\langle \delta N^c, \eta\cdot X^\circ\rangle + \langle \delta N^d, \rho\star\widetilde{\mu}^\natural\rangle\big)=0,
\end{align*}
where we have used \cite[Theorem III.6.4.b)]{jacod2003limit} and \cite[Theorem 13.3.16]{cohen2015stochastic}.
The last term of \eqref{eq:comp} can be written as
\begin{align*}
&\mathcal E(v)_{-}\gamma\cdot[C^d,(\delta U) \star\widetilde{\mu}^{\natural}+\delta N^d\big]\\
&=\underbrace{\Big[ \mathcal E(v)_{-}\gamma [C^d,(\delta U) \star\widetilde{\mu}^{\natural}+\delta N^d\big] 
	 - \mathcal E(M)_-^{-1}\cdot \langle \mathcal E(v)_{-}\gamma \cdot[C,(\delta U) \star\widetilde{\mu}^{\natural}+\delta N^d\big], \mathcal E(M)\rangle\Big]}_{=:P \text{ which is a }\mathbb Q-\text{martingale}}\\
	 &\hspace{0.9em}+\mathcal E(v)_{-}\gamma\cdot \langle[C,(\delta U) \star\widetilde{\mu}^{\natural}+\delta N^d\big], M\rangle\\
&= P +\mathcal E(v)_{-}\gamma \Delta C\cdot  \langle (\delta U) \star\widetilde{\mu}^{\natural}+\delta N^d, M^d\rangle\\
&= P +\mathcal E(v)_{-}\gamma \Delta C\cdot  \langle (\delta U) \star\widetilde{\mu}^{\natural}, \rho \star \widetilde{\mu}^\natural\rangle \\
&= P +\int_0^\cdot \mathcal E(v)_{s-}\gamma_s \Delta C_s\Big(\widehat{K}_s\Big( (\delta U_s(\cdot) - \widehat{\delta U_s})(\rho_s(\cdot)-\widehat{\rho}_s)\Big)+ (1-\zeta_s)\Delta C_s\widehat{K}_s\big( \delta U_s(\cdot) - \widehat{\delta U_s}\big)\widehat{K}_s\big( \rho_s(\cdot)-\widehat{\rho_s}\big) \Big)\ud C_s.
\end{align*}

Using Girsanov's theorem, \emph{i.e.} the above $\mathbb Q-$canonical decompositions, as well as the assumption on $\delta \xi$ and $\delta f^{1,2}$,we deduce from \eqref{eq:comp} after applying the $\mathbb Q-$conditional expectation with respect to $\mathcal G_t$
\begin{align*}
\mathcal E(v)_t\delta Y_t 
	\le&\  \mathbb E^{\mathbb Q}\bigg[\int_t^T\mathcal E(v)_{s-} \Big( \big[\lambda_s(1+\gamma_s\Delta C_s) - \gamma_s\big] \delta Y_{s-} \Big)\ud C_s\bigg|\mathcal G_t\bigg]\\ 
	& +\mathbb E^{\mathbb Q}\bigg[\int_t^T  \mathcal E(v)_{s-}(1+\gamma_s \Delta C_s)\Big(f^1\big(s,Y_{s-}^2,Z_s^{2}, U_s^1(\cdot)\big)-f^1\big(s,Y_{s-}^2,Z_s^{2}, U_s^2(\cdot)\big)\\
	&\hspace{6em}-\widehat{K}_s\big( (\delta U_s(\cdot)- \widehat{\delta U_s})(\rho_s(\cdot)-\widehat{\rho}_s)\big) - (1-\zeta_s)\Delta C_s\widehat{K}_s\big( \delta U_s(\cdot) - \widehat{\delta U_s}\big)\widehat{K}_s\big( \rho_s(\cdot)-\widehat{\rho_s}\big)\Big)
				   \ud C_s\bigg|\mathcal G_t\bigg].
\end{align*}
Recall our choice for $\gamma$, \emph{i.e.} $\gamma=\frac{\lambda}{1-\lambda \Delta C}$, and we have that the first conditional expectation on the right--hand side vanishes. 
In view of Assumption \ref{Ctwo}, we can conclude if the stochastic exponential $\mathcal E(v)$ (as well as the process $\mathcal E(v)_-$) remains strictly positive, which is true if and only if
$$\Delta v_s>-1\Longleftrightarrow (\lambda_s\Delta C_s)^2<1,$$
which is automatically satisfied since we assumed that $r_s^1\Delta C_s^2<1$. 
\end{proof}

\begin{remark}
The Assumption \ref{Comp2} is the natural generalization of the Assumption \textbf{(A$_{\mathbf{\gamma}}$)} in \cite{royer2006backward}, where we have abstained from assuming that the predictable function $\gamma$ (we follow at this point the notation of \cite{royer2006backward}) may depend on $y,z,u,u'$.
In order to verify our statement, let us assume that the martingale $X^\natural$ exhibits jumps only on totally inaccessible times\footnote{In \cite{royer2006backward} the filtration is quasi--left--continuous, which means that the uniformly integrable purely discontinuous martingale jumps only on totally inaccessible times.}.
In this case, by \cite[Corollary 13.3.17]{cohen2015stochastic} and the polarization identity, we have 
$$\langle (\delta U)\star \widetilde{\mu}^\natural, \rho \star \widetilde{\mu}^\natural\rangle_\cdot = \big((\delta U)\rho\big)\star\nu^\natural_\cdot = \Delta C_\cdot\widehat{K}_\cdot\big((\delta U_\cdot(\cdot))\rho_\cdot(\cdot)\big).$$
In other words, Assumption \ref{Comp2} can be simplified to: for\ $\ud\mathbb P\otimes \ud C-a.e.\ (\omega,s)$ holds
$$f\big(\omega,s,y,z, u(\cdot)\big)-f\big(\omega,s,y,z, u^\prime(\cdot)\big)\leq \widehat{K}_s\big( \delta u_s(\cdot)\rho_s(\cdot)\big),$$
for any $(s,y,z,u,u^\prime)\in\mathbb R_+\times \mathbb R\times\mathbb R^m\times\mathfrak H\times\mathfrak H$.
\end{remark}

\begin{remark}
Recall that $|\eta|^2\le \theta^{1,\circ}$ $\ud \mathbb P\otimes \ud C-a.e$. 
Observe, now, that there exists a constant $D_m$, which depends only on the dimension of $\mathbb R^m$, 
such that 
$$\eta^\top \frac{\ud \langle X^\circ\rangle}{\ud C} \eta \le D_m \theta^{1,\circ} \sum_{j=1}^m\sum_{i=1}^m \left(\frac{\ud \langle X^\circ\rangle}{\ud C}\right)^{ij},\; \ud \mathbb P\otimes \ud C-a.e.$$ 
Therefore, given Assumption \ref{Comp3}, a sufficient condition for the process $\eta/(1-\lambda \Delta C)\cdot X^\circ$ to be well--defined is the existence of the 
(local) martingale $\sqrt{\theta^{1,\circ}}/(1-\sqrt{r^1} \Delta C)\cdot X^\circ$. 
For the comment which comes after the following lines, the reader may recall some criteria which guarantee the (true) martingale property of a stochastic exponential and which involve an integrability condition on the predictable covariation of the continuous part of the density 
(in our case $\langle \eta/(1-\lambda \Delta C)\cdot X^\circ \rangle$), for instance \cite[Theorem 15.4.2 -- Theorem 15.4.6]{cohen2015stochastic}.
Having these criteria in mind, one may check if the process $\mathcal E\big(\frac{\sqrt{\theta^{1,\circ}}}{1-\sqrt{r^1} \Delta C}\cdot X^\circ + \rho\star\widetilde{\mu}^\natural\big)$ satisfies any of them. 
If the answer is affirmative, then we have a sufficient condition for the process $\mathcal E\big(\frac{\eta}{1-\lambda \Delta C}\cdot X^\circ + \rho\star\widetilde{\mu}^\natural\big)$ to be a uniformly integrable martingale for any choice of $\lambda$ and $\eta$.
\end{remark}

%% file: Applications.tex

As an application of the main theorem, we show that a BSDE driven by an extended Grigelionis process, which is, roughly speaking, a superposition of a time--inhomogeneous L\'evy process with a (discrete-time) random walk, admits a unique solution under appropriate conditions. 
The main point here is that when $C$ is allowed to have jumps, there is a subtle interplay between the size of the jumps of $C$ and the strength of the dependence of the generator of the BSDE, measured by the value of the Lipschitz coefficients, in the sense that their product has to remain small.

\begin{definition}
A square--integrable $\mathbb R^m-$valued martingale $X$ is called \emph{$K-$almost} \emph{quasi--left--conti\-nuous} if there exists a constant $K\geq 0$ such that $|\Delta\langle X^{i,j}\rangle|_t\leq K$ for every $i,j=1,\dots,\ell$ and for every $t\in\mathbb R_+,$ $\Pm-a.s.$
In other words, the predictable quadratic covariation $\langle X\rangle$ of $X$ has jumps uniformly bounded by $K$.
\end{definition}

The next result follows directly from the definition above and Theorems \ref{BSDEMainTheorem} and \ref{th:main2}.

\begin{corollary}
Let $(X, \Fil, T, \xi, f, C)$ be standard data under $\hat{\beta}$, $X$ be $K-$almost quasi--left--continuous and the process $\alpha^2$ $($defined in \ref{Ffour}$)$ be bounded by $1/(18 {\rm e} m K)$, $\mathbb P-a.s.$, where $m$ is the dimension of $X.$ 	
Then, for $C={\rm Tr}[\langle X\rangle]$ and for $\hat{\beta}$ large enough, there exists a unique solution $(Y,Z,U,N)$ to the BSDE \eqref{BSDE}. Similarly, if Conditions \ref{Hone}--\ref{Hsix} are satisfied, then there is a unique solution in the sense of Definition \ref{solutiondef-newapriori}.

\end{corollary}

\begin{example}
Let $(X, \Fil, T, \xi, f, C)$ be standard data under $\hat{\beta}$ such that $X=\lambda G$, for some $\lambda\in\mathbb R$ and some \emph{extended Grigelionis} martingale $G$. 
In other words, $C$ can be chosen to be of the form
\begin{align*}
C_t=\lambda^2\Big(t + \sum_{s\leq t}{\mathds{1}_{\Theta}(s)}\Big),
\end{align*}
where $\Theta\subset(0,+\infty)$ is at most countable, see {\rm \cite[Definition 2.15]{kallsen1998semimartingale}}.
Then, since $X$ is $\lambda^2-$almost quasi--left--continuous, for $\alpha^2$ bounded by $1/(18{\rm e}\lambda^2)$ and for $\hat\beta$ large enough, there exists a unique solution to the BSDE \eqref{BSDE}. Similarly, if Conditions \ref{Hone}--\ref{Hsix} are satisfied, then there is a unique solution in the sense of Definition \ref{solutiondef-newapriori}.
\end{example}

Another interesting application of this result consists in ensuring the existence and uniqueness of the solution of the BSDE \eqref{BSDE} when $X$ is the continuous--time extension of a discrete time martingale $\hat{X}$. In particular, when $\hat{X}^n$ is the discretization of a square integrable, quasi--left continuous martingale with independent increments. Then, as the mesh of the grid tends to $0$, the bound $K^n$ of the jumps of $\langle\hat{X}^n\rangle$ tends to $0$. Hence, we have that the sequence of BSDEs is, for $n$ large enough, well--posed, given that the associated $\alpha^2$ is bounded $\mathbb P-a.s.$ 
We emphasize again that in general, Theorem \ref{th:main2} cannot be applied in this setting, recall Remark \ref{rem:H6}, and the only result available in the literature in such a general framework is our Theorem \ref{BSDEMainTheorem}.

%% file: Auxilliary_Results_Proof_of_Propositions_Preliminaries.tex

Let us fix a pair $(X^\circ,X^\natural)\in\mathcal H^2(\mathbb R^m)\times\mathcal H^{2,d}(\mathbb R^n)$ such that $M_{\mu^{\natural}}[\Delta X^\circ|\widetilde{\mathcal P}]=0.$
We will adopt the notation which follows \ref{Hsix}.
Moreover, $\mathbb F^Y$ denotes the natural filtration of the process $Y$.

\begin{proof}[Proof of Lemma \ref{lem:StochIntOrthog}]
For every $Y^\natural\in\mathcal K^2(\mu^{X^\natural})$ holds
\begin{equation}\label{lem:widebarX-Orthog-Parts1}
	\langle X^{\circ},Y^\natural\rangle = \langle X^{\circ,c},Y^\natural\rangle + \langle X^{\circ,d},Y^\natural\rangle = 0,
\end{equation}
where the summand $\langle X^{\circ,c},Y^\natural\rangle$ vanishes because $X^{\circ,c}\in\mathcal H^{2,c}(\mathbb R^m)$ and $Y^{\natural}\in\mathcal H^{2,d}(\mathbb R^d)$, 
while the summand $\langle X^{\circ,d},Y^\natural\rangle$ equals $0$ in view of the assumption  $M_{\mu^{X^{\natural}}}[\Delta X^{\circ}|\widetilde{\mathcal P}]=0$, of \cite[Theorem III.4.20]{jacod2003limit} and of \cite[Theorem 13.3.16]{cohen2015stochastic}.
The equality \eqref{lem:widebarX-Orthog-Parts1} already proves the second statement.
Indeed, for $j=1,\ldots,n$ define $\mathbb R^n\ni x\overset{\pi_j}{\longmapsto} (x^j,0,\ldots,0)\in\mathbb R^d\footnotemark$. \footnotetext{ $x^j$ is the $j-$component of the vector $x.$ In other words, $\pi^j$ behaves as the canonical $j-$projection.}
Then,  we have that $\pi_j\in\mathbb H^2(\mu^{X^\natural})$ and $X^{\natural,j}= \pi_j\star\widetilde{\mu}^{X^{\natural}}$, since $X^\natural\in\mathcal H^{2,d}(\mathbb R^n).$

\vspace{0.2em}
Assume now the factorization 
\begin{equation*}
	\langle X^\circ\rangle = \int_{(0,\cdot]}\frac{\textup{d}\langle X^\circ\rangle_s}{\textup{d}F_s}\textup{d}F_s.
\end{equation*}
In view of \eqref{lem:widebarX-Orthog-Parts1} we obtain\footnote{One can follow similar arguments to those following \cite[Statement III.4.3]{jacod2003limit}.} for the predictable process $r^{Y^\natural} := \frac{\textup{d}\langle Y^\natural, X^{\circ}\rangle_s}{\textup{d}F_s} =0.$
Consequently, by \cite[Theorem III.6.4.b)]{jacod2003limit} we have for every $Z\in\mathbb H^2(X^\circ)$ that
\begin{equation*}
	\langle Y^\natural, Z\cdot X^\circ\rangle = \int_{(0,\cdot]}r^{Y^\natural}Z^\top\,\textup{d}F_s=0, 
\end{equation*}
where the equality is understood componentwise.
\end{proof}
\begin{proof}[Proof of Proposition \ref{prop:OrthogDecomp}]
By the Galtchouk--Kunita--Watanabe Decomposition, see \cite[Chapitre IV, Section 2]{jacod1979calcul}, there exists $Z\in\mathbb H^2(X^\circ)$ such that 
\begin{equation}\label{equal:orthog-decomp-GKW}
Y-Y_0=Z\cdot X^\circ + \overline{N} 
\end{equation}
with $\overline{N}\in\mathcal H^2(\mathbb R^d)$ and $\langle X^{\circ}, \overline{N}\rangle=0$.
Moreover, by \cite[Theorem III.4.20]{jacod2003limit}, there exists a unique $U\in\mathbb H^2(\mu^{X^\natural})$ and $N\in\mathcal H^2(\mathbb R^d)$ with $M_{\mu^{X^{\natural}}}[\Delta N|\widetilde{\mathcal P}]=0$ such that 
\begin{equation}\label{equal:orthog-decomp-2}
\overline{N} = U\star \widetilde{\mu}^{(X^\natural,\mathbb G)} + N.
\end{equation}
In total, we have determined $Z\in\mathbb H^2(X^\circ)$, $U\in\mathbb H^2(\mu^{X^\natural})$ and $N\in\mathcal H^2(\mathbb R^d)$ such that 
\begin{equation}\label{Orthog-Decomp-Y}
Y=Y_0+Z\cdot X^\circ + U\star\widetilde{\mu}^{(X^\natural,\mathbb G)}+  N.
\end{equation}

We have to verify that this decomposition satisfies the properties of Definition \ref{def:OrthogDecomp} and, moreover, that it does not depend on the way we have determined $Z,U$ and $N$.

\vspace{0.5em}
We will prove initially that the $\mathbb G-$predictable function $U$ is the one characterized by the triplet $(\mathbb G,\mu^{X^\natural},Y)$.
To this end, we are going to prove that $M_{\mu^{X^\natural}}\big[\Delta (Z\cdot X^\circ)\big|\widetilde{\mathcal P}\big]=0$.
By \cite[Proposition III.6.9]{jacod2003limit}, we can write 
$$\Delta (Z\cdot X^\circ)^i = \sum_{j=1}^m Z^{ij} \Delta X^{\circ,j}.$$
Therefore, for every positive and bounded $\mathbb G-$predictable function $W$ holds for every $i=1,\ldots,d$
\begin{align*}
M_{\mu^{X^\natural}}\big[ W \Delta (Z\cdot X^\circ)^i\big] 
	= M_{\mu^{X^\natural}}\Big[ W \sum_{j=1}^mZ^{ij} \Delta X^{\circ,j} \Big]
	= \sum_{j=1}^m M_{\mu^{X^\natural}}\Big[ W  Z M_{\mu^{X^\natural}}\big[\Delta X^{\circ,j}\big|\widetilde{\mathcal P}\big]\Big]  
	= 0,
\end{align*}
where, in order to conclude, we used the assumption $M_{\mu^{X^{\natural}}}[\Delta X^\circ|\widetilde{\mathcal P}]=0$ and that $W$, resp. $Z$, is a $\mathbb G-$predictable function, resp. $\mathbb G-$predictable process.
The above, after using standard monotone class arguments, allows us to conclude the required property $M_{\mu^{X^\natural}}\big[\Delta (Z\cdot X^\circ)\big|\widetilde{\mathcal P}\big]=0$.
By Equality \eqref{Orthog-Decomp-Y}, the equality
$M_{\mu^{X^{\natural}}}[\Delta (Z\cdot X^\circ)|\widetilde{\mathcal P}]=0$
 and the linearity of the Dol\'eans-Dade measure $M_{\mu^\natural}$ we obtain that the following hold $M_{\mu^{X^\natural}}-$almost everywhere
\begin{equation*}
M_{\mu^{X^\natural}}\big[\Delta Y\big|\widetilde{\mathcal P}\big]
	=M_{\mu^{X^\natural}}\big[\Delta \overline{N}\big|\widetilde{\mathcal P}\big]
	=M_{\mu^{X^\natural}}\big[\Delta (U\star\widetilde{\mu}^{(X^{\natural},\mathbb G)})\big|\widetilde{\mathcal P}\big].
\end{equation*}
Hence the $\mathbb G-$predictable function $U$ is uniquely determined $M_{\mu^{X^\natural}}-$almost everywhere, see \cite[Theorem III.4.20]{jacod2003limit} and \cite[Lemma III.4.19]{jacod2003limit}.

\vspace{0.5em}
We need to prove now that $\langle Z\cdot X^\circ, U\star\widetilde{\mu}^{X^\natural}\rangle =0$ as well as $\langle N,X^{\circ}\rangle = 0$.
But the former is immediate by Lemma \ref{lem:StochIntOrthog}. 
We proceed to prove the  $\langle N,X^{\circ}\rangle = 0$.
By the Galtchouk--Kunita--Watanabe decomposition \eqref{equal:orthog-decomp-GKW} and the orthogonality of the stochastic integrals, we obtain
\begin{equation*}
\langle N, X^{\circ}\rangle \overset{\text{\eqref{equal:orthog-decomp-2}}}{=} \langle\overline{N}, X^{\circ}\rangle - \langle U\star\widetilde{\mu}^{X^\natural}, X^{\circ}\rangle = 0. 
\numberthis\label{equal:orthog-decomp-4}
\end{equation*}

\vspace{0.5em}
To sum up,
\begin{enumerate}[label=(\roman*)]
	\item $Z\in\mathbb H^2(X^\circ)$ and $U\in\mathbb H^2(\mu^{X^\natural})$, by Decompositions \eqref{equal:orthog-decomp-GKW} and \eqref{equal:orthog-decomp-2}. 
	Moreover, $Z\cdot X^\circ$ and $U\star\widetilde{\mu}^{(X^\natural,\mathbb G)}$ are unique up to indistinguishability.
	Therefore, also $N$ is unique up to indistinguishability.
	\item $\langle Z\cdot X^\circ , U\star \widetilde{\mu}^{(X^\natural,\mathbb G)}\rangle = 0$ and
	\item $\langle N, X^{\circ}\rangle = 0$ and $M_{\mu^{X^{\natural}}}[\Delta N|\widetilde{\mathcal P}]=0$.\qedhere
\end{enumerate}
\end{proof}

\begin{proof}[Proof of Proposition \ref{prop:CharacterOrthogSpace}]
Let us define 
\begin{align*}
\mathcal N:=\big\{ L\in\mathcal H^2(\mathbb R^d), \langle X^\circ,L\rangle=0 \text{ and } M_{\mu^{X^\natural}}[\Delta L|\widetilde{\mathcal P}]=0\big\}.
\end{align*}
It is immediate that $\mathcal N$ is a linear subspace of $\mathcal H^2(\mathbb R^d)$.
By the properties of the stochastic integrals, see \cite[Theorem III.6.4]{jacod2003limit} and \cite[Theorem 13.3.16]{cohen2015stochastic}, we have that $\mathcal N\subset \mathcal H^2(\oX^\perp)$.
Therefore, we need to prove the inverse inclusion.

\vspace{0.2em}
Let $L\in\mathcal H^2(\oX^\perp)$. 
Then, we have immediately  $\langle X^{\circ},L\rangle=0$.
We need, now, to prove that $M_{\mu^{X^\natural}}[\Delta L|\widetilde{\mathcal P}]=0.$
By Proposition \ref{prop:OrthogDecomp} and due to the fact that $\langle X^{\circ},L\rangle=0$ we can assume that
\begin{align*}
L &= W\star\widetilde{\mu}^{X^\natural} + N,
\numberthis\label{equal:test}
\end{align*}
where $W\in\mathbb H^2(\mu^{X^\natural})$ and $N\in\mathcal H^{2}(\mathbb R^d)$ be such that $M_{\mu^{X^\natural}}[\Delta N|\widetilde{\mathcal P}]=0.$
The last property implies that $$\langle W\star \widetilde{\mu}^{X^\natural},N^d\rangle = \langle W\star \widetilde{\mu}^{X^\natural},N\rangle=0;$$ see \cite[Theorem 13.3.16]{cohen2015stochastic}.
On the other hand, since $\langle L, U\star\widetilde{\mu}^{X^\natural}\rangle=0$ for every $U\in\mathbb H^2(\mu^{X^\natural})$, we have    
\begin{align*}
0	= \langle L,W\star \widetilde{\mu}^{X^\natural}\rangle 
	= \langle L^d,W\star \widetilde{\mu}^{X^\natural}\rangle 
\overset{\eqref{equal:test}}{=} 
	\langle W\star \widetilde{\mu}^{X^\natural},W\star \widetilde{\mu}^{X^\natural}\rangle + \langle N^d,W\star \widetilde{\mu}^{X^\natural}\rangle 	
	= \langle W\star \widetilde{\mu}^{X^\natural},W\star \widetilde{\mu}^{X^\natural}\rangle.
\end{align*}
Therefore, the stochastic integral $W\star\widetilde{\mu}^{X^\natural}$ in Decomposition \eqref{equal:test} is indistinguishable from the zero martingale.
This further implies that $L^d$ is indistinguishable from $N^d$ or in other words the processes $\Delta L$ and $\Delta N$ are indistinguishable.
Therefore, for $L$ holds also $M_{\mu^{X^\natural}}[\Delta L|\widetilde{\mathcal P}]=0,$ which proves that $L\in \mathcal N.$ 

\vspace{0.5em}
It is only left to prove that the space $\big(\mathcal H^2(\oX^\perp),\Vert\cdot\Vert_{\mathcal H^2(\mathbb R^d)}\big)$ is closed.
To this end, assume the sequence $(L^k)_{k\in\mathbb N}\subset \mathcal N$ and $L^\infty\in\mathcal H^2(\mathbb R^d)$ be such that $L^k\xrightarrow[k\to\infty]{\hspace{0.1em} \Vert\cdot\Vert_{\mathcal H^2(\mathbb R^d)}\hspace{0.1em}} L^\infty\footnote{It is well-known that $\mathcal H^2(\mathbb R^d)$ is closed.}.$
By Proposition \ref{prop:OrthogDecomp} we have that there exist $Z^\infty\in\mathbb H^2(X^\circ)$, $U^\infty\in\mathbb H^2(\mu^{X^\natural})$ and $N^\infty\in\mathcal H^2(\mathbb R^d)$ such that 
\begin{align}\label{decomp:L_infty}
L^\infty = Z^\infty \cdot X^\circ + U^\infty \star \widetilde{\mu}^{X^\natural} + N^\infty,
\end{align}
where $\langle X^\circ, N^\infty\rangle= 0$ and $M_{\mu^{X^\natural}}[\Delta N^\infty |\widetilde{\mathcal P}]=0.$

\vspace{0.2em}

By Vitali's Convergence Theorem, we have that the sequence $\big(\textrm{Tr}([L^k]_\infty\big)_{k\in\mathbb N\cup\{\infty\}}$ is uniformly integrable and, consequently, $\big(\textrm{Tr}([L^k]_\infty\big)_{k\in\mathbb N\cup\{\infty\}}\cup\{\textrm{Tr}([ X^\circ]_\infty)\}$ as well.
Therefore, by the Kunita--Watanabe and Young inequalities we have that the family $\big\{ \textrm{Var}([L^k, X^\circ])_\infty \big\}_{k\in\mathbb N\cup\{\infty\}}$ is uniformly integrable. 
For more details in the last argument, one can consult \cite[Lemma A.2]{saplaouras2017backward}.  
Consequently, from \cite[Theorem VI.6.26]{jacod2003limit} we have for every $A\in\mathcal G_\infty$ that 
\begin{align*}
\big([L^k, X^\circ], \mathbb E[\mathds{1}_A|\mathcal G_\cdot]\big) \xrightarrow[k\to\infty]{} \big([L^\infty, X^\circ], \mathbb E[\mathds{1}_A|\mathcal G_\cdot]\big)
\end{align*}
and by \cite[Proposition IX.1.12]{jacod2003limit} that $[L^\infty,  X^\circ]$ is a uniformly integrable martingale with respect to the natural filtration of $Y^A:=(L^\infty,X^\circ, \mathbb E[\mathds{1}_A|\mathcal G_\cdot])$; the latter is denoted by $\mathbb F^{Y^A}$. 
It is easy to verify now that $[L^\infty,X^\circ]$ is a uniformly integrable $\mathbb G-$martingale.
Indeed, fix $0\le s<t$. 
For every $A\in\mathcal G_s$ holds
\begin{align*}
	\int_A \mathbb E\big[ [L^\infty,X^\circ]_t \big|\mathcal G_s\big] \ud \mathbb P
	=\int_A \mathbb E\Big[\mathbb E\big[ [L^\infty,X^\circ]_t \big|\mathcal G_s\big]\Big| \mathcal F^{Y^A}_s\Big] \ud \mathbb P
	=\int_A \mathbb E\big[ [L^\infty,X^\circ]_t \big|\mathcal F^{Y^A}_s\big] \ud \mathbb P
	=\int_A  [L^\infty,X^\circ]_s \ud \mathbb P.
\end{align*} 
Therefore, $\langle L^\infty, X^\circ\rangle$ is well-defined and, in view of the previous information, it is equal to $0.$
The properties of the It\=o Integral allow us to conclude that $\langle L^\infty, Z\cdot X^\circ\rangle =0$ for every $Z\in\mathbb H^2(X^\circ)$.
Following similar arguments, we can prove that $\langle L^\infty, U\star\widetilde{\mu} \rangle=0$ for every $U\in\mathbb H^2(\mu^{X^\natural})$.
\end{proof}

\begin{proof}[Proof of Corollary \ref{cor:SpaceDecomposition}]
In view of the previous results we need only to justify the closedness of the spaces $\mathcal L^2(X^\circ)$ and $\mathcal K^2(\mu^{X^\natural})$.
For the former see \cite[Theorem III.6.26]{jacod2003limit} and use the fact that the topology induced by $\Vert \cdot\Vert_{\mathcal H^2(\mathbb R^d)}$ is stronger than the Emery topology.
For the latter, we can follow arguments analogous to the proof of Proposition \ref{prop:OrthogDecomp}. 

Let $(U^k)_{k\in\mathbb N}\subset \mathbb H^2(\mu^{X^\natural})$ such that $U^k\star \widetilde{\mu}^{X^\natural}\xrightarrow[\hspace{1em}k\to\infty\hspace{1em}]{\mathcal H^2(\mathbb R^d)}L^\infty$, for some $L^\infty\in\mathcal H^2(\mathbb R^d)$ with orthogonal decomposition
\begin{align*}
L^\infty = Z^\infty\cdot X^\circ + U^\infty\star\widetilde{\mu}^{X^\natural} + N^\infty,
\end{align*} 
where $Z^\infty\in\mathbb H^2(X^\circ)$, $U^\infty \in\mathbb H^2(\mu^{X^\natural})$ and $N^\infty\in\mathcal H^2(\oX^\perp)$.
We need to prove that the martingales $Z^\infty\cdot X^\circ$ and $N^\infty$ are indistinguishable from the zero process.
Using the convergence 
\begin{align*}
\big([U^k\star \widetilde{\mu}^{X^\natural}, X^\circ], \mathbb E[\mathds{1}_A|\mathcal G_\cdot]\big) \xrightarrow[k\to\infty]{} \big([L^\infty, X^\circ], \mathbb E[\mathds{1}_A|\mathcal G_\cdot]\big)
\shortintertext{and}
\big([U^k\star \widetilde{\mu}^{X^\natural}, L], \mathbb E[\mathds{1}_A|\mathcal G_\cdot]\big) \xrightarrow[k\to\infty]{} \big([L^\infty, L], \mathbb E[\mathds{1}_A|\mathcal G_\cdot]\big),
\end{align*}
which are true for every $A\in\mathcal G_\infty$ and $L\in\mathcal H^2(\oX^\perp)$, we can obtain the martingale property of $[L^\infty, X^\circ]$ and of the elements of the family $([L^\infty,L])_{L\in\mathcal H^2(\oX^\perp)}.$
In particular, we obtain that $\langle Z^\infty\cdot X^\circ\rangle =0$, by making use of the usual properties of the It\=o integral, and $\langle N^\infty\rangle=0,$ which provide the required result.
\end{proof}

%% file: Auxilliary_Results_Proof_of_Lemma_of_generalised_inverses.tex

The proof of Lemma \ref{mainlemma} heavily relies on Lemma \ref{GeneralisedInverseLemma}.\ref{GILemmaSeven}, which is a complement result to the already known ones about generalized inverses.

\begin{lemma}\label{GeneralisedInverseLemma}
	Let $A:\Rp\longrightarrow\Rp$ be a c\`adl\`ag and increasing function with $A_0=0$. Denote by $L:\Rp\to\Rp\cup\{\infty\}$ the  c\`agl\`ad generalized inverse of $A$, i.e.
	\[
		L(s):=\inf \left\{ t\in\Rp, \ A(t)\geq s\right\}
	\]
	and by
	$R:\Rp\to\Rp\cup\{\infty\}$ the  c\`adl\`ag generalized inverse of $A$, i.e.
	\[
		R(s):=\inf \left\{ t\in\Rp, \ A(t)> s\right\}.
	\]
We have
	\begin{enumerate}[label={\rm(\roman*)},leftmargin=*]
		\item 	$L,R$ are increasing.
		\item  	$L(s)=R(s-)$ and $L(s+)=R(s)$.
		\item  	$s\leq A(t)$ if and only if $L(s)\leq t$ and $s<A(t)$ if and only if $R(s)<t$.
		\item  	$A(t)< s$ if and only if $t<L(s)$ and $A(t)\leq s$ if and only if $t\leq R(s)$.
		\item\label{GILemmaFive}  
				$A\left( R(s)\right)\geq A\left(L(s)\right) \geq s$, for $s\in \Rp$,
			   	and at most one of the inequalities can be strict.
		\item\label{GILemmaSix}  
				For $s\in A(\mathbb R_+)$, $A(L(s))=s$.
		\item\label{GILemmaSeven}
			   For $s$ such that $ L(s)<\infty$, we have 
			   $$s\leq A\left(L(s)\right) \leq s + \Delta A(L(s)),$$
			  	where $\Delta A(L(s))$ is the jump of the function $A$ at the point $L(s)$. 
	\end{enumerate}
\end{lemma}

\begin{proof}
	Here we will prove only inequality \ref{GILemmaSeven}. However, we present the other properties, since we will make 
	use of them in the proof. The interested reader can find their proofs in a slightly more general 
	framework in \cite{embrechts2013note} and the references therein.
	
	\vspace{0.5em}
	We need to prove that $A(L(s)) - s\leq \Delta A( L(s))$ for any $s$ such that $L(s)<\infty$.	By \ref{GILemmaSix}, when $s\in A(\mathbb R_+)$, we have since $A$ is increasing
	$$A(L(s))-s =0\leq \Delta A( L(s)).$$
	
	Now if $s\notin A(\mathbb R_+)$ and $s>A_\infty$, then $L(s)=\infty$, so that this case is automatically excluded. Therefore, we now assume that $s\notin A(\mathbb R_+)$ and $s\leq A_\infty$. Since $s\notin A(\mathbb R_+)$, there exists some $t\in\mathbb R_+$ such that $s\in[A(t-),A(t))$. Then, we immediately have $L(s)=t$. Hence
	$$s + \Delta A(L(s))=s+\Delta A(t)\geq A(t)=A(L(s)),$$
	since $s\geq A(t-)$.
\end{proof}

\begin{proof}[Proof of Lemma \ref{mainlemma}]
	Using a change of variables, Lemma \ref{GeneralisedInverseLemma}.\ref{GILemmaSeven} and that $g$ is non-decreasing and sub-multiplicative, we have
	\begin{align*}
		\int_0^t g(A_s) dA_s &= \int_{A_0}^{A_t} g(A_{L_s}) ds\leq \int_{A_0}^{A_t} g(s+\Delta A_{L_s}) ds\\
		&\leq \int_{A_0}^{A_t} g\Big(s+\max_{\set{s,\ L_s<\infty}}\Delta A_{L_s}\Big) ds
		 \leq c g\Big(\max_{\set{s,\ L_s<\infty}}\Delta A_{L_s}\Big)\int_{A_0}^{A_t} g(s) ds.
	\end{align*}
\end{proof}

%% file: Auxilliary_Results_Proof_of_Lemma_of_constants.tex

\begin{proof}
Let $(\gamma,\delta)\in\mathcal C_\beta.$
We shall begin by obtaining the critical points of the map $\Pi^\Phi.$
We have 
\begin{align*}
\frac{\partial}{\partial \gamma}\Pi^\Phi(\gamma,\delta) 
	& = (2+9\delta) \, \e{(\delta-\gamma)\Phi} \, \frac{\Phi\gamma^2 + (2-\delta\Phi)\gamma - \delta}{\gamma^2(\delta-\gamma)^2}, \\
\frac{\partial}{\partial \delta}\Pi^\Phi(\gamma,\delta) 
	& = -\frac{9}{\delta^2} + \e{(\delta-\gamma)\Phi} \left\{\frac{\big[ 9 + (2+9\delta)\Phi \big](\delta-\gamma)}{\gamma(\delta-\gamma)^2} 
	    - \frac{2+9\delta}{\gamma(\delta-\gamma)^2} \right\}.
\end{align*}
The only possible critical points for $\Pi^\Phi$ are therefore such that $\delta=-2/9$ or $\gamma=\frac{\delta\Phi - 2\pm\sqrt{4+\delta^2\Phi^2}}{2\Phi}$.
However, the values $\delta=-2/9$ and $\gamma=\frac{\delta\Phi - 2-\sqrt{4+\delta^2\Phi^2}}{2\Phi}$ are ruled out as negative.
For $0<\delta\le\beta$ we have 
$$\Big(\frac{\delta\Phi - 2+\sqrt{4+\delta^2\Phi^2}}{2\Phi}, \delta\Big)\in\mathcal C_\beta.$$
Let us define $\overline{\gamma}^\Phi(\delta):=\frac{\delta\Phi - 2+\sqrt{4+\delta^2\Phi^2}}{2\Phi}$, for $0<\delta\le \beta.$
It is easy to verify that $\overline{\gamma}^\Phi(\delta)\in(0,\delta).$
Then, some tedious calculations yield that
\begin{align*}
\frac{\partial \Pi^\Phi}{\partial \delta}\big(\bar\gamma^\Phi(\delta),\delta\big)
	&= -\frac{9}{\delta^2} 
	- \frac{\exp[\big(\delta-\bar\gamma^\Phi(\delta)\big)\Phi]}{\bar\gamma^\Phi(\delta)\big(\delta-\bar\gamma^\Phi(\delta)\big)^2}\cdot\frac{2\bar\gamma^\Phi(\delta)\Phi + 9\bar\gamma^\Phi(\delta) + 2}{(\bar\gamma^\Phi(\delta)\Phi+1)^2}<0
\end{align*}
therefore $\Pi^\Phi$ does not admit any critical point on $\mathcal C_\beta$, for which $0<\gamma<\delta<\beta.$
Hence, the infimum on this set is necessarily attained on its boundary. 
The cases where at least one among $\delta$ and $\gamma$ goes to $0$, or where their difference goes to $0$, lead to the value $\infty$. 
The only remaining case is therefore $0<\gamma<\delta=\beta$, where $\beta$ is fixed. 
Then we get 
\begin{align*}
\frac{\ud}{\ud \gamma} \Pi^\Phi\big(\gamma,\beta\big)
	&= (2+9\beta) \, \e{(\beta-\gamma)\Phi} \, \frac{\Phi\gamma^2 + (2-\beta\Phi)\gamma - \beta}{\gamma^2(\beta-\gamma)^2},
\end{align*}
and $\Pi^\Phi(\gamma,\beta)$ viewed as a function of $\gamma$ attains its minimum at its critical point given by $\overline{\gamma}^\Phi(\beta),$
since $\frac{\ud \Pi^\Phi}{\ud \gamma}\big(\gamma,\beta\big)<0$ on $(0,\overline{\gamma}^\Phi(\beta))$ and $\frac{\ud \Pi^\Phi}{\ud \gamma}\big(\gamma,\beta\big)>0$ on $(\overline{\gamma}^\Phi(\beta),\beta)$.

\vspace{0.5em}
Now, we proceed to the second case, and start by determining the critical points of $\Pi^\Phi_\star$.
It holds 
\begin{align*}
\frac{\partial}{\partial\gamma}\Pi_\star^\Phi(\gamma,\delta)
	& = -\frac{8}{\gamma^2}  
		+ 9\delta\e{(\delta-\gamma)\Phi} \frac{\Phi \gamma^2  - (\delta \Phi -2)\gamma - \delta}{\gamma^2(\delta-\gamma)^2},\\
\frac{\partial }{\partial\delta}\Pi_\star^\Phi(\gamma,\delta)
	& = -\frac{9}{\delta^2} + 9 \e{(\delta-\gamma)\Phi} \frac{(1+\delta\Phi)(\delta-\gamma)-\delta}{\gamma(\delta-\gamma)^2}.
\end{align*}
Following analogous computations as above we can prove that, for $(\gamma,\delta)\in\mathcal C_\beta$, the equation
\begin{align*}
\frac{\partial }{\partial \gamma}\Pi^\Phi_\star(\gamma,\delta)=0 
	&\Leftrightarrow 
	P_\delta(\gamma):= 8(\delta-\gamma)^2 - 9\delta\e{(\delta-\gamma)\Phi}\big(\Phi \gamma^2 - (\delta\Phi-2)\gamma - \delta\big)=0
\end{align*}  
has a unique root, say $\overline{\gamma}_\star^\Phi(\delta),$ which moreover satisfies $\overline{\gamma}_\star^\Phi(\delta)\in(\overline{\gamma}^\Phi(\delta),\delta).$
This can be proved because the function $P_\delta:(0,\delta)\to\mathbb R$ is decreasing, for each fixed $\delta\in(0,\beta)$, with $P_\delta\big(\overline{\gamma}^\Phi(\delta)\big)>0 $ and $P_\delta(\delta)<0$.
Now observe that for $\gamma>\frac{\delta^2\Phi}{1+\delta\Phi}$ it holds $\frac{\partial }{\partial\delta}\Pi_\star^\Phi(\gamma,\delta)<0 $ and that $P_\delta\big(\frac{\delta^2\Phi}{1+\delta\Phi}\big)>0.$
Using the monotonicity of $P_\delta$ we have that $\overline{\gamma}_\star^\Phi(\delta)>\frac{\delta^2\Phi}{1+\delta\Phi}$ and therefore also  $\frac{\partial }{\partial\delta}\Pi_\star^\Phi(\overline{\gamma}_\star^\Phi(\delta),\delta)<0.$
Arguing as above we can conclude that the infimum is attained for $\delta=\beta$ at the point $\big(\overline{\gamma}_\star^\Phi(\beta),\beta\big).$ 

\vspace{.5em}
Finally, the limiting statements follow by straightforward but tedious computations.
\end{proof}

%% file: Auxilliary_Results_Optional_Measures.tex
Let $(\Omega, \mathcal G, \mathbb G, \mathbb P)$ be a filtered probability space and $Y=\{Y_t\}_{t\in[0,\infty]}$  be a uniformly integrable measurable process. Then, thanks to the uniform integrability, we have by a clear adaptation of \cite[Theorem 5.1]{he1992semimartingale} that there exists a unique optional process, denoted by $^oY$, such that for every $\mathbb G-$stopping time $\tau$, we have $\mathbb E_{\tau}[Y_{\tau}]=\! ^oY_{\tau},$ $\Pm-a.s.$ Observe that $\tau$ is allowed to take infinite values, since $Y_{\infty}$ is well-defined and integrable.
For any increasing, c\`adl\`ag and $\mathbb G-$adapted process $A$, the measure $\mu_A:(\Omega\times[0,\infty], \G\otimes\mathcal B([0,\infty]) )\longrightarrow (\mathbb R,\mathcal B(\mathbb R) )$ defined as
\begin{align*}
\mu_A(H)=\mathbb E\bigg[\int_0^\infty \mathds{1}_{H} \dA t\bigg], \textrm{ for } H\in\G\otimes\mathcal B([0,\infty]),
\end{align*}
is optional, see \cite[Definition 5.10, Definition 5.12 and Theorem 5.13]{he1992semimartingale}.
For convenience, we state the following well-known result as a lemma.

\begin{lemma}\label{OptFubini}
Let $A$ be an increasing, \cadlag~ and adapted process and $Y$ be a uniformly integrable and measurable process, then it holds
$\mu_A( Y ) =\mu_A( ^oY ).$
\end{lemma}
\begin{proof}
Let $L$ be the c\`agl\`ad generalized inverse of $A$ (see Lemma \ref{mainlemma} for the definition). We have
\begin{align*}
\mathbb E\bigg[\int_0^{\infty}{Y_t} \dA t\bigg] 
&=\mathbb E\bigg[ \int_0^{\infty}Y_{L_s} \mathds{1}_{[L_s<\infty]} \ds \bigg]
=\int_0^{\infty}\mathbb E[ Y_{L_s} \mathds{1}_{[L_s<\infty]}] \ds\\
&=\int_0^{\infty}\mathbb E[ ^oY_{L_s} \mathds{1}_{[L_s<\infty]}] \ds
=\mathbb E\bigg[\int_0^{\infty}{^oY_t} \dA t\bigg],
\end{align*}
where for the change of variables we used \cite[Lemma 1.38]{he1992semimartingale}, and for the third equality the definition of the optional projection, see \cite[Theorem 5.1]{he1992semimartingale}.
\end{proof}

%% file: Auxilliary_Results_New_estimates.tex
\vspace{0.5em}
\textbullet\ \textbf{Case A. :} $X^{\circ,c}\ne 0$.

\vspace{0.5em}
In order to obtain a contraction, we need to choose positive predictable processes $\varepsilon$ and $\gamma$ such that
$$\begin{cases}
\displaystyle \varepsilon_s>\Delta C_s,\\[0.5em]
\displaystyle 1-(1+\gamma_s\Delta C_s)(\varepsilon_s-\Delta C_s)\theta^\circ_s> 0,\\[0.5em]
\displaystyle 1-\theta^\natural_s(\varepsilon_s-\Delta C_s)> 0,\\[0.5em]
\displaystyle 1+C_s\big(\gamma_s-(1+\gamma_s\Delta C_s)\varepsilon_s^{-1}\big)>0,\\[0.5em]
\displaystyle \frac{C_sr_s(1+\gamma_s\Delta C_s)(\varepsilon_s-\Delta C_s)}{1+C_s\big(\gamma_s-(1+\gamma_s\Delta C_s)\varepsilon_s^{-1}\big)}<1,
\end{cases}\Longleftrightarrow \begin{cases}
\displaystyle\Delta C_s<\varepsilon_s<\Delta C_s+\frac{1}{\theta^\circ_s\vee\theta^{\natural}_s} ,\\
\displaystyle \frac{\big(C_s-\varepsilon_s\big)^+}{C_s(\varepsilon_s-\Delta C_s)}< \gamma_s<\frac{1-\theta^\circ_s(\varepsilon_s-\Delta C_s)}{\theta^\circ_s\Delta C_s(\varepsilon_s-\Delta C_s)},\\[0.5em]
\displaystyle (\varepsilon_sr_s\Delta C_s-1)\gamma_s< \frac{-r_sC_s\varepsilon_s^2+(1+r_sC_s\Delta C_s)\varepsilon_s-C_s}{C_s(\varepsilon_s-\Delta C_s)}.
\end{cases}$$
We need then to distinguish two cases. 

\vspace{0.5em}
\hspace{2em}\textbullet\ {\bf Case A.i. :} $\varepsilon_sr_s\Delta C_s>1$.
Then, it can be proven that the system is not compatible.

\vspace{0.5em}
\hspace{2em}\textbullet\ {\bf Case A.ii. :} $\varepsilon_sr_s\Delta C_s<1$.
We need to consider two sub--cases

\vspace{0.5em}
\hspace{4em} \textbullet\ {\bf Case A.ii.a. :} $-r_sC_s\varepsilon_s^2+(1+r_sC_s\Delta C_s)\varepsilon_s-C_s\geq 0.$

\vspace{0.5em}
The above condition, after some computations, implies that
$$\begin{cases}
r_s\in\bigg(0,\frac{(\sqrt{C_s}-\sqrt{C_{s-}})^2}{C_s(\Delta C_s)^2}\bigg)\bigcup\bigg(\frac{(\sqrt{C_s}+\sqrt{C_{s-}})^2}{C_s(\Delta C_s)^2},+\infty\bigg),\\[0.8em]
\displaystyle\frac{1+r_sC_s\Delta C_s-\sqrt{(1+r_sC_s\Delta C_s)^2-4C_s^2r_s}}{2r_sC_s}<\varepsilon_s<\frac{1+r_sC_s\Delta C_s+\sqrt{(1+r_sC_s\Delta C_s)^2-4C_s^2r_s}}{2r_sC_s}.
\end{cases}$$
This will be compatible with $\varepsilon_sr_s\Delta C_s<1$ if and only if
\begin{align*}
&\frac{1+r_sC_s\Delta C_s-\sqrt{(1+r_sC_s\Delta C_s)^2-4C_s^2r_s}}{2r_sC_s}<\frac{1}{r_s\Delta C_s}
\Longleftrightarrow r_s< \frac{C_s+C_{s-}}{C_s\Delta C_s^2}.
\end{align*}

Therefore, the system now becomes
$$\begin{cases}
\displaystyle\max\bigg\{\Delta C_s,\frac{1+r_sC_s\Delta C_s-\sqrt{(1+r_sC_s\Delta C_s)^2-4C_s^2r_s}}{2r_sC_s}\bigg\}<\varepsilon_s,\\[1em]
\displaystyle\ \varepsilon_s<\min\bigg\{\frac{1}{r_s\Delta C_s},\Delta C_s+\frac{1}{\theta^\circ_s\vee\theta^{\natural}_s},\frac{1+r_sC_s\Delta C_s+\sqrt{(1+r_sC_s\Delta C_s)^2-4C_s^2r_s}}{2r_sC_s}\bigg\},\\[1em]
\displaystyle \frac{\big(C_s-\varepsilon_s\big)^+}{C_s(\varepsilon_s-\Delta C_s)}\leq \gamma_s<\frac{1-\theta^\circ_s(\varepsilon_s-\Delta C_s)}{\theta^\circ_s\Delta C_s(\varepsilon_s-\Delta C_s)},
\end{cases}$$
with the requirement that 
$$r_s< \frac{(\sqrt{C_s}-\sqrt{C_{s-}})^2}{C_s(\Delta C_s)^2}.$$
For the system to have solutions, we necessarily need to have
\begin{align*}
&\frac{1+r_sC_s\Delta C_s+\sqrt{(1+r_sC_s\Delta C_s)^2-4C_s^2r_s}}{2r_sC_s}>\Delta C_s\Longleftrightarrow r_s< \frac{1}{C_s\Delta C_s},
\end{align*}
as well as
$$\frac{1}{r_s\Delta C_s}>\Delta C_s\Longleftrightarrow r_s<\frac1{\Delta C_s^2},$$
and in addition
\begin{align*}
&\Delta C_s+\frac{1}{\theta^\circ_s\vee\theta^{\natural}_s}>\frac{1+r_sC_s\Delta C_s-\sqrt{(1+r_sC_s\Delta C_s)^2-4C_s^2r_s}}{2r_sC_s}\\
\Longleftrightarrow&\; r_s> \frac{\theta^\circ_s\vee\theta^{\natural}_s}{C_s(2+(\theta^\circ_s\vee\theta^{\natural}_s)\Delta C_s)}\;
\text{or}\; r_s<\frac{\theta^\circ_s\vee\theta^{\natural}_s(1-C_{s-}(\theta^\circ_s\vee\theta^{\natural}_s))}{C_s(1+(\theta^\circ_s\vee\theta^{\natural}_s)\Delta C_s)}.
\end{align*}
The last four conditions on $r$ are then equivalent to
$$r_s<\min\bigg\{\frac{(\sqrt{C_s}-\sqrt{C_{s-}})^2}{C_s(\Delta C_s)^2},\frac{\theta^\circ_s\vee\theta^{\natural}_s(1-C_{s-}(\theta^\circ_s\vee\theta^{\natural}_s))}{C_s(1+(\theta^\circ_s\vee\theta^{\natural}_s)\Delta C_s)}\bigg\},\ \text{or}\ \frac{\theta^\circ_s\vee\theta^{\natural}_s}{C_s(2+(\theta^\circ_s\vee\theta^{\natural}_s)\Delta C_s)}<r_s<\frac{(\sqrt{C_s}-\sqrt{C_{s-}})^2}{C_s(\Delta C_s)^2}.$$
We then need to distinguish two further sub--cases

\vspace{0.5em}
\hspace{6em} \textbullet\ {\bf Case A.ii.a.1. :} $\varepsilon_s>C_s$

\vspace{0.5em}
Under this additional condition, we necessarily need to have
\begin{align*}
&\frac{1+r_sC_s\Delta C_s+\sqrt{(1+r_sC_s\Delta C_s)^2-4C_s^2r_s}}{2r_sC_s}> C_s\Longleftrightarrow r_s<\frac{1}{C_s(C_s+C_{s-})},
\end{align*}
as well as
$$\Delta C_s+\frac1{\theta^\circ_s\vee\theta^{\natural}_s}>C_s\Longleftrightarrow(\theta^\circ_s\vee\theta^{\natural}_s)C_{s-}<1,$$
and
$$\frac{1}{r_s\Delta C_s}>C_s\Longleftrightarrow r_s<\frac1{C_s\Delta C_s},$$
so that in this case the final system is
$$\begin{cases}
\displaystyle\max\bigg\{ C_s,\frac{1+r_sC_s\Delta C_s-\sqrt{(1+r_sC_s\Delta C_s)^2-4C_s^2r_s}}{2r_sC_s}\bigg\}<\varepsilon_s,\\[1em]
\displaystyle\ \varepsilon_s<\min\bigg\{\frac{1}{r_s\Delta C_s},\Delta C_s+\frac{1}{\theta^\circ_s\vee\theta^{\natural}_s},\frac{1+r_sC_s\Delta C_s+\sqrt{(1+r_sC_s\Delta C_s)^2-4C_s^2r_s}}{2r_sC_s}\bigg\},\\[1em]
\displaystyle 0\leq \gamma_s<\frac{1-\theta^\circ_s(\varepsilon_s-\Delta C_s)}{\theta^\circ_s\Delta C_s(\varepsilon_s-\Delta C_s)}, 
\end{cases}$$
which has solutions if and only if
$$\begin{cases}
\displaystyle C_{s-}(\theta^\circ_s\vee\theta^{\natural}_s)<1,\\[1em]
\displaystyle r_s\in\bigg(\frac{\theta^\circ_s\vee\theta^{\natural}_s}{C_s(2+(\theta^\circ_s\vee\theta^{\natural}_s)\Delta C_s)},\frac{(\sqrt{C_s}-\sqrt{C_{s-}})^2}{C_s(\Delta C_s)^2}\bigg),
\end{cases}$$
or
$$
\begin{cases}
\displaystyle C_{s-}(\theta^\circ_s\vee\theta^{\natural}_s)<1,\\[1em]
\displaystyle r_s<\min\bigg\{\frac{(\theta^\circ_s\vee\theta^{\natural}_s)(1-C_{s-}(\theta^\circ_s\vee\theta^{\natural}_s))}{C_s(1+(\theta^\circ_s\vee\theta^{\natural}_s)\Delta C_s)},\frac{(\sqrt{C_s}-\sqrt{C_{s-}})^2}{C_s(\Delta C_s)^2}\bigg\}.
\end{cases}$$
Therefore, we will discard the condition that imposes a lower bound on $r$ and we will keep only the right one.

\vspace{0.5em}
\hspace{6em} \textbullet\ {\bf Case A.ii.a.2 :} $\varepsilon_s\leq C_s$

\vspace{0.5em}

Under this additional condition, we necessarily need to have
\begin{align*}
&\frac{1+r_sC_s\Delta C_s-\sqrt{(1+r_sC_s\Delta C_s)^2-4C_s^2r_s}}{2r_sC_s}< C_s\Longleftrightarrow r_s>\frac{1}{C_s(C_s+C_{s-})}.
\end{align*}
However, this is not compatible with the constraint $r_s<(\sqrt{C_s}-\sqrt{C_{s-}})^2/(C_s(\Delta C_s)^2)$, and the system does not admit any solutions in this case.

\vspace{0.5em}
\hspace{4em} \textbullet\ {\bf Case A.ii.b. :} $-r_sC_s\varepsilon_s^2+(1+r_sC_s\Delta C_s)\varepsilon_s-C_s< 0.$

\vspace{0.5em}
This requires either that 
$$r_s\in\bigg(\frac{(\sqrt{C_s}-\sqrt{C_{s-}})^2}{C_s(\Delta C_s)^2},\frac{(\sqrt{C_s}+\sqrt{C_{s-}})^2}{C_s(\Delta C_s)^2}\bigg),$$
and no further restrictions on $\varepsilon_s$, or 
$$r_s\in\bigg(0,\frac{(\sqrt{C_s}-\sqrt{C_{s-}})^2}{C_s(\Delta C_s)^2}\bigg)\bigcup\bigg(\frac{(\sqrt{C_s}+\sqrt{C_{s-}})^2}{C_s(\Delta C_s)^2},+\infty\bigg),$$
and either
$$\varepsilon_s<\frac{1+r_sC_s\Delta C_s-\sqrt{(1+r_sC_s\Delta C_s)^2-4C_s^2r_s}}{2r_sC_s},\; \text{or}\; \varepsilon_s>\frac{1+r_sC_s\Delta C_s+\sqrt{(1+r_sC_s\Delta C_s)^2-4C_s^2r_s}}{2r_sC_s}.$$
Let us distinguish between these two cases.

\vspace{0.5em}
\hspace{6em} \textbullet\ {\bf Case A.ii.b.1. :} $r_s\in\bigg(\frac{(\sqrt{C_s}-\sqrt{C_{s-}})^2}{C_s(\Delta C_s)^2},\frac{(\sqrt{C_s}+\sqrt{C_{s-}})^2}{C_s(\Delta C_s)^2}\bigg).$
We will not examine this case since $r$ will be bounded from below.

\vspace{0.5em}
\hspace{8em} \textbullet\ {\bf Case A.ii.b.1.$\boldsymbol{\alpha}$. :} $\varepsilon_s<\min\Big\{C_s,\frac{C_s}{\theta^\circ_sC_{s-}}\Big\}.$

\vspace{0.5em}
Then the system becomes
$$\begin{cases}
\displaystyle\Delta C_s<\varepsilon_s<\min\bigg\{\frac{1}{r_s\Delta C_s},C_s,\frac{C_s}{\theta^\circ_sC_{s-}+r_sC_s\Delta C_s},\Delta C_s+\frac{1}{\theta^\circ_s\vee\theta^{\natural}_s}\bigg\} ,\\[0.8em]
\displaystyle\max\bigg\{\frac{r_sC_s\varepsilon_s^2-(1+r_sC_s\Delta C_s)\varepsilon_s+C_s}{C_s(\varepsilon_s-\Delta C_s)(1-r_s\Delta C_s\varepsilon_s)},\frac{C_s-\varepsilon_s}{C_s(\varepsilon_s-\Delta C_s)}\bigg\}< \gamma_s<\frac{1-\theta^\circ_s(\varepsilon_s-\Delta C_s)}{\theta^\circ_s\Delta C_s(\varepsilon_s-\Delta C_s)},
\end{cases}$$
which admits solutions if and only if
{$$\begin{cases}
\displaystyle \theta^\circ_s\Delta C_s<{2}\sqrt{\frac{C_s}{C_{s-}}}-1,\\[1em]
\displaystyle \frac{(\sqrt{C_s}-\sqrt{C_{s-}})^2}{C_s(\Delta C_s)^2}<r_s<\min\bigg\{\frac1{\Delta C_s^2},\frac{\Delta C_s+C_{s-}(1-\theta^\circ_s\Delta C_s)}{C_s(\Delta C_s)^2}\bigg\}.
\end{cases}$$}
Therefore, we will not take into account this case.

\vspace{0.5em}
\hspace{8em} \textbullet\ {\bf Case A.ii.b.1.$\boldsymbol{\beta}$. :} $\varepsilon_s\geq C_s.$

\vspace{0.5em}
Then the system becomes
$$\begin{cases}
\displaystyle C_s\leq \varepsilon_s<\min\bigg\{\frac{1}{r_s\Delta C_s},\frac{C_s}{\theta^\circ_sC_{s-}+r_sC_s\Delta C_s},\Delta C_s+\frac{1}{\theta^\circ_s\vee\theta^{\natural}_s}\bigg\} ,\\[0.8em]
\displaystyle \frac{r_sC_s\varepsilon_s^2-(1+r_sC_s\Delta C_s)\varepsilon_s+C_s}{C_s(\varepsilon_s-\Delta C_s)(1-r_s\Delta C_s\varepsilon_s)}< \gamma_s<\frac{1-\theta^\circ_s(\varepsilon_s-\Delta C_s)}{\theta^\circ_s\Delta C_s(\varepsilon_s-\Delta C_s)},
\end{cases}$$
which admits solutions if and only if
{$$\begin{cases}
\displaystyle (\theta^\circ_s\vee\theta^{\natural}_s) C_{s-}<1,\\
\displaystyle \theta^\circ_s<\min\bigg\{\frac{2\sqrt{C_s}-\sqrt{C_{s-}}}{\Delta C_s\sqrt{C_{s-}}},\frac1{C_{s-}},\frac{2(\sqrt{C_s}-\sqrt{C_{s-}})}{\Delta C_s\sqrt{C_{s-}}}\bigg\}=\frac{2(\sqrt{C_s}-\sqrt{C_{s-}})}{\Delta C_s\sqrt{C_{s-}}},\\[1em]
\displaystyle \frac{(\sqrt{C_s}-\sqrt{C_{s-}})^2}{C_s(\Delta C_s)^2}<r_s<\min\bigg\{\frac1{C_s\Delta C_s},\frac{\Delta C_s+C_{s-}(1-\theta^\circ_s\Delta C_s)}{C_s(\Delta C_s)^2},\frac{1-\theta^\circ_sC_{s-}}{C_s\Delta C_s}\bigg\}=\frac{1-\theta^\circ_sC_{s-}}{C_s\Delta C_s}.
\end{cases}$$}
Therefore, we will not take into account this case as well.

\vspace{0.5em}
\hspace{6em}\textbullet\ {\bf Case A.ii.b.2. :} $r_s\in\bigg(0,\frac{(\sqrt{C_s}-\sqrt{C_{s-}})^2}{C_s(\Delta C_s)^2}\bigg)\bigcup\bigg(\frac{(\sqrt{C_s}+\sqrt{C_{s-}})^2}{C_s(\Delta C_s)^2},+\infty\bigg).$

\vspace{0.5em}
The system becomes either
$$\begin{cases}
\displaystyle\Delta C_s<\varepsilon_s<\min\bigg\{\frac{1}{r_s\Delta C_s},\frac{1+r_sC_s\Delta C_s-\sqrt{(1+r_sC_s\Delta C_s)^2-4C_s^2r_s}}{2r_sC_s},\Delta C_s+\frac{1}{\theta^\circ_s\vee\theta^{\natural}_s}\bigg\} ,\\[0.8em]
\displaystyle\max\bigg\{\frac{r_sC_s\varepsilon_s^2-(1+r_sC_s\Delta C_s)\varepsilon_s+C_s}{C_s(\varepsilon_s-\Delta C_s)(1-r_s\Delta C_s\varepsilon_s)},\frac{\big(C_s-\varepsilon_s\big)^+}{C_s(\varepsilon_s-\Delta C_s)}\bigg\}< \gamma_s<\frac{1-\theta^\circ_s(\varepsilon_s-\Delta C_s)}{\theta^\circ_s\Delta C_s(\varepsilon_s-\Delta C_s)},
\end{cases}$$
or
$$\begin{cases}
\displaystyle\max\bigg\{\Delta C_s,\frac{1+r_sC_s\Delta C_s+\sqrt{(1+r_sC_s\Delta C_s)^2-4C_s^2r_s}}{2r_sC_s}\bigg\}<\varepsilon_s<\min\bigg\{\frac{1}{r_s\Delta C_s},\Delta C_s+\frac{1}{\theta^\circ_s\vee\theta^{\natural}_s}\bigg\} ,\\[0.8em]
\displaystyle\max\bigg\{\frac{r_sC_s\varepsilon_s^2-(1+r_sC_s\Delta C_s)\varepsilon_s+C_s}{C_s(\varepsilon_s-\Delta C_s)(1-r_s\Delta C_s\varepsilon_s)},\frac{\big(C_s-\varepsilon_s\big)^+}{C_s(\varepsilon_s-\Delta C_s)}\bigg\}< \gamma_s<\frac{1-\theta^\circ_s(\varepsilon_s-\Delta C_s)}{\theta^\circ_s\Delta C_s(\varepsilon_s-\Delta C_s)}.\end{cases}$$
In both cases, this imposes that 
$$r_s\Delta C_s^2<1,$$
as well as
$$\frac{(C_s-\varepsilon_s)^+}{C_s(\varepsilon_s-\Delta C_s)}<\frac{1-\theta^\circ_s(\varepsilon_s-\Delta C_s)}{\theta^\circ_s\Delta C_s(\varepsilon_s-\Delta C_s)}\Longleftrightarrow \varepsilon_s<\min\bigg\{C_s,\frac{C_s}{\theta^\circ_sC_{s-}}\bigg\}, \text{ or }\varepsilon_s\geq C_s,$$
and
\begin{align*}
&\frac{1-\theta^\circ_s(\varepsilon_s-\Delta C_s)}{\theta^\circ_s\Delta C_s(\varepsilon_s-\Delta C_s)}>\frac{r_sC_s\varepsilon_s^2-(1+r_sC_s\Delta C_s)\varepsilon_s+C_s}{C_s(\varepsilon_s-\Delta C_s)(1-r_s\Delta C_s\varepsilon_s)}\Longleftrightarrow\varepsilon_s<\frac{C_s}{\theta^\circ_sC_{s-}+r_sC_s\Delta C_s}.
\end{align*}
This requires in turn that we must necessarily have
\begin{align*}
&\frac{C_s}{\theta^\circ_sC_{s-}+r_sC_s\Delta C_s}>\Delta C_s
\Longleftrightarrow r_s<\frac{\Delta C_s+C_{s-}(1-\theta^\circ_s\Delta C_s)}{C_s(\Delta C_s)^2}\;\text{and}\; \theta^\circ_s<\frac{C_s}{C_{s-}\Delta C_s}.
\end{align*}
Furthermore, the first system requires in addition that
$$\frac{1+r_sC_s\Delta C_s-\sqrt{(1+r_sC_s\Delta C_s)^2-4C_s^2r_s}}{2r_sC_s}>\Delta C_s\Longleftrightarrow r_s<\frac{1}{C_s\Delta C_s},$$
while the second one requires
$$\frac{1+r_sC_s\Delta C_s+\sqrt{(1+r_sC_s\Delta C_s)^2-4C_s^2r_s}}{2r_sC_s}<\frac{1}{r_s\Delta C_s}\Longleftrightarrow r_s<\frac{C_s+C_{s-}}{C_s\Delta C_s^2},$$
as well as
\begin{align*}
&\Delta C_s+\frac{1}{\theta^\circ_s\vee\theta^{\natural}_s}>\frac{1+r_sC_s\Delta C_s+\sqrt{(1+r_sC_s\Delta C_s)^2-4C_s^2r_s}}{2r_sC_s}\\
\Longleftrightarrow&\; r_s> \max\bigg\{\frac{\theta^\circ_s\vee\theta^{\natural}_s}{C_s(2+(\theta^\circ_s\vee\theta^{\natural}_s)\Delta C_s)},\frac{(\theta^\circ_s\vee\theta^{\natural}_s)(1-C_{s-}(\theta^\circ_s\vee\theta^{\natural}_s))}{C_s(1+(\theta^\circ_s\vee\theta^{\natural}_s)\Delta C_s)}\bigg\},
\end{align*}
and
$$\frac{1+r_sC_s\Delta C_s+\sqrt{(1+r_sC_s\Delta C_s)^2-4C_s^2r_s}}{2r_sC_s}<\frac{C_s}{\theta^\circ_sC_{s-}+r_sC_s\Delta C_s}.$$
We do not proceed, since the process $r$ will be bounded from below.

\vspace{0.5em}
\hspace{8em} \textbullet\ {\bf Case A.ii.b.2.$\boldsymbol{\alpha}$. :} $\varepsilon_s<\min\Big\{C_s,\frac{C_s}{\theta^\circ_sC_{s-}}\Big\}.$

\vspace{0.5em}
Then the system becomes
$$\begin{cases}
\displaystyle\Delta C_s<\varepsilon_s<\min\bigg\{\frac{1}{r_s\Delta C_s},C_s,\frac{C_s}{\theta^\circ_sC_{s-}+r_sC_s\Delta C_s},\frac{1+r_sC_s\Delta C_s-\sqrt{(1+r_sC_s\Delta C_s)^2-4C_s^2r_s}}{2r_sC_s},\Delta C_s+\frac{1}{\theta^\circ_s\vee\theta^{\natural}_s}\bigg\} ,\\[0.8em]
\displaystyle\max\bigg\{\frac{r_sC_s\varepsilon_s^2-(1+r_sC_s\Delta C_s)\varepsilon_s+C_s}{C_s(\varepsilon_s-\Delta C_s)(1-r_s\Delta C_s\varepsilon_s)},\frac{C_s-\varepsilon_s}{C_s(\varepsilon_s-\Delta C_s)}\bigg\}< \gamma_s<\frac{1-\theta^\circ_s(\varepsilon_s-\Delta C_s)}{\theta^\circ_s\Delta C_s(\varepsilon_s-\Delta C_s)},
\end{cases}$$
which admits solutions if and only if
$$\begin{cases}
\displaystyle \Delta C_s\theta^\circ_s<\frac{C_s}{C_{s-}},\\
\displaystyle(\Delta C_s)^2 r_s<\min\bigg\{\frac{\Delta C_s+C_{s-}(1-\theta^\circ_s\Delta C_s)}{C_s},\frac{(\sqrt{C_s}-\sqrt{C_{s-}})^2}{C_s}\bigg\}.
\end{cases}$$

\vspace{0.5em}
\hspace{8em} \textbullet\ {\bf Case A.ii.b.2.$\boldsymbol{\beta}$. :} $\varepsilon_s\geq C_s.$

\vspace{0.5em}
Then the system becomes
$$\begin{cases}
\displaystyle C_s<\varepsilon_s<\min\bigg\{\frac{1}{r_s\Delta C_s},\frac{C_s}{\theta^\circ_sC_{s-}+r_sC_s\Delta C_s},\frac{1+r_sC_s\Delta C_s-\sqrt{(1+r_sC_s\Delta C_s)^2-4C_s^2r_s}}{2r_sC_s},\Delta C_s+\frac{1}{\theta^\circ_s\vee\theta^{\natural}_s}\bigg\} ,\\[0.8em]
\displaystyle\frac{r_sC_s\varepsilon_s^2-(1+r_sC_s\Delta C_s)\varepsilon_s+C_s}{C_s(\varepsilon_s-\Delta C_s)(1-r_s\Delta C_s\varepsilon_s)}< \gamma_s<\frac{1-\theta^\circ_s(\varepsilon_s-\Delta C_s)}{\theta^\circ_s\Delta C_s(\varepsilon_s-\Delta C_s)},
\end{cases}$$
which admits solutions if and only if
{$$\begin{cases}
\displaystyle (\theta^\circ_s\vee\theta^{\natural}_s) C_{s-}<1,\\[1em]
\displaystyle r_s<\min\bigg\{ \frac{(\sqrt{C_s}-\sqrt{C_{s-}})^2}{C_s(\Delta C_s)^2},\frac{1}{C_s(C_s+C_{s-})},\frac{1-\theta^\circ_sC_{s-}}{C_s\Delta C_s}\bigg\}=\min\bigg\{ \frac{(\sqrt{C_s}-\sqrt{C_{s-}})^2}{C_s(\Delta C_s)^2},\frac{1-\theta^\circ_sC_{s-}}{C_s\Delta C_s}\bigg\}.
\end{cases}$$}

\vspace{0.5em}
\textbullet\ \textbf{Case B.} $X^{\circ,c} = 0$.

\vspace{0.5em}
In order to obtain a contraction, we need to choose positive predictable processes $\varepsilon$ and $\gamma$ such that
$$\begin{cases}
\displaystyle \varepsilon_s>\Delta C_s,\\[0.5em]
\displaystyle 1-\theta^\circ_s(\varepsilon_s-\Delta C_s)> 0,\\[0.5em]
\displaystyle 1-\theta^\natural_s(\varepsilon_s-\Delta C_s)> 0,\\[0.5em]
\displaystyle 1+C_s\big(\gamma_s-(1+\gamma_s\Delta C_s)\varepsilon_s^{-1}\big)>0,\\[0.5em]
\displaystyle \frac{C_sr_s(1+\gamma_s\Delta C_s)(\varepsilon_s-\Delta C_s)}{1+C_s\big(\gamma_s-(1+\gamma_s\Delta C_s)\varepsilon_s^{-1}\big)}<1,
\end{cases}\Longleftrightarrow \begin{cases}
\displaystyle\Delta C_s<\varepsilon_s<\Delta C_s+\frac{1}{\theta^\circ_s\vee\theta^{\natural}_s} ,\\[0.5em]
\displaystyle  \gamma_s<\frac{1-\theta^\circ_s(\varepsilon_s-\Delta C_s)}{\theta^\circ_s\Delta C_s(\varepsilon_s-\Delta C_s)},\\[0.5em]
\displaystyle (\varepsilon_sr_s\Delta C_s-1)\gamma_s< \frac{-r_sC_s\varepsilon_s^2+(1+r_sC_s\Delta C_s)\varepsilon_s-C_s}{C_s(\varepsilon_s-\Delta C_s)}.
\end{cases}$$
This is exactly the same system as in {\bf Case A}, except that we no longer need the inequality $(C_s-\varepsilon_s)^+/(C_s(\varepsilon_s-\Delta C_s))<\gamma_s$. Hence, the exact same reasoning as before will tell us that the system admits solutions if one of the following set of conditions is satisfied

$$r_s<\min\bigg\{\frac{(\sqrt{C_s}-\sqrt{C_{s-}})^2}{C_s(\Delta C_s)^2},\frac{\theta^\circ_s\vee\theta^{\natural}_s(1-C_{s-}(\theta^\circ_s\vee\theta^{\natural}_s))}{C_s(1+(\theta^\circ_s\vee\theta^{\natural}_s)\Delta C_s)}\bigg\},\ \text{or}\ \frac{\theta^\circ_s\vee\theta^{\natural}_s}{C_s(2+(\theta^\circ_s\vee\theta^{\natural}_s)\Delta C_s)}<r_s<\frac{(\sqrt{C_s}-\sqrt{C_{s-}})^2}{C_s(\Delta C_s)^2},$$
or
{$$\begin{cases}
\displaystyle \theta^\circ_s\Delta C_s<{2}\sqrt{\frac{C_s}{C_{s-}}}-1,\\[1em]
\displaystyle \frac{(\sqrt{C_s}-\sqrt{C_{s-}})^2}{C_s(\Delta C_s)^2}<r_s<\min\bigg\{\frac1{\Delta C_s^2},\frac{\Delta C_s+C_{s-}(1-\theta^\circ_s\Delta C_s)}{C_s(\Delta C_s)^2}\bigg\}
\end{cases},$$}
or
$$\begin{cases}
\displaystyle \Delta C_s\theta^\circ_s<\frac{C_s}{C_{s-}},\\
\displaystyle(\Delta C_s)^2 r_s<\min\bigg\{\frac{\Delta C_s+C_{s-}(1-\theta^\circ_s\Delta C_s)}{C_s},\frac{(\sqrt{C_s}-\sqrt{C_{s-}})^2}{C_s}\bigg\}
\end{cases}.$$